\documentclass[11pt,oneside]{article}
\usepackage{amssymb}
\usepackage{amsmath}
\usepackage{latexsym}
\usepackage{amsfonts}
\usepackage{amscd}
\usepackage{amsthm}
\usepackage{fancyhdr}

\usepackage[headings]{fullpage}

\newcommand{\gt}[1]{\mathfrak{#1}}

\newcommand{\mc}[1]{\mathcal{#1}}

\newcommand{\RR}{{\mathbb R}}%Reals
\newcommand{\CC}{{\mathbb C}}%Complex
\newcommand{\ZZ}{{\mathbb Z}}%Integers
\newcommand{\ZZh}{\ZZ+\tfrac{1}{2}}%Integers plus half
\newcommand{\hZZ}{\tfrac{1}{2}\ZZ}%half integers

\newcommand{\QQ}{{\mathbb Q}}%Rationals
\newcommand{\hh}{{\mathbf h}}%hyperbolic space
\newcommand{\FF}{{\mathbb F}}%Field

\newcommand{\fset}{{\Sigma}}%A finite set
\newcommand{\hfset}{{\Delta}}    %Half a finite set

\newcommand{\ii}{{\bf i}}%Sqrt of -1
\newcommand{\lab}{{\langle}}    %Left angle brackets
\newcommand{\rab}{{\rangle}}    %Right angle brackets
\newcommand{\tcolon}{\!:\!} %Tight colon

\newcommand{\LL}{{\Lambda}}

\newcommand{\End}{\operatorname{End}}
\newcommand{\Aut}{\operatorname{Aut}}
\newcommand{\Id}{\operatorname{Id}}
\newcommand{\tr}{\operatorname{{\sf tr}|}}
\newcommand{\Fix}{\operatorname{Fix}}
\newcommand{\Res}{\operatorname{Res}}
\newcommand{\Sing}{\operatorname{Sing}}
\newcommand{\Lie}{\operatorname{Lie}}
\newcommand{\Img}{\operatorname{Im}}
\newcommand{\Span}{\operatorname{Span}}
\newcommand{\rank}{\operatorname{rank}}
\newcommand{\Gal}{\operatorname{Gal}}
\newcommand{\Inn}{\operatorname{Inn}}

\newcommand{\Kumgp}{\lab\!\pm\!\Id\!\rab}
\newcommand{\Kumz}{\lab\ogz\rab}

\newcommand{\PSL}{\operatorname{\textsl{PSL}}}    %PSL group
\newcommand{\SL}{\operatorname{\textsl{SL}}}
\newcommand{\GL}{\operatorname{\textsl{GL}}}

\newcommand{\Sp}{\operatorname{\textsl{Spin}}}    %Spin group
\newcommand{\Pin}{\operatorname{\textsl{Pin}}}
\newcommand{\SO}{\operatorname{\textsl{SO}}}    %SO group
\newcommand{\Or}{\operatorname{\textsl{O}}}       %Orthog group
\newcommand{\U}{\operatorname{\textsl{U}}}        %Unitary group
\newcommand{\SU}{\operatorname{\textsl{SU}}}    %SU group
\newcommand{\Spql}[1]{\Sp(#1)/\Kumz}

\newcommand{\SUql}[1]{\SU(#1)/\Kumgp}

\newcommand{\MM}{\mathbb{M}}    %Monster
\newcommand{\Co}{\operatorname{\textsl{Co}}}    %Conway 1

\newcommand{\Ru}{\operatorname{\textsl{Ru}}}
\newcommand{\cRu}{2.\!\Ru}
\newcommand{\Rtd}{\widetilde{R}}

\newcommand{\Cl}{\operatorname{Cliff}}
\newcommand{\Cm}{\operatorname{CM}}
\newcommand{\atw}[1]{\tilde{A}(#1)}     %Twisted Clifford mod VOSA

\newcommand{\vn}{V^{\natural}}     %Moonshine VOA
\newcommand{\afn}{{A^{f\natural}}}    %Fermionic \cdot 1 VOSA
\newcommand{\aru}{A_{{\Ru}}}  %Rudvalis VOhmA
\newcommand{\au}{A(\gt{u})}

\newcommand{\cgset}{\Omega}   %Set of conformal generators
\newcommand{\cgU}{\cgset_{\U}}
\newcommand{\cgRu}{\cgset_{\Ru}}
\newcommand{\cfalg}{\mc{A}}   %algebra determined by the enhanced
\newcommand{\vemp}{\nu_{\circ}}     %Conformal element
\newcommand{\vful}{\nu_{\bullet}}   %Conformal element
\newcommand{\vp}{\nu_+}
\newcommand{\Vp}{V_+}
\newcommand{\vm}{\nu_-}
\newcommand{\Vm}{V_-}

\newcommand{\Vpm}{V_{\pm}}

\newcommand{\cej}{\jmath}
\newcommand{\ceJ}{J}

\newcommand{\spsp}{\varrho}

\newcommand{\ogi}{\theta}   %involution orbifold element
\newcommand{\ogz}{\gt{z}}
\newcommand{\ainv}{\vartheta} %Anti-involution
\newcommand{\ogih}{\ogi^{1/2}}%\ii on \gt{a} and -\ii on \gt{a}^*

\newcommand{\vac}{{\bf 1}}      %Vacuum idem
\newcommand{\cas}{{\bf \omega}} %Conformal element
\newcommand{\scas}{{\bf \tau}}  %Super conformal element

\newcommand{\uu}{{\bf 1}}   %Bold unit
\newcommand{\trip}[3]{   \begin{tabular}{cc}
                         \multicolumn{2}{c}{$#1$}\\
                         $#2$ & $#3$
                         \end{tabular}}

\newtheorem{thm}{Theorem}[section]
\newtheorem{cor}[thm]{Corollary}
\newtheorem{lem}[thm]{Lemma}
\newtheorem{prop}[thm]{Proposition}
\newtheorem*{conj}{Conjecture}

\theoremstyle{definition}
\newtheorem*{defn}{Definition}

\theoremstyle{remark}
\newtheorem*{rmk}{Remark}

\numberwithin{equation}{subsection}
\newtheorem*{eg}{Example}

\theoremstyle{plain}
\newtheorem*{sthm}{Theorem}

\begin{document}

\title{
    \textsc{Moonshine for Rudvalis's sporadic group I}
     \footnote{{\it Mathematics Subject Classification (2000)}
                %Primary
                17B69, %; Secondary
                20D08.
                }
          }

\author{John F. Duncan
          \footnote{
          Harvard University,
          Department of Mathematics,
          One Oxford Street,
          Cambridge, MA 02138,
          U.S.A.}
          {}\footnote{
          Email: {\tt duncan@math.harvard.edu};\;
          homepage: {\tt
          http://math.harvard.edu/\~{}jfd/}
               }
               }

\date{October 18, 2008}
%\date{\today}

\setcounter{section}{-1}

%\subjclass[2000]{Primary 17B69; Secondary 20D08.}

\maketitle

\begin{abstract}
We introduce the notion of vertex operator superalgebra with
enhanced conformal structure, which is a refinement of the notion of
vertex operator superalgebra. We exhibit several examples, including
a particular one which is self-dual, and whose full symmetry group
is a direct product of a cyclic group of order seven with the
sporadic simple group of Rudvalis. We thus obtain an analogue of
Monstrous Moonshine for a sporadic group not involved in the
Monster. Two variable analogues of the usual McKay--Thompson series
are naturally associated to the action of the Rudvalis group on this
object, and we provide explicit expressions for all the series
arising.
\end{abstract}

\tableofcontents
%\listoftables

\section{Introduction}\label{sec:intro}

We are interested in the problem of realizing sporadic simple groups
as symmetry groups of vertex operator algebras (VOAs).

Perhaps the most striking example of this phenomena is the Moonshine
VOA $\vn$ of \cite{FLM} which provides a realization of the Monster
sporadic group $\MM$. A feature of this example is that the
Moonshine VOA conjecturally admits a characterization as the unique
self-dual VOA of rank $24$ with no (non-trivial) small (degree $1$)
elements. Assuming the conjecture is true, we then obtain an
interesting description of the Monster group: as the symmetry group
of this natural object in the category of VOAs, the unique self-dual
VOA of rank $24$ with no small elements.

Another feature of this example is the extent to which it elucidates
a connection between the Monster group and the theory of modular
forms; what has been referred to as Monstrous Moonshine since the
historic article \cite{ConNorMM}. In particular, the fact that the
Monster group acts as automorphisms of the Moonshine VOA allows one
to attach a power series, a kind of $\ZZ$-graded character called
the McKay--Thompson series, to each conjugacy class of $\MM$. It is
a non-trivial property of VOAs that these characters can be regarded
as holomorphic functions on the upper half plane with certain
modular properties \cite{ZhuPhd}. The fact that the Moonshine VOA is
self-dual (that is, has no irreducible modules other than itself)
and has rank $24$, entails that these functions are actually modular
invariants with very special properties, although the theory of VOAs
does perhaps not yet provide a full explanation for all the
properties these series have been conjectured \cite{ConNorMM}, and
subsequently proven \cite{BorMM} to satisfy.

A second example was given in \cite{HohnPhD} where it was shown how
one may modify the Moonshine VOA $\vn$ so as to arrive at a vertex
operator superalgebra (VOSA) $V\!B^{\natural}$ whose full
automorphism group is the direct product of a group of order two
with the Baby Monster group $\mathbb{B}$. In this case also, one has
a conjectural characterization for the VOSA $V\!B^{\natural}$, as
the unique self-dual VOSA of rank $23\tfrac{1}{2}$ with no
(non-trivial) small (degree $1$ or degree ${1}/{2}$) elements, and
thus, conjecturally, we obtain an analogous description of the Baby
Monster group: as the symmetry group of this object modulo its
center. One may also consider the McKay--Thompson series associated
to the Baby Monster group via its action on $V\!B^{\natural}$, so
that one obtains another point of interaction between a sporadic
simple group and the theory of modular forms.

Further examples of VOAs with finite simple symmetry groups have
appeared in the literature. A VOA whose full automorphism group is
$O^+_{10}(2)$ is studied in \cite{GriO}, and the automorphism groups
of closely related VOAs are determined in \cite{ShiAutmGpVLpnorts}
and \cite{ShiAutmGpVLpgen}. Results regarding the realization of
$3$-transposition groups as VOA symmetries have appeared in
\cite{KitMiy3transVOAauts} and \cite{Mat3TransSympVOAs}. Passing to
positive characteristic, modular representations for certain
sporadic groups involved in the Monster are constructed via VOAs
over finite fields in \cite{BorRybModuMII}. This work arose from
modular analogues of the Moonshine conjectures that were introduced
in \cite{RybModuM}, and the proof of these conjectures was completed
in \cite{BorModuMIII}.

\medskip

Evidence for a third example involving a sporadic simple group in
characteristic zero appeared earlier in \cite{FLMBerk} (see also
\cite{BorRybModuMII}) and suggested a possible realization by vertex
operators of the largest sporadic group of Conway $\Co_1$. We
considered this example in some detail in \cite{Dun_VACo} and found
that although the Conway group does indeed act as symmetries, the
full automorphism group of the VOSA structure alone is not finite.
On the other hand, for a particular action of the Conway group on
this VOSA there is a unique vector in the degree ${3}/{2}$ space
that is invariant for this group. What is more, the vertex operator
associated to this element is such that its Fourier components,
along with those of the usual conformal element (Virasoro element),
generate an action of the $N=1$ Neveu--Schwarz super Virasoro
algebra, a natural super analogue of the Virasoro algebra, acting on
this VOSA. In \cite{Dun_VACo} we prove that the symmetry group
fixing both the usual conformal element and this new superconformal
element is precisely Conway's largest sporadic group. There we
denote this object by $\afn$, and we call it an $N=1$ VOSA in order
to indicate that we regard it as one example from a family of vertex
operator superalgebras for which the super Virasoro action is to be
taken as axiomatic.

In direct analogy with $\vn$ and $V\!B^{\natural}$, one is quickly
lead to conjecture the following characterization of $\afn$: the
unique self dual $N=1$ VOSA of rank $12$ with no (non-trivial) small
(degree ${1}/{2}$) elements. We establish this conjecture, up to
some technical conditions, in \cite{Dun_VACo}.

Beyond ordinary VOSAs, $\afn$ is our first example of what we now
refer to as VOSAs with enhanced conformal structure, or more
briefly, enhanced VOSAs. The idea is that we refine the notion of
VOSA by imposing certain extra algebraic structure arising from
distinguished vectors beyond the usual Virasoro and vacuum elements,
and a precise definition is given in \S\ref{sec:eVOAs}. One purpose
of the present article is to motivate the notion of enhanced VOSA by
furnishing further examples.

Our main example is a self-dual enhanced VOSA $\aru$ of rank $28$
whose full automorphism group is a direct product of a cyclic group
of order seven with the sporadic simple group of Rudvalis. In
particular, the enhanced VOSA $\aru$ furnishes a realization by
vertex operators of one of the six sporadic groups not involved in
the Monster.

In analogy with case of the Moonshine VOA and the Monster group we
obtain a kind of moonshine for the Rudvalis group by considering the
graded traces of elements of this group acting on $\aru$. In fact,
the enhanced conformal structure on $\aru$ is such that $\aru$
admits a structure of (what we call) $\U(1)$--VOSA, and to such
objects one can naturally associate McKay--Thompson series in two
variables, which specialize to the ordinary McKay--Thompson series
after taking a suitable limit. The action of the Rudvalis group on
$\aru$ is sufficiently transparent that it is a straight forward
task to provide explicit expressions for all the two variable
McKay--Thompson series arising. We find that they are Jacobi forms
for congruence subgroups of the modular group $\SL_2(\ZZ)$.

It will certainly be interesting to study the series associated to
the Rudvalis group in more detail, and compare them to those
associated to the Monster group via the Moonshine VOA, and to
Conway's group via the enhanced VOSA studied in \cite{Dun_VACo}.
%The
%functions of Monstrous Moonshine have been characterized in
%\cite{ConMcKSebDiscGpsM}.

\medskip

The main result of this article is the following theorem.
\begin{sthm}[\ref{thm:Ru:main}]
The quadruple $\aru=(\aru,Y,\vac,\cgRu)$ is a self-dual enhanced
$\U(1)$--VOSA of rank $28$. The full automorphism group of
$(\aru,\cgRu)$ is a direct product of a cyclic group of order seven
with the sporadic simple group of Rudvalis.
\end{sthm}
In light of the above discussion, it is natural to hope that VOA
theory might be used to furnish convenient characterizations of
certain finite simple groups. In the case of the Conway group and
$\afn$, the enhanced conformal structure was essential to the
characterization, and we regard this as another motivation for the
notion of enhanced VOSA.

Later in the article we conjecture a uniqueness result for $\aru$,
which is analogous to those which hold for the Golay code, the Leech
lattice (see \cite{Con69}), and the $N=1$ VOSA for Conway's group,
and those which are conjectured to hold for the Moonshine VOA and
the Baby Monster VOSA. All these objects have sporadic automorphism
groups. The object $\aru$ is a first example with non-Monstrous
sporadic automorphism group.

\subsection{Outline}\label{sec:intro:outline}

In \S\ref{sec:VOAs} we review some facts about VOAs\footnote{Here,
and from here on, we will usually suppress the ``super'' in
super-objects, so that unless extra clarification is necessary,
superspaces and superalgebras will be referred to as spaces and
algebras, respectively, and the term VOA for example, will be used
even when the underlying vector space comes equipped with a
$\ZZ/2$-grading.}. We recall some basics from the formal calculus in
\S\ref{sec:VOAs:FmlCalc}, and then the definition of VOAs in
\S\ref{sec:VOAs:struc}. We also recall from \cite{FLM} the higher
order brackets defined on a VOA in \S\ref{sec:VOAs:nbrackets}. These
higher order brackets are generalizations of the ordinary Lie
bracket, and are natural in the context of VOA theory, and to some
extent they motivate the notion of vertex Lie algebra which we
review in \S\ref{sec:VOAs:VLAs}. The notion of vertex Lie algebra
plays a role in the definition of enhanced VOA, and this in turn
appears in \S\ref{sec:eVOAs}.

In \S\ref{sec:cliffalgs:struc} we review our conventions regarding
Clifford algebras; in \S\ref{sec:cliffalgs:spin} the Spin groups,
and in \S\ref{sec:cliffalgs:mods}, Clifford algebra modules. In
\S\ref{sec:cliffalgs:VOAs} we review the well-known construction of
VOA structure on certain infinite dimensional Clifford algebra
modules constructed in turn from a finite dimensional vector space
with symmetric bilinear form, and in \S\ref{sec:cliffalgs:Herm} we
make some conventions regarding the Hermitian structures arising in
the case that the initial vector space comes equipped also with a
suitable Hermitian form.

In \S\ref{sec:GLnq} we present a family of enhanced VOAs whose full
automorphism groups are the general linear groups $\GL_N(\CC)$. We
later require to realize symmetry groups no larger than
$\SL_N(\CC)$, and for this it is useful to consider a twisted
analogue of the Clifford module VOA construction given in
\S\ref{sec:cliffalgs:VOAs}. We review such a construction in
\S\ref{sec:GLnq:tw}, and then in \S\ref{sec:GLnq:SLgps} we consider
the specific example of an enhanced VOA with symmetry group of the
form $\SL_{28}(\CC)/\Kumgp$.

In \S\ref{sec:Ru} we construct an enhanced VOA whose automorphism
group is a sevenfold cover of the sporadic simple group of Rudvalis.
We give two constructions of the enhanced conformal structure. The
first construction, in \S\ref{sec:Ru:geom}, arises directly from the
geometry of the Conway--Wales lattice, which was introduced in
\cite{ConRu}. This approach is more brief and more conceptual than
the second construction, in \S\ref{sec:Ru:mnml}, which is
nonetheless more convenient for (computer aided) computations. The
approach of \S\ref{sec:Ru:mnml} involves analysis of a certain
monomial group that is closely related to the Cayley algebra (or
integral octonion algebra) structure with which one may equip the
$E_8$ lattice. We review this structure in
\S\ref{sec:Ru:mnml:cayleyalg}, then the monomial group is
constructed in \S\ref{sec:Ru:mnml:monomialgp}, and we determine in
\S\ref{sec:Ru:mnml:monomialinvs} the invariants of this monomial
group under its action on a certain space. Sections
\ref{sec:Ru:geom} and \ref{sec:Ru:mnml} are independent. Picking up
from the end of either one, the construction of the enhanced VOA for
the Rudvalis group is given in \S\ref{sec:Ru:cnst}. The symmetry
group is considered in \S\ref{sec:Ru:symms}, and we discuss a
conjectural characterization of the enhanced VOA for the Rudvalis
group at the very end of this section.

In \S\ref{sec:series} we furnish explicit expressions for the two
variable McKay--Thompson series associated to $\aru$. In particular,
we introduce the notion of two variable McKay--Thompson series in
\S\ref{sec:series:2var}, and we consider the two variable series
arising for the Rudvalis group in \S\ref{sec:series:MbeyM}. The
Appendix contains tables displaying the first few coefficients of
each two variable McKay--Thompson series arising from the action of
the Rudvalis group on $\aru$.

%%%%%%%%%%%%%%%%%%%%%%%%%%%%%%%%%%%%%%%%%%%%%%%%%%%%%%%%%

\subsection{Notation}\label{sec:intro:notation}

We choose a square root of $-1$ in $\CC$ and denote it by $\ii$. For
$q$ a prime power $\FF_q$ shall denote a field with $q$ elements. We
use $\FF^{\times}$ to denote the non-zero elements of a field $\FF$.
More generally, $A^{\times}$ shall denote the set of invertible
elements in an algebra $A$. For $G$ a group and $\FF$ a field we
write $\FF G$ for the group algebra of $G$ over $\FF$. For the
remainder we shall use $\FF$ to denote either $\RR$ or $\CC$.

%\medskip

For $\fset$ a finite set, we denote the power set of $\fset$ by
$\mc{P}(\fset)$. The set operation of symmetric difference equips
$\mc{P}(\fset)$ with a structure of $\FF_2$-vector space, and with
this structure in mind, we sometimes write $\FF_2^{\fset}$ in place
of $\mc{P}(\fset)$. Suppose that $\fset$ has $N$ elements. The space
$\FF_2^{\fset}$ comes equipped with a function $\FF_2^{\fset}\to
\{0,1,\ldots,N\}$ called {\em weight}, which assigns to an element
$\gamma\in\FF_2^{\fset}$ the cardinality of the corresponding
element of $\mc{P}(\fset)$. An $\FF_2$-subspace of $\FF_2^{\fset}$
is called a {\em binary linear code of length $N$}. %We will also
%encounter {\em ternary codes}, which are linear subspaces of
%$\FF_3^{\fset}$. In this later case we use the term weight to mean
%the number of non-zero entries in a given codeword
%$C\in\FF_3^{\fset}$.
A binary linear code $\mc{C}$ is called {\em even} if all its
codewords have even weight, and is called {\em doubly even} if all
codewords have weight divisible by $4$. The space $\FF_2^{\fset}$
carries a bilinear form defined by setting $\lab
C,D\rab=\sum_{i}\gamma_i\delta_i$ for $C=(\gamma_i)$ and
$D=(\delta_i)$, so that given a binary linear code $\mc{C}$ we may
consider the dual code $\mc{C}^{\circ}$ defined by
\begin{gather}
     \mc{C}^{\circ}=\left\{D\in\FF_2^{\fset}\mid \lab C,D\rab=0
          ,\,\forall C\in\mc{C}\right\}
\end{gather}
We say that $\mc{C}$ is {\em self-dual} in the case that
$\mc{C}=\mc{C}^{\circ}$, and we say that $\mc{C}$ is {\em
self-orthogonal} when $\mc{C}\subset\mc{C}^{\circ}$. We write
$\mc{C}^*$ for the {\em co-code} $\mc{C}^*=\FF_p^{\fset}/\mc{C}$.

%\medskip

A {\em superspace} is a vector space with a grading by
$\ZZ/2=\{\bar{0},\bar{1}\}$. When $M$ is a superspace, we write
$M=M_{\bar{0}}\oplus M_{\bar{1}}$ for the superspace decomposition,
and for $u\in M$ we set $|u|=\gamma\in\{\bar{0},\bar{1}\}$ when $u$
is $\ZZ/2$-homogeneous and $u\in M_{\gamma}$. The dual space $M^*$
has a natural superspace structure such that
$(M^*)_{\gamma}=(M_{\gamma})^*$ for $\gamma\in \{\bar{0},\bar{1}\}$.
The space $\End(M)$ admits a structure of Lie superalgebra when
equipped with the Lie superbracket $[\cdot\,,\cdot]$ which is
defined so that $[a,b]=ab-(-1)^{|a||b|}ba$ for $\ZZ/2$-homogeneous
$a,b$ in $\End(M)$.

Almost all vector spaces in this article will be most conveniently
regarded as superspaces (with possibly trivial odd parts) and
similarly for algebras, and so from this point onwards, and unless
otherwise qualified we will use the terms ``space'' and ``algebra''
as inclusive of the notions of superspace and superalgebra,
respectively. Thus we will speak of vertex operator algebras, and
Lie algebras, and so on, even when the underlying vector spaces are
$\ZZ/2$ graded.

When $z$ denotes a formal variable, we write $M[[z]]$ for the space
of formal Taylor series with coefficients in $M$, and $M((z))$ for
the space of formal Laurent series with coefficients in $M$. All
formal variables in this article will be regarded as even, so that
the superspace decomposition of $M[[z]]$ is $M_{\bar{0}}[[z]]\oplus
M_{\bar{1}}[[z]]$, and similarly for $M((z))$. We write
$\bigwedge(M)$ for the full exterior algebra of a vector space $M$.
We write $\bigwedge(M)=\bigwedge(M)^0\oplus\bigwedge(M)^1$ for the
parity decomposition of $\bigwedge(M)$, and we write
$\bigwedge(M)=\bigoplus_{k\geq 0}\bigwedge^k(M)$ for the natural
$\ZZ$-grading on $\bigwedge(M)$. We denote by $D_z$ the operator on
formal power series which is formal differentiation in the variable
$z$, so that if $f(z)=\sum c_mz^{-m-1}\in M((z))$ is a formal
power-series with coefficients in some space $M$, we have
$D_zf(z)=\sum(-m)f_{m-1}z^{-m-1}$. For $m$ a non-negative integer,
we set $D_z^{(m)}=\tfrac{1}{m!}D_z^m$.

%\medskip

As is customary, we use $\eta(\tau)$ to denote the Dedekind eta
function.
\begin{gather}\label{eqn:intro:notation:eta}
        \eta(\tau)=q^{1/24}\prod_{n= 1}^{\infty}(1-q^n)
\end{gather}
Here $q=e^{2\pi\ii\tau}$ and $\tau$ is a variable in the upper half
plane, which we denote by $\hh$. Recall also the Jacobi theta
function $\vartheta_3(\xi|\tau)$ defined by
\begin{gather}\label{eqn:intro:notation:theta}
    \vartheta_3(\xi|\tau)=\sum_{m\in\ZZ}
        e^{2\ii\xi m+\pi\ii\tau m^2}
\end{gather}
for $\tau\in\hh$ and $\xi\in\CC$. According to the Jacobi Triple
Product Identity we have
\begin{gather}\label{eqn:intro:tripprod}
    \vartheta_3(\xi|\tau)
    =\prod_{m\geq 0}(1-q^{m+1})
        (1+e^{2\ii\xi}q^{m+1/2})
        (1+e^{-2\ii\xi}q^{m+1/2})
\end{gather}
with $q=e^{2\pi\ii\tau}$ as before.

The modular group $\bar{\Gamma}=\PSL_2(\ZZ)$ acts in a natural way
on the upper half plane $\hh$. This action is generated by the
modular transformations $\tau\mapsto-1/\tau$ and $\tau\mapsto
\tau+1$, which we denote by $S$ and $T$, respectively, and these
generators are subject to the relations $S^2=(ST)^3=\Id$. We write
$\bar{\Gamma}_{\ogi}$ for the theta group: the subgroup of
$\bar{\Gamma}$ generated by $S$ and $T^2$.

\medskip

The most specialized notations arise in \S\ref{sec:cliffalgs}. We
include here a list of them, with the relevant subsections indicated
in brackets. They are grouped roughly according to similarity of
appearance, rather than by order of appearance, so that the list may
be easier to search through, whenever the need might arise.
\begin{small}
\begin{list}{}{     \itemsep -3pt
                    \topsep 3pt
                         }
\item[$\gt{a}$]     A complex vector space with non-degenerate
Hermitian form (\S\ref{sec:cliffalgs:Herm}).

\item[$\gt{a}^*$]   The dual space to $\gt{a}$
(\S\ref{sec:cliffalgs:Herm}).

\item[$\gt{u}$]  A real or complex vector space of even dimension with
non-degenerate bilinear form, assumed to be positive definite in the
real case (\S\ref{sec:cliffalgs:struc}). In the case that
$\gt{u}=\gt{a}\oplus \gt{a}^*$, the bilinear form is assumed to be
$1/2$ times the symmetric linear extension of the natural pairing
between $\gt{a}$ and $\gt{a}^*$ (\S\ref{sec:cliffalgs:Herm}).

\item[$\{a_i\}_{i\in\hfset}$]  A basis for $\gt{a}$, orthonormal in
the sense that $(a_i,a_j)=\delta_{ij}$ for $i,j\in\hfset$
(\S\ref{sec:cliffalgs:Herm}).

\item[$\{a_i^*\}_{i\in\hfset}$]    The dual basis to
$\{a_i\}_{i\in\hfset}$ (\S\ref{sec:cliffalgs:Herm}).

\item[$\{e_i\}_{i\in\fset}$]     A basis for $\gt{u}$,
orthonormal in the sense that $\lab e_i,e_j\rab=\delta_{ij}$ for
$i,j\in\fset$ (\S\ref{sec:cliffalgs:struc}). In the case that
$\gt{u}=\gt{a}\oplus\gt{a}^*$ we set $\fset=\hfset\cup\hfset'$, and
we insist that $e_i=a_i+a_i^*$ and $e_{i'}=\ii(a_i-a_i^*)$ for
$i\in\hfset$ (\S\ref{sec:cliffalgs:Herm}).

\item[$e_I$]   We write $e_I$ for
$e_{i_1}\cdots e_{i_k}\in\Cl(\gt{u})$ when $I=\{i_1,\ldots,i_k\}$ is
a subset of $\fset$ and $i_1<\cdots<i_k$
(\S\ref{sec:cliffalgs:struc}).

\item[$g(\cdot)$]   We write $g\mapsto
g(\cdot)$ for the natural homomorphism $\Sp(\gt{u})\to\SO(\gt{u})$.
Regarding $g\in\Sp(\gt{u})$ as an element of $\Cl(\gt{u})^{\times}$
we have $g(u)=gug^{-1}$ in $\Cl(\gt{u})$ for $u\in\gt{u}$. More
generally, we write $g(x)$ for $gxg^{-1}$ when $x$ is any element of
$\Cl(\gt{u})$.

\item[$e_I(\cdot)$] When $I$ is even, $e_I\in\Cl(\gt{u})$ lies also
in $\Sp(\gt{u})$, and $e_I(x)$ denotes $e_Ixe_I^{-1}$ when
$x\in\Cl(\gt{u})$.

\item[$\ogz$]  We denote $e_{\fset}\in\Sp(\gt{u})$ also by $\ogz$
(\S\ref{sec:cliffalgs:spin}).

\item[$\ogi$]  The map which is $-\Id$ on $\gt{u}$, or the parity
involution on $\Cl(\gt{u})$ (\S\ref{sec:cliffalgs:struc}), or the
parity involution on $A(\gt{u})_{\Theta}$
(\S\ref{sec:cliffalgs:VOAs}).

%\item[$\ainv$] An anti-linear involution defined on $\gt{u}$ and on
%$A(\gt{u})_{\Theta}$ in the case that $\gt{u}=\gt{a}\oplus\gt{a}^*$
%(\S\ref{sec:cliffalgs:Herm}).

\item[$\ogih$] The map which is $\ii\Id$ on $\gt{a}$ and $-\ii\Id$
on $\gt{a}^*$, or a lift of this map to $\Cl(\gt{u})$ or
$A(\gt{u})_{\Theta}$ (\S\ref{sec:cliffalgs:Herm}).

\item[$1_E$]   A vector in $\Cm(\gt{u})_E$ such that $x1_E=1_E$ for
$x$ in $E$ (\S\ref{sec:cliffalgs:mods}).

\item[$\vac_E$]     The vector corresponding to $1_E$ under the
identification between $\Cm(\gt{u})_{E}$ and
$(A(\gt{u})_{\ogi,E})_{N/8}$ when $\gt{u}$ has dimension $2N$
(\S\ref{sec:cliffalgs:VOAs}).

\item[$\hfset$]     A finite ordered set indexing an orthonormal
basis for $\gt{a}$ (\S\ref{sec:cliffalgs:Herm}).

\item[$\hfset'$]    A second copy of $\hfset$ with the natural
identification $\hfset\leftrightarrow\hfset'$ denoted
$i\leftrightarrow i'$ for $i\in\hfset$ (\S\ref{sec:cliffalgs:Herm}).

\item[$\fset$] A finite ordered set indexing an orthonormal basis
for $\gt{u}$ (\S\ref{sec:cliffalgs:struc}). In the case that
$\gt{u}=\gt{a}\oplus\gt{a}^*$, we set $\fset=\hfset\cup\hfset'$ and
insist that $\fset$ be ordered according to the ordering on $\hfset$
and the rule $i<j'$ for $i,j\in\hfset$ (\S\ref{sec:cliffalgs:Herm}).

\item[$\mc{E}$]     A label for the basis
$\{e_i\}_{i\in\fset}$ (\S\ref{sec:cliffalgs:struc}).

\item[$E$]     A subgroup of $\Cl(\gt{u})^{\times}$ homogeneous with
respect to the $\FF_2^{\mc{E}}$ grading on $\Cl(\gt{u})$
(\S\ref{sec:cliffalgs:mods}).

\item[$X$]     In the case that $\gt{u}=\gt{a}\oplus\gt{a}^*$, we
take $E=X$ in the above (and the below), where $X$ is the subgroup
of $\Cl(\gt{u})$ generated by elements $\ii e_ie_{i'}$ for
$i\in\hfset$ (\S\ref{sec:cliffalgs:Herm}).

%\item[${_{\RR}V}$]  The real vector space consisting of
%the $\ainv$-fixed points in $V$ when $V$ is a vector space with an
%action by $\ainv$ (\S\ref{sec:cliffalgs:Herm}).

\item[$A(\gt{u})$]  The Clifford module VOA associated to the
vector space $\gt{u}$ (\S\ref{sec:cliffalgs:VOAs}).

\item[$\atw{\gt{u}}$]    The twisted Clifford module VOA
associated to $\gt{u}$, obtained by taking $\ogi$-fixed points of
$A(\gt{u})_{\Theta}$ (\S\ref{sec:GLnq:tw}).

\item[$A(\gt{u})_{\ogi}$]     The canonically $\ogi$-twisted module
over $A(\gt{u})$ (\S\ref{sec:cliffalgs:VOAs}).

\item[$A(\gt{u})_{\ogi,E}$]   The $\ogi$-twisted module
$A(\gt{u})_{\ogi}$ realized in such a way that the subspace of
minimal degree may be identified with $\Cm(\gt{u})_E$
(\S\ref{sec:cliffalgs:VOAs}).

\item[$A(\gt{u})_{\Theta}$]   The direct sum of $A(\gt{u})$-modules
$A(\gt{u})\oplus A(\gt{u})_{\ogi}$ (\S\ref{sec:cliffalgs:VOAs}).

\item[$\mc{C}(E)$]  The binary linear code on $\fset$ consisting of
elements $I$ in $\FF_2^{\fset}$ for which $E$ has non-trivial
intersection with $\FF e_I\subset\Cl(\gt{u})$
(\S\ref{sec:cliffalgs:mods}).

\item[$\Cm(\gt{u})_E$]   The module over $\Cl(\gt{u})$ induced from
a trivial module over $E$ (\S\ref{sec:cliffalgs:mods}).

\item[$\Cl(\gt{u})$]     The Clifford algebra associated to the vector space
$\gt{u}$ (\S\ref{sec:cliffalgs:struc}).

\item[$\Sp(\gt{u})$]     The spin group associated to the vector space
$\gt{u}$ (\S\ref{sec:cliffalgs:spin}).

\item[$\lab\cdot\,,\cdot\rab$]     A non-degenerate symmetric
bilinear form on $\gt{u}$ or on $\Cl(\gt{u})$
(\S\ref{sec:cliffalgs:struc}), or on the $\Cl(\gt{u})$-module
$\Cm(\gt{u})_E$ (\S\ref{sec:cliffalgs:mods}). In the case that
$\gt{u}$ is real all of these forms will be positive definite.

\item[$(\cdot\,,\cdot)$] A non-degenerate Hermitian form on
$\gt{a}$, or on $\Cl(\gt{u})$ or $\Cm(\gt{u})_X$
 in the case that $\gt{u}=\gt{a}\oplus\gt{a}^*$
 (\S\ref{sec:cliffalgs:Herm}). The Hermitian forms arising will
always be anti-linear in the right hand slot.%, and will always
%satisfy $(u,v)=\lab u,\ainv v\rab$.

\item[$\lab\cdot|\cdot\rab$]  A non-degenerate symmetric bilinear
form on $A(\gt{u})_{\Theta}$ (\S\ref{sec:cliffalgs:VOAs}).

%\item[$(\cdot\,|\cdot)$]   In the case that
%$\gt{u}=\gt{a}\oplus\gt{a}^*$, a Hermitian form is defined on
%$A(\gt{u})_{\Theta}$ by setting $(u|v)=\lab u|\ainv v\rab$
%(\S\ref{sec:cliffalgs:Herm}).
\end{list}
\end{small}

\section{Vertex operator algebras}\label{sec:VOAs}

In \S\ref{sec:VOAs:struc} we review the definition of VOA, after
recalling some facts from the formal calculus in
\S\ref{sec:VOAs:FmlCalc}. In \S\ref{sec:VOAs:nbrackets} we recall
from \cite{FLM} the higher order generalizations of the Lie bracket,
which arise naturally in the context of VOAs, and in
\S\ref{sec:VOAs:VLAs} we recall some facts about vertex Lie algebras
(also known as conformal algebras).

\subsection{Formal calculus}\label{sec:VOAs:FmlCalc}

Given a rational function $f(z,w)$ with poles only at $z=0$, $w=0$,
and $z=w$, we write $\iota_{z,w}f(z,w)$, $\iota_{w,z}f(z,w)$, and
$\iota_{w,z-w}f(z,w)$ for the power series expansions of $f(z,w)$ in
the respective domains: $|z|>|w|>0$, $|w|>|z|>0$, and $|w|>|z-w|>0$.
Then in the case that $f(z,w)=z^mw^n(z-w)^l$ for some $m,n,l\in\ZZ$
for example, we have
\begin{align}
     \iota_{z,w}f(z,w)&=\sum_{k\geq 0}(-1)^k\binom{l}{k}
          z^{m+l-k}w^{n+k}\\
     \iota_{w,z}f(z,w)&=\sum_{k\geq 0}(-1)^{l+k}\binom{l}{k}
          z^{m+k}w^{n+l-k}\\
     \iota_{w,z-w}f(z,w)&=\sum_{k\geq 0}\binom{m}{k}w^{m+n-k}(z-w)^{l+k}
\end{align}
For $f(z_1,\ldots,z_k)=\sum c_{m_1,\ldots,m_k}z_1^{-m_1-1}\cdots
z_k^{-m_k-1}$ a formal power series in variables $z_i$, we write
$\Res_{z_i}f$ for the coefficient of $z_i^{-1}$ in
$f(z_1,\ldots,z_k)$, so that
\begin{gather}
     \Res_{z_1}f(z_1,\ldots,z_k)=\sum
          c_{0,m_2,\ldots,m_k}z^{-m_2-1}\cdots z^{-m_k-1}
\end{gather}
for example, and we write $\Sing f$ for the sum of terms involving
only negative powers of all the variables $z_i$.
\begin{gather}
     \Sing f(z_1,\ldots,z_k)=\sum_{m_1,\ldots,m_k\geq 0}
          c_{m_1,\ldots,m_k}z^{-m_1-1}\cdots z^{-m_k-1}
\end{gather}
Suppose again that $f(z,w)$ is a rational function with possible
poles at $z=0$, $w=0$, and $z=w$. Then we have the following special
case of the Cauchy Theorem
\begin{gather}\label{eqn:VOAs:FmlCalc:AnCauchy}
     \int_{C_R(0)}f(z,w){\rm d}z
     -\int_{C_r(0)}f(z,w){\rm d}z=
     \int_{C_{\epsilon}(w)}f(z,w){\rm d}z
\end{gather}
where $C_a(z_0)$ denotes a positively oriented circular contour (in
the $z$-plane) of radius $a$ about $z_0$, and $R$, $r$ and
$\epsilon$ are chosen so that $R>|w|>r>0$, and $\epsilon<{\rm
min}\{R-|w|,|w|-r\}$. The following identity is then an algebraic
reformulation of (\ref{eqn:VOAs:FmlCalc:AnCauchy}).
\begin{gather}\label{eqn:VOAs:FmlCalc:AlgCauchy}
     \Res_{z}\iota_{z,w}f(z,w)-\Res_{z}\iota_{w,z}f(z,w)
     =\Res_{z-w}\iota_{w,z-w}f(z,w)
\end{gather}

\subsection{Structure}\label{sec:VOAs:struc}

For a {\em vertex algebra} structure on a (super)space
$U=U_{\bar{0}}\oplus U_{\bar{1}}$ over a field $\FF$ we require the
following data.
\begin{itemize}
\item{\em Vertex operators:}   an even morphism $Y:U\otimes U\to
U((z))$ such that when we write
$Y(u,z)v=\sum_{n\in\ZZ}u_{(n)}vz^{-n-1}$, we have $Y(u,z)=0$ only
when $u=0$.
\item{\em Vacuum:}   a distinguished vector $\vac\in U_{\bar{0}}$
such that $Y(\vac,z)u=u$ for $u\in U$, and $Y(u,z)\vac|_{z=0}=u$.
%\item{\em Translation operator:}   a superspace morphism $T:U\to
%U$ such that $Y(Tu,z)=D_zY(u,z)$ for all $u\in U$.
\end{itemize}
This data furnishes a vertex algebra structure on $U$ just when the
following identity is satisfied.
\begin{itemize}
\item{\em Jacobi identity:}   for $\ZZ/2$ homogeneous $u,v\in U$,
and for any $m,n,l\in\ZZ$ we have
\begin{gather}\label{eqn:VOAs:struc:Jac}
     \begin{split}
     &\Res_{z}Y(u,z)Y(v,w)\iota_{z,w}F(z,w)\\
     &-\Res_{z}(-1)^{|u||v|}Y(v,w)Y(u,z)\iota_{w,z}F(z,w)\\
     &=\Res_{z-w}Y(Y(u,z-w)v,w)\iota_{w,z-w}F(z,w)
     \end{split}
\end{gather}
where $F(z,w)=z^mw^n(z-w)^l$.
\end{itemize}
We denote such an object by the triple $(U,Y,\vac)$.

A {\em vertex operator algebra (VOA)} is vertex algebra $(U,Y,\vac)$
equipped with a distinguished element $\cas\in\U_{\bar{0}}$ called
the {\em Virasoro element}, such that $L(-1):=\cas_{(0)}$ satisfies
\begin{gather}
     [L(-1),Y(u,z)]=D_zY(u,z)
\end{gather}
for all $u\in U$, and such that the operators $L(n):=\cas_{(n+1)}$
furnish a representation of the Virasoro algebra on the vector space
underlying $U$, so that we have
\begin{gather}\label{eqn:VOAs:struc:VirRel}
    \left[L(m),L(n)\right]=(m-n)L(m+n)
    +\frac{m^3-m}{12}\delta_{m+n,0}{ c}{\rm Id}
\end{gather}
for some $c\in\FF$. We also require the following grading condition.
\begin{itemize}
\item{\em $L(0)$-grading:}  the action of $L(0)$ on $U$ is
diagonalizable with rational eigenvalues bounded from below, and is
such that the $L(0)$-homogeneous subspaces $U_n:=\{u\in U\mid
L(0)u=nu\}$ are finite dimensional.
\end{itemize}
When these conditions are satisfied we write $(U,Y,\vac,\cas)$ in
order to indicate the particular data that constitutes the VOA
structure on $U$. %In the case that $U=U_{\bar{0}}$ we are speaking
%of ordinary vertex operator algebras (VOAs).
%By definition the coefficients of $Y(\cas,z)$ define a
%representation of the Virasoro algebra on the VOA $U$.
The value $c$ in (\ref{eqn:VOAs:struc:VirRel}) is called the {\em
rank} of $U$, and we denote it by $\rank(U)$.

Following \cite{HohnPhD} we say that a VOA $U$ is {\em nice} when
the eigenvalues of $L(0)$ are non-negative and contained in $\hZZ$,
and the degree zero subspace $U_0$ is spanned by the vacuum vector
$\vac$. All the VOAs we consider in this paper will be nice VOAs.

We refer the reader to \cite{FHL} %or \cite{Dun_VACo}
for a discussion of VOA modules, twisted modules, adjoint operators,
and intertwining operators. A VOA is said to be {\em self-dual} in
the case that it has no non-trivial irreducible modules other than
itself.

\subsection{Higher order brackets}\label{sec:VOAs:nbrackets}

Let $(U,Y,\vac,\cas)$ be a VOA, and let $u,v\in U$. Then taking
$F(z,w)=z^m$ in the Jacobi identity for VOAs
(\ref{eqn:VOAs:struc:Jac}) we obtain
\begin{gather}
     [u_{(m)},Y(v,w)]=\sum_{k\geq 0}
          \binom{m}{k}Y(u_{(k)}v,w)w^{m-k}
\end{gather}
and from this we may derive the following formula for the Lie
(super)bracket of two vertex operators $Y(u,z)$ and $Y(v,w)$.
\begin{gather}\label{eqn:0bracket}
%     \begin{split}
     [Y(u,z),Y(v,w)]=\sum_{m}\sum_{k\geq 0}\binom{m}{k}
          Y(u_{(k)}v,w)w^{m-k}z^{-m-1}
%     \end{split}
\end{gather}
More generally, there are a family of products on vertex operators
$Y(u,z)$ and $Y(v,w)$ defined by setting
\begin{gather}
     \left[Y(u,z)\times_n Y(v,w)\right]
          =(z-w)^n[Y(u,z),Y(v,w)]
\end{gather}
for $n$ a non-negative integer. We call these products the {\em
$n^{\text{th}}$-order brackets}, so that the $0^{\text{th}}$-order
bracket is the usual Lie bracket. Taking $F(z,w)=z^m(z-w)^n$ in the
Jacobi identity (\ref{eqn:VOAs:struc:Jac}) we obtain the following
identity for the $n^{\text{th}}$-order brackets.
\begin{gather}\label{eqn:nbracket}
     %\begin{split}
     \left[Y(u,z)\times_n Y(v,w)\right]=%&\;
     \sum_{m}\sum_{k\geq 0}\binom{m}{k}
          Y(u_{(k+n)}v,w)w^{m-k}z^{-m-1}
     %\end{split}
\end{gather}
On the other hand, the Jacobi identity for VOAs also entails that
for any $u,v,a\in U$, the expressions $Y(u,z)Y(v,w)a$ and
$Y(Y(u,z-w)v,w)a$ may be viewed as expansions of a common element of
the space
\begin{gather}
     U[[z,w]][z^{-1},w^{-1},(z-w)^{-1}]
\end{gather}
in the domains $|z|>|w|>0$ and $|w|>|z-w|>0$, respectively (c.f
\cite{FHL}). With this understanding of ``equality'' we might write
\begin{gather}\label{eqn:assoc}
     Y(u,z)Y(v,w)=\sum_{k}
          {Y(u_{(k)}v,w)}{(z-w)^{-k-1}}.
\end{gather}
Then, by (\ref{eqn:0bracket}), the data on the right hand side of
(\ref{eqn:assoc}) suffices for the purpose of computing the
commutator $[Y(u,z),Y(v,w)]$. In fact, only the terms with positive
$k$ contribute to $[Y(u,z),Y(v,w)]$, and thus it is common to write
simply
\begin{gather}\label{eqn:OPE}
     Y(u,z)Y(v,w)=\sum_{k\geq 0}
          \frac{Y(u_{(k)}v,w)}{(z-w)^{k+1}}+{\rm reg.}
\end{gather}
where ``${\rm reg.}$'' denotes regular terms in $(z-w)$. This
expression is called the {\em operator product expansion (OPE)} for
$Y(u,z)Y(v,w)$, and we have seen that the bracket $[Y(u,z),Y(v,w)]$
can be easily computed once this OPE is given. Similarly for the
$n^{\text{th}}$-order bracket, it is clear from (\ref{eqn:nbracket})
that it is sufficient to specify the expression
\begin{gather}\label{eqn:nthOPE}
     (z-w)^nY(u,z)Y(v,w)=\sum_{k\geq 0}
          \frac{Y(u_{(k+n)}v,w)}{(z-w)^{k+1}}+{\rm reg.}
\end{gather}
which we call the $n^{\text{th}}$-order OPE, in order to compute
$[Y(u,z)\times_n Y(v,w)]$.

We may define products $[u\times_nv]_{kl}$ for $u,v\in U$ and
$k,l\in\ZZ$ by setting
\begin{gather}
     [Y(u,z)\times_n Y(v,w)]=
          \sum_{k,l}\;[u\times_n v]_{kl}\;
          z^{-k-1}w^{-l-1}
\end{gather}
Then one may check directly from the definition that the component
operators $[u\times_nv]_{kl}$ may be expressed as follows in terms
of the usual Lie bracket.
\begin{gather}
     [u\times_n v]_{kl}=\sum_{m\geq 0}(-1)^m\binom{n}{m}
          [u_{(k+n-m)},v_{(l+m)}]
\end{gather}

\subsection{Vertex Lie algebras}\label{sec:VOAs:VLAs}

From \S\ref{sec:VOAs:nbrackets} we see that much of the algebra
structure on a VOA (and similarly for vertex algebras) is encoded in
the singular terms of the operators $Y(u,z)v$. The notion of vertex
Lie algebra is an axiomatic formulation of the kind of object one
obtains by replacing $Y(u,z)v$ with $\Sing Y(u,z)v$ in the
definition of vertex algebra (recall \S\ref{sec:VOAs:FmlCalc}). More
precisely, for a structure of vertex Lie algebra on a superspace
$R=R_{\bar{0}}\oplus R_{\bar{1}}$ we require morphisms $Y_-:R\otimes
R\to z^{-1}R[z^{-1}]$ and $T:R\to R$ such that the following axioms
are satisfied.
\begin{itemize}
\item{\em Translation:} $Y_-(Tu,z)=D_zY(u,z)$.
\item{\em Skew--symmetry:} $Y_-(u,z)v=\Sing e^{zT}Y_-(v,-z)u$.
\item{\em Jacobi identity:} for $\ZZ/2$ homogeneous $u,v\in R$,
and for any $m,n,l\in\ZZ$ we have
\begin{gather}\label{eqn:VOAs:VLAs:Jac}
     \begin{split}
     &\Sing\Res_{z}Y_-(u,z)Y_-(v,w)\iota_{z,w}F(z,w)\\
     &-\Sing\Res_{z}(-1)^{|u||v|}Y_-(v,w)Y_-(u,z)\iota_{w,z}F(z,w)\\
     &=\Sing\Res_{z-w}Y_-(Y_-(u,z-w)v,w)\iota_{w,z-w}F(z,w)
     \end{split}
\end{gather}
where $F(z,w)=z^mw^n(z-w)^l$.
\end{itemize}
We denote such an object by $R=(R,Y_-,T)$. We write the image of
$u\otimes v$ under $Y_-$ as $Y_-(u,z)v=\sum_{n\geq
0}u_{(n)}vz^{-n-1}$. We then obtain an equivalent definition by
replacing the Jacobi identity with the following commutativity
requirement (c.f. \cite{Pri_VAsGenByLieAlgs},\cite{Li_VAsVPAs}).
\begin{gather}
     [u_{(m)},Y_-(v,z)]=\sum_{n\geq 0}
          \binom{m}{n}\Sing Y_-(u_{(n)}v,w)w^{m-n}
\end{gather}
Vertex Lie algebras were introduced independently by V. Kac
\cite{Kac_VAsBegin} (where they are called {\em conformal algebras})
and M. Primc \cite{Pri_VAsGenByLieAlgs}. They have been studied
extensively by V. Kac and his collaborators, and we refer the reader
to \cite{Kac_FrmlDistAlgsCnfmlAlgs}, \cite{BakKacVor_CohCnflAlgs},
\cite{BakDanKac_FnPsAlgs}, and the references therein for many
interesting results.

As hinted at above, for any vertex algebra $U=(U,Y,\vac)$ we obtain
a vertex Lie algebra $(U,Y_-,T)$ by setting $Y_-=\Sing Y$, and by
setting $T$ to be the {\em translation operator:} the morphism on
$U$ defined so that $Tu=u_{-2}\vac$ for $u\in U$. (If
$(U,Y,\vac,\cas)$ is a VOA, then $T$ so defined satisfies
$T=L(-1)$.) We abuse notation somewhat to write $\Sing U$ for this
object $(U,\Sing Y,T)$. Conversely, to any vertex Lie algebra one
may canonically associate a vertex algebra called the {\em
enveloping vertex algebra}, and this construction plays an analogous
role for vertex Lie algebras as universal enveloping algebras do for
ordinary Lie algebras.

On the other hand, to each vertex Lie algebra $R$ is canonically
associated a Lie algebra $\Lie(R)$ called the {\em local Lie
algebra} of $R$. As a vector space, we have
$\Lie(R)=R[t,t^{-1}]/\Img \partial$ where $\partial$ is the operator
$T\otimes 1+\Id_R\otimes D_t$ on $R[t,t^{-1}]=
R\otimes\FF[t,t^{-1}]$. Writing $u_{[m]}$ for the image of $u\otimes
t^m$ in $\Lie(R)$ we have
\begin{prop}\label{prop:VOAs:VLAs:loclielag}
$\Lie(R)$ is a Lie algebra under the Lie bracket
\begin{gather}
     \left[u_{[m]},v_{[n]}\right]=\sum_{k\geq 0}\binom{m}{k}
          (u_{(k)}v)_{[m+n-k]}
\end{gather}
and $(Tu)_{[m]}=-mu_{[m-1]}$ for all $u$ in $R$.
\end{prop}
\begin{eg}
In the case that $R=\Sing U$ for $U$ a vertex algebra, Proposition~
\ref{prop:VOAs:VLAs:loclielag} furnishes a Lie algebra structure on
the abstract space spanned by Fourier coefficients $u_{(m)}$ for
$u\in U$, $m\in\ZZ$, and subject to the relations
$(Tu)_{(m)}=-mu_{(m-1)}$. One can define a kind of enveloping
algebra for the vertex algebra $U$ by taking a certain quotient of
(a suitable completion of) the universal enveloping algebra of
$\Lie(R)$.
\end{eg}
Given a vertex Lie algebra $R$ and a subset $\cgset\subset R$ we may
consider the {\em vertex Lie subalgebra generated by $\cgset$}. This
is by definition just the intersection of all vertex Lie subalgebras
of $R$ that contain $\cgset$. When $\cgset\subset U$ for some vertex
algebra $(U,Y,\vac)$, we will write $[\cgset]$ for the vertex Lie
subalgebra of $\Sing U$ generated by $\cgset$.
\begin{eg}
Suppose $(U,Y,\vac,\cas)$ is a VOA. Set $\cgset=\{\cas\}$ and let
$R=[\cgset]$ be the vertex Lie subalgebra of $\Sing U$ generated by
$\cgset$. Then we have
\begin{gather}
     R=\Span\{\vac,L(-1)^k\cas\mid k\geq 0\}
\end{gather}
and $\Lie(R)$ is a copy of the Virasoro algebra.
\end{eg}
The following result is straightforward.
\begin{prop}\label{prop:eVOAs:locLieRepOnVA}
If $U=(U,Y,\vac)$ is a vertex algebra and $R$ is a vertex Lie
subalgebra of $\Sing U$, then the map $\Lie(R)\to\End(U)$ given by
$u_{[m]}\mapsto u_{(m)}$ is a (well-defined) homomorphism of Lie
algebras. In particular, $\Lie(R)$ is canonically represented on
$U$.
\end{prop}

%%%%%%%%%%%%%%%%%%%%%%%%%%%%%%%%%%%%%%%%%%%%%%%%%%%%%%%%%%%%%%%%%%%%%

\section{Enhanced vertex algebras}\label{sec:eVOAs}

We begin with the
\begin{defn}
An {\em enhanced vertex algebra} is a quadruple $(U,Y,\vac,R)$ such
that $(U,Y,\vac)$ is a vertex algebra, and $R$ is a vertex Lie
subalgebra of $\Sing(U,Y,\vac)$. We say that $(U,Y,\vac,R)$ is an
{\em enhanced vertex operator algebra (enhanced VOA)} if there is a
unique $\cas\in R$ such that
\begin{enumerate}
\item     $(U,Y,\vac,\cas)$ is a VOA, and
\item     $\cas_{(n)}u=0$ for all $n\geq 2$ whenever $u\in R$ and
$\cas_{(1)}u=u$.
\end{enumerate}
When $(U,Y,\vac,R)$ is an enhanced VOA, we call $(U,Y,\vac,\cas)$
the {\em underlying VOA}. The element $\cas$ is called the {\em
Virasoro element} of $(U,Y,\vac,R)$.
\end{defn}
In many instances, the vertex Lie algebra $R$ will be of the form
$R=[\cgset]$ (recall \S\ref{sec:VOAs:VLAs}) for some (finite) subset
$\cgset\subset U$, and in such a case we may write
$(U,Y,\vac,\cgset)$ in place of $(U,Y,\vac,[\cgset])$ since there is
no loss of information. We then regard $(U,Y,\vac,\cgset)$ and
$(U,Y,\vac,\cgset')$ as {\em identical} enhanced vertex algebras
just when $[\cgset]=[\cgset']$; that is, when $\cgset$ and $\cgset'$
generate the same vertex Lie subalgebra of $\Sing(U,Y,\vac)$. Also,
we will write $(U,\cgset)$ or even $U$ in place of
$(U,Y,\vac,\cgset)$ when no confusion will arise. When
$U=(U,Y,\vac,R)$ is an enhanced vertex algebra we say that $R$
determines the {\em conformal structure} on $U$. (Then there is a
convenient coincidence with the terminology of \cite{Kac_VAsBegin},
where the objects we refer to as vertex Lie algebras are called {\em
conformal algebras}.) The automorphism group of an enhanced vertex
algebra $(U,\cgset)$ is the subgroup of the group of vertex algebra
automorphisms of $U$ that fixes each element of $\cgset$. In
practice, we may write $\Aut(U,\cgset)$ in order to emphasize this.
A pair $(M,Y^M)$ is a module over an enhanced VOA
$(U,Y,\vac,\cgset)$ just when it is a module over
$(U,Y,\vac,\omega)$. Other notions associated to VOAs (such as rank,
rationality, simplicity, \&c.) carry over directly to enhanced VOAs
in a similar way, via the underlying VOA. A morphism of enhanced
vertex algebras $(U_1,\cgset_1)\to (U_2,\cgset_2)$ is a morphism of
the underlying VOAs $(U_1,\cas_1)\to (U_2,\cas_2)$ that restricts to
a morphism of vertex Lie algebras $[\cgset_1]\to [\cgset_2]$.

For $U=(U,Y,\vac,\cgset)$ an enhanced vertex algebra, the set
$\cgset$ is called a {\em conformal generating set} for $U$, and the
elements of $\cgset$ are called {\em conformal generators}. Given a
conformal generating set $\cgset$ for an enhanced vertex algebra
$U$, we say that $\cgset$ has {\em defect $d$} if $d$ is the minimal
non-negative integer such that
\begin{gather}
     \cas^1_{(k)}\cas^2\in\Span\{\vac,T^m\nu
          \mid m\geq 0,\;\nu\in\cgset\}
\end{gather}
for all $k\geq d$ and all $\cas^1,\cas^2\in\Omega$. (The operator
$T$ here is the translation operator, defined as in
\S\ref{sec:VOAs:VLAs}.) If $\cgset$ has defect $d$ then the subspace
of $\End(U)$ spanned by the $\nu_{(n)}$ for
$\nu\in\Omega\cup\{\vac\}$ and $n\in\ZZ$
\begin{gather}
     {\rm Span}
          \{{\rm Id}_U,\nu_{(n)}\mid \nu\in\cgset, n\in\ZZ\}
          \subset \End(U)
\end{gather}
is closed under the $d^{\text{th}}$-order bracket
$[\,\cdot\,\times_d\cdot\,]$ (see \S\ref{sec:VOAs:nbrackets}). In
particular, the case that $\cgset$ has defect $0$ is just the case
that the vertex Lie algebra $[\cgset]$ generated by $\cgset$
coincides with the $\FF[T]$ module generated by
$\cgset\cup\{\vac\}$.
\begin{eg} Any VOA $(U,Y,\vac,\cas)$ furnishes an enhanced VOA (with
conformal generating set of defect $0$) when we set
$\cgset=\{\cas\}$.
\end{eg}
\begin{defn} Let $(U,\cgset)=(U,Y,\vac,\cgset)$ be an enhanced
vertex algebra and set $\cfalg(\cgset)=\Lie([\cgset])$, so that
$\cfalg(\cgset)$ is the local Lie algebra of the vertex Lie
subalgebra of $\Sing(U,Y,\vac)$ generated by $\cgset$. We call
$\cfalg(\cgset)$ the {\em local algebra} associated to $(U,\cgset)$.
If $(U,\cgset)$ is an enhanced VOA, then we regard $\cfalg(\cgset)$
as a $\QQ$-graded Lie algebra.
\end{defn}
For $\mc{L}$ a Lie algebra with $\QQ$-grading, we say that an
enhanced VOA $(U,\cgset)$ admits a {\em conformal structure of type
$\mc{L}$}, or simply an {\em $\mc{L}$--structure}, if there is a
vertex Lie subalgebra $R<[\cgset]$ such that $\Lie(R)$ is isomorphic
to $\mc{L}$, as graded Lie algebras. Usually we are interested in
the case that the choice of vertex Lie subalgebra $R$ here is
unique, in which case we say the enhanced VOA $(U,\cgset)$ admits a
{\em unique $\mc{L}$--structure}.
\begin{eg}
In \cite{Dun_VACo} we studied an $N=1$ VOA whose full automorphism
group is Conway's largest sporadic group. In the present setting, an
$N=1$ VOA is an enhanced VOA which admits a conformal generating set
(of defect $0$) consisting of two elements $\cgset=\{\cas,\scas\}$
such that the associated local algebra $\cfalg(\cgset)$ is a copy of
the $N=1$ Virasoro superalgebra. (Indeed, the vector space
underlying an $N=1$ VOA must be a superspace with non-trivial odd
part.)
\end{eg}
Another case of $\mc{L}$--structure that we would like to single out
is the case that $\mc{L}$ is the (purely even) Lie algebra spanned
by symbols $J_m$, $L_m$, and ${\bf{c}}$, for $m\in\ZZ$, and subject
to the following relations
\begin{gather}\label{eqn:eVOAs:UVirRels}
\begin{split}
    [L_m,L_n]&=(m-n)L_{m+n}+\frac{m^3-m}{12}\delta_{m+n,0}{\bf{c}},
     \quad\left[L_m,{\bf{c}}\right]=0,\\
    [L_m,J_n]&=-nJ_{m+n},\quad
    [J_m,J_n]=-m\delta_{m+n,0}{\bf c},\quad
    \left[J_m,{\bf{c}}\right]=0.
\end{split}
\end{gather}
We call this algebra the {\em $\U(1)$ Virasoro algebra}. The grading
is such that $J_m$, and $L_m$ have degree $-m$, and ${\bf{c}}$ has
degree $0$. The subalgebra generated by the $L_m$ is a copy of the
Virasoro algebra, and the subalgebra generated by the $J_m$ is a
copy of the rank $1$ Heisenberg algebra. A representation of the
$\U(1)$ Virasoro algebra is said to have {\em rank $c$} if the
central element ${\bf{c}}$ acts as multiplication by $c$, for some
$c\in\FF$.
\begin{defn}
We say that an enhanced VOA $(U,Y,\vac,\cgset)$ is an {\em enhanced
$\U(1)$--VOA} if there is a unique (up to sign) $\cej\in [\cgset]$
such that the Fourier components of the operators $Y(\cas,z)=\sum
L(m)z^{-m-2}$ and $Y(\cej,z)=\sum J(m)z^{-m-1}$ furnish a
representation of the $\U(1)$ Virasoro algebra
(\ref{eqn:eVOAs:UVirRels}) under the assignment $L_m\mapsto L(m)$,
$J_m\mapsto J(m)$.
\end{defn}
Note that for an enhanced $\U(1)$--VOA the element
$\cas_{\alpha}=\cas+\alpha T\cej$ may render
$(U,Y,\vac,\cas_{\alpha})$ a VOA for many choices of $\alpha\in\FF$.
On the other hand, there is only one choice ($\alpha=0$) for which
$(\cas_{\alpha})_{(2)}\cej=0$. If $(U,Y,\vac,\cgset)$ is an enhanced
$\U(1)$--VOA and $[\cgset]=[\cas,\cej]$ then we may call
$(U,Y,\vac,\cgset)$ simply a {\em $\U(1)$--VOA}. In
\S\ref{sec:series} we will see that $\U(1)$--VOAs (and enhanced
$\U(1)$--VOAs) admit a richer character theory than do ordinary
VOAs.

%%%%%%%%%%%%%%%%%%%%%%%%%%%%%%%%%%%%%%%%%%%%%%%%%%%%%%%%%%%%%%%%%

\section{Clifford algebras}\label{sec:cliffalgs}

The construction of VOAs that we will use arises from certain
infinite dimensional Clifford algebra modules. In this section we
recall some basic properties of Clifford algebras, and also a
construction of modules over finite dimensional Clifford algebras
using even binary linear codes, which is just a slight
generalization of the method for doubly even codes used in
\cite{Dun_VACo}. We discuss briefly the group $\Sp_{2N}$ in
\S\ref{sec:cliffalgs:spin}, and in \S\ref{sec:cliffalgs:VOAs} we
recall the construction of VOA module structure on modules over
certain infinite dimensional Clifford algebras. In
\S\ref{sec:cliffalgs:Herm} we review the Hermitian structure that
arises naturally on these objects given the existence of a suitable
Hermitian form.

\subsection{Clifford algebra structure}\label{sec:cliffalgs:struc}

Recall that $\FF$ denotes either $\RR$ or $\CC$. We suppose that
$\gt{u}$ is an $\FF$-vector space of even dimension with
non-degenerate symmetric bilinear form $\langle
\cdot\,,\cdot\rangle$, and in the case that $\FF=\RR$ we will assume
this form to be positive definite.

We write $\Cl(\gt{u})$ for the {\em Clifford algebra over $\FF$
generated by $\gt{u}$}. More precisely, we set
$\Cl(\gt{u})=T(\gt{u})/I(\gt{u})$ where $T(\gt{u})$ is the tensor
algebra of $\gt{u}$ over $\FF$ with unit denoted ${\bf 1}$, and
$I(\gt{u})$ is the ideal of $T(\gt{u})$ generated by all expressions
of the form $u\otimes u+\langle u,u\rangle{\bf 1}$ for $u\in
\gt{u}$. The natural algebra structure on $T(\gt{u})$ induces an
associative algebra structure on $\Cl(\gt{u})$. The vector space
$\gt{u}$ embeds in $\Cl(\gt{u})$, and when it is convenient we
identify $\gt{u}$ with its image in $\Cl(\gt{u})$. We also write $a$
in place of $a{\bf 1}+I(\gt{u})\in\Cl(\gt{u})$ for $a\in\FF$ when no
confusion will arise. For $u\in \gt{u}$ we have the relation
$u^2=-|u|^2$ in $\Cl(\gt{u})$. Polarization of this identity yields
$uv+vu=-2\langle u,v\rangle$ for $u,v\in \gt{u}$.

The linear transformation on $\gt{u}$ which is $-1$ times the
identity map lifts naturally to $T(\gt{u})$ and preserves
$I(\gt{u})$, and hence induces an involution on $\Cl(\gt{u})$ which
we denote by $\ogi$. The map $\ogi$ is often referred to as the {\em
parity involution}. We have $\ogi(u_1\cdots u_k)=(-1)^ku_1\cdots
u_k$ for $u_1\cdots u_k\in\Cl(\gt{u})$ with $u_i\in \gt{u}$, and we
write $\Cl(\gt{u})=\Cl(\gt{u})^0\oplus \Cl(\gt{u})^1$ for the
decomposition into eigenspaces for $\ogi$. Define a bilinear form on
$\Cl(\gt{u})$, denoted $\langle\cdot\,,\cdot\rangle$, by setting
$\langle{\bf 1},{\bf 1}\rangle=1$, and requiring that for $u\in
\gt{u}$, the adjoint of left multiplication by $u$ is left
multiplication by $-u$. Then the restriction of
$\langle\cdot\,,\cdot\rangle$ to $\gt{u}$ agrees with the original
form on $\gt{u}$. The {\em main anti-automorphism} of $\Cl(\gt{u})$
is the map we denote $\alpha$, which acts by sending $u_1\cdots u_k$
to $u_k\cdots u_1$ for $u_i\in\gt{u}$.

Suppose that $\mc{E}=\{e_i\}_{i\in{\fset}}$ is an orthonormal basis
for $\gt{u}$, indexed by a finite set ${\fset}$, and suppose that
${\fset}$ is equipped with some ordering. For
$I=\{i_1,\ldots,i_k\}\subset\fset$ we write $e_I$ for the element
$e_{i_1}e_{i_2}\cdots e_{i_k}$ in $\Cl(\gt{u})$ just when
$i_1<\cdots <i_k$. In this way we obtain an element $e_{I}$ in
$\Cl(\gt{u})$ for any $I\subset{\fset}$. (We set $e_{\emptyset}={\bf
1}$.) This correspondence depends on the choice of ordering, but our
discussion will be invariant with respect to this choice. Note that
$e_{I}e_{J}=\pm e_{I+J}$ for any $I,J\subset{\fset}$, and the set
$\{e_I\mid I\subset{\fset}\}$ furnishes an orthonormal basis for
$\Cl(\gt{u})$.

\subsection{Spin groups}\label{sec:cliffalgs:spin}

Let us write $\Cl(\gt{u})^{\times}$ for the group of invertible
elements in $\Cl(\gt{u})$. For $x\in\Cl(\gt{u})^{\times}$ and
$a\in \Cl(\gt{u})$, we set $x(a)=xax^{-1}$. %Note that $u\in
%\gt{u}$ is invertible in $\Cl(\gt{u})$ with inverse $-u/|u|^2$ so
%long as $u\neq 0$.
We will define the Pinor and Spinor groups associated to $\gt{u}$
slightly differently according as $\gt{u}$ is real or complex: in
the case that $\gt{u}$ is real, we define the {\em Pinor group}
$\textsl{Pin}(\gt{u})$ to be the subgroup of $\Cl(\gt{u})^{\times}$
comprised of elements $x$ such that $x(\gt{u})\subset \gt{u}$ and
$\alpha(x)x=\pm 1$; in the case that $\gt{u}$ is complex we define
$\textsl{Pin}(\gt{u})$ to be the set of $x\in\Cl(\gt{u})^{\times}$
such that $x(\gt{u})\subset\gt{u}$ and $\alpha(x)x=1$. In both cases
we define the {\em Spinor group} by setting
$\Sp(\gt{u})=\textsl{Pin}(\gt{u})\cap\Cl(\gt{u})^0$.

Let $x\in \textsl{Pin}(\gt{u})$. Then we have $\langle
x(u),x(v)\rangle=\langle u,v\rangle$ for $u,v\in \gt{u}$, and thus
the map $x\mapsto x(\cdot)$, which has kernel $\pm {\bf 1}$,
realizes the Pinor group as a double cover of ${O}(\gt{u})$. (If
$u\in\gt{u}$ and $\langle u,u\rangle=1$, then $u(\cdot)$ is the
orthogonal transformation of $\gt{u}$ which is minus the reflection
in the hyperplane orthogonal to $u$.) The image of $\Sp(\gt{u})$
under the map $x\mapsto x(\cdot)$ is just the group $SO(\gt{u})$.

In the case that $\gt{u}$ is real with definite bilinear form, we
have $\alpha(x)x=1$ for all $x\in\Sp(\gt{u})$, and the group
$\Sp(\gt{u})$ is generated (as a subgroup of $\Cl(\gt{u})^{\times}$)
by the expressions $\exp(\lambda e_i e_j)$ for $\lambda\in\RR$ and
$\{e_i\}$ an orthonormal basis of $\gt{u}$. The Spinor group of the
complexified space $_{\CC}\gt{u}$ is then generated by the
$\exp(\lambda e_ie_j)$ for $\lambda\in\CC$.

\subsection{Clifford algebra modules}\label{sec:cliffalgs:mods}

We obtain an $\FF_2^{{\fset}}$-grading on $\Cl(\gt{u})$ by decreeing
that for $I\subset{\fset}$, the homogeneous subspace of
$\Cl(\gt{u})$ with degree $I$ is just the $\FF$-span of the vector
$e_{I}$.
\begin{gather}
    \Cl(\gt{u})=\bigoplus_{I\subset{\fset}}
        \Cl(\gt{u})^I,\quad
    \Cl(\gt{u})^I=\FF e_I.
    %,\\
    %\Cl(\gt{u})^S\Cl(\gt{u})^R\subset\Cl(\gt{u})^{S+R}.
\end{gather}
Since this grading depends on the choice of orthonormal basis
$\mc{E}$, we will refer to it as the $\FF_2^{\mc{E}}$-grading, and
we refer to the homogeneous elements $a e_I$ for $a\in\FF$ as
$\FF_2^{\mc{E}}$-homogeneous elements. A given subset of
$\Cl(\gt{u})$ is called $\FF_2^{\mc{E}}$-homogeneous if all of its
elements are $\FF_2^{\mc{E}}$-homogeneous.

Suppose that $E$ is an $\FF_2^{\mc{E}}$-homogeneous subgroup of
$\Sp(\gt{u})$ such that the natural map $E\to \FF_2^{\fset}$ is
injective. Then $E$ is a union of elements of the form $\pm e_C$ or
$\pm\ii e_C$ for $C\subset{\fset}$, and $-{\bf 1}\notin E$, and the
image of $E$ in $\FF_2^{\fset}$ is a binary linear code on
${\fset}$. For $E$ such a subgroup of $\Sp(\gt{u})$, we write
$\mc{C}(E)$ for the associated code, and we call $E$ an
$\FF_2^{\mc{E}}$-homogeneous lift of $\mc{C}(E)$.

Suppose now that $\mc{C}(E)$ is a self-dual even code. (In the case
that $\FF=\RR$ this forces $\mc{C}(E)$ to be a doubly even code, and
this in turn forces ${\rm dim}(\gt{u})$ to be a multiple of eight.)
Then we write $\Cm(\gt{u})_E$ for the $\Cl(\gt{u})$-module defined
by $\Cm(\gt{u})_{E}=\Cl(\gt{u})\otimes_{\FF E}\FF_{1}$ where $\FF_1$
denotes a trivial $E$-module. Let us set $1_E={\bf 1}\otimes 1\in
\Cm(\gt{u})_E$. Then $\Cm(\gt{u})_E$ admits a bilinear form defined
so that $\langle 1_{E},1_E\rangle=1$, and the adjoint to left
multiplication by $u\in\gt{u}\hookrightarrow \Cl(\gt{u})$ is left
multiplication by $-u$. %In the case that $\FF=\CC$ and
%$\gt{u}=\gt{a}\oplus \gt{a}^*$ as in \S\ref{sec:cliffalgs:unitary}
%we may define a Hermitian form $(\cdot\,,\cdot)$ on $\Cm(\gt{u})_E$
%by setting $(1_E,1_E)=1$ and by decreeing that the adjoint to left
%multiplication by $\alpha a_i$ for $\alpha\in\CC$ is left
%multiplication by $-\bar{\alpha}a_i^*$. In each case the form is
%non-degenerate.
\begin{prop}
The $\Cl(\gt{u})$-module $\Cm(\gt{u})_E$ is irreducible, and a
vector-space basis for $\Cm(\gt{u})_E$ is naturally indexed by the
elements of the co-code $\mc{C}(E)^{*}$.
\end{prop}
\begin{proof}
We have $e_{S+C}{1}=\pm e_{S}{1}$ for any $S\subset{\fset}$ when
$C\in\mc{C}(E)$. This shows that a basis for $\Cm(\gt{u})_E$ is
indexed by the elements of the co-code
$\mc{C}(E)^{*}=\FF_2^{\fset}/\mc{C}(E)$, and it follows that the
irreducible submodules of $\Cm(\gt{u})_E$ are indexed by the cosets
of $\mc{C}(E)$ in its dual code $\mc{C}(E)^{\circ}$ (see
\S\ref{sec:intro:notation}). Since $\mc{C}(E)$ is assumed to be
self-dual, $\Cm(\gt{u})_E$ is irreducible.
\end{proof}
Note that the vector $1_E\in\Cm(\gt{u})_E$ is such that $g1_E=1_E$
for all $g\in E$, and $\Cm(\gt{u})_E$ is spanned by the $a1_E$ for
$a\in \Cl(\gt{u})$.

\subsection{Clifford module VOAs}\label{sec:cliffalgs:VOAs}

%Let $\gt{u}$ be again as in \S\ref{sec:cliffalgs:struc} so that
%$\gt{u}$ is a vector space of even dimension over $\FF$ with a
%non-degenerate bilinear form $(\cdot\,,\cdot)$, and
%$\mc{E}=\{e_i\}_{i\in{\fset}}$ is an orthonormal basis for $\gt{u}$
%indexed by a finite set ${\fset}$.
In this section we review the construction of VOA structure on
modules over a certain infinite dimensional Clifford algebra
associated to $\gt{u}$. The construction is quite standard and one
may refer to \cite{FFR} for example, for all the details. Our setup
is somewhat different from that in \cite{FFR} in that we also wish
to be able to handle the case that $\FF=\RR$, and we must therefore
use an alternative construction of the canonically twisted VOA
module since a polarization of $\gt{u}$ does not exist in this case.
On the other hand, all one requires is an irreducible module over
the (finite dimensional) Clifford algebra $\Cl(\gt{u})$, and the
arguments of \cite{FFR} then go through with only cosmetic changes.

\medskip

Let $\hat{\gt{u}}$ and $\hat{\gt{u}}_{\ogi}$ denote the infinite
dimensional inner product spaces described as follows.
\begin{gather}
    \hat{\gt{u}}=\coprod_{m\in\ZZ}\gt{u}\otimes t^{m+1/2},\quad
    \hat{\gt{u}}_{\ogi}=\coprod_{m\in\ZZ}\gt{u}\otimes t^m,\\
    \langle u\otimes t^r,v\otimes t^s\rangle
        =\langle u,v\rangle \delta_{r+s,0},
    \; \text{ for $u,v\in\gt{u}$ and $r,s\in\hZZ$.}
\end{gather}
We write $u(r)$ for $u\otimes t^r$ when $u\in\gt{u}$ and $r\in\hZZ$.
We consider the Clifford algebras $\Cl(\hat{\gt{u}})$ and
$\Cl(\hat{\gt{u}}_{\ogi})$. The inclusion of $\gt{u}$ in
$\hat{\gt{u}}_{\ogi}$ given by $u\mapsto u(0)$ induces an embedding
of algebras $\Cl(\gt{u})\hookrightarrow \Cl(\hat{\gt{u}}_{\ogi})$.
For $S=(i_1,\ldots,i_k)$ an ordered subset of ${\fset}$ we write
$e_S(r)$ for the element $e_{i_1}(r)\cdots e_{i_k}(r)$, which lies
in $\Cl(\hat{\gt{u}})$ or $\Cl(\hat{\gt{u}}_{\ogi})$ according as
$r$ is in $\ZZh$ or $\ZZ$. With this notation $e_S(0)$ coincides
with the image of $e_S$ under the embedding
$\Cl(\gt{u})\hookrightarrow \Cl(\hat{\gt{u}}_{\ogi})$.

Let $\mc{C}$ be an even self-dual code on $\fset$, and let
$E<\Sp(\gt{u})$ be an $\FF_2^{\mc{E}}$-homogeneous lift of $\mc{C}$
(see \S\ref{sec:cliffalgs:mods}). Note that in the case $\FF=\RR$
this forces $\mc{C}$ to be a doubly even code, and that in turn
forces $\gt{u}$ to have dimension divisible by $8$. We write
$\mc{B}(\hat{\gt{u}})$ for the subalgebra of $\Cl(\hat{\gt{u}})$
generated by the $u(m+\tfrac{1}{2})$ for $u\in\gt{u}$ and
$m\in\ZZ_{\geq 0}$. We write $\mc{B}(\hat{\gt{u}}_{\ogi})_E$ for the
subalgebra of $\Cl(\hat{\gt{u}}_{\ogi})$ generated by
$E\subset\Cl(\gt{u})$, and the $u(m)$ for $u\in\gt{u}$ and
$m\in\ZZ_{> 0}$. Let $\FF_1$ denote a one-dimensional module for
either $\mc{B}(\hat{\gt{u}})$ or $\mc{B}(\hat{\gt{u}}_{\ogi})_E$,
spanned by a vector $1_E$, such that $u(r)1_E=0$ whenever
$r\in\hZZ_{>0}$, and such that $g(0)1_{E}=1_{E}$ for $g\in E$. We
write $A(\gt{u})$ (respectively $A(\gt{u})_{\ogi,E}$) for the
$\Cl(\hat{\gt{u}})$-module (respectively
$\Cl(\hat{\gt{u}}_{\ogi})$-module) induced from the
$\mc{B}(\hat{\gt{u}})$-module structure (respectively
$\mc{B}(\hat{\gt{u}}_{\ogi})_E$-module structure) on $\FF_{1}$.
\begin{gather}
    A(\gt{u})=\Cl(\hat{\gt{u}})
        \otimes_{\mc{B}(\hat{\gt{u}})}\FF_{1},\qquad
    A(\gt{u})_{\ogi,E}=\Cl(\hat{\gt{u}}_{\ogi})
        \otimes_{\mc{B}(\hat{\gt{u}}_{\ogi})_E}\FF_{1}.
\end{gather}
We write $\vac$ for the vector $1\otimes 1_{E}$ in $A(\gt{u})$, and
we write $\vac_{\ogi}$ or $\vac_E$ for the vector $1\otimes 1_{E}$
in $A(\gt{u})_{\ogi,E}$. %When no confusion will arise, we simply
%write $A(\gt{u})_{\ogi}$ in place of $A(\gt{u})_{\ogi,E}$. We often
%write $u_1(-r_1)\cdots u_k(-r_k)$ for $u_1(-r_1)\cdots u_k(-r_k){\bf
%1}$ when the $r_i$ are in $\ZZh_{>0}$, and similarly, we write
%$u_1(-m_1)\cdots u_k(-m_k)$ for $u_1(-m_1)\cdots u_k(-m_k){\bf
%1}_{\ogi}$ when the $m_i$ are in $\ZZ_{\geq 0}$.

\medskip

The space $A(\gt{u})$ supports a structure of VOA. In order to
define the vertex operators we require the notion of fermionic
normal ordering for elements in $\Cl(\hat{\gt{u}})$ and
$\Cl(\hat{\gt{u}}_{\ogi})$. The fermionic normal ordering on
$\Cl(\hat{\gt{u}})$ is the multi-linear operator defined so that for
$u_i\in\gt{u}$ and $r_i\in\ZZh$ we have
\begin{equation}
    :\!u_1(r_1)\cdots u_k(r_k)\!:=\,
        {\rm sgn}(\sigma)u_{\sigma 1}(r_{\sigma 1})
        \cdots u_{\sigma k}(r_{\sigma k})
\end{equation}
where $\sigma$ is any permutation of the index set $\{1,\ldots,k\}$
such that $r_{\sigma 1}\leq\cdots\leq r_{\sigma k}$. For elements in
$\Cl(\hat{\gt{u}}_{\ogi})$ the fermionic normal ordering is defined
in steps by first setting
\begin{equation}
    :\!u_1(0)\cdots u_k(0)\!:\,=\frac{1}{k!}\sum_{\sigma\in S_k}
        {\rm sgn}(\sigma)u_{\sigma 1}(0)
        \cdots u_{\sigma k}(0)
\end{equation}
for $u_i\in \gt{u}$. Then in the situation that $n_i\in\ZZ$ are such
that $n_{i}\leq n_{i+1}$ for all $i$, and there are some $s$ and $t$
(with $1\leq s\leq t\leq k$) such that $n_j=0$ for $s\leq j\leq t$,
we set
\begin{equation}
     \begin{split}
    &:\!u_1(n_1)\cdots u_k(n_k)\!:\\
    =u_{1}(n_{1})\cdots u_{ s-1}(n_{s-1})&
    \!:\!u_s(0)\cdots u_t(0)\!:\!
    u_{t+1}(n_{t+1})\cdots u_{k}(n_k)
    \end{split}
\end{equation}
Finally, for arbitrary $n_i\in\ZZ$ we set
\begin{equation}
    :\!u_1(n_1)\cdots u_k(n_k)\!:\,=\,
        {\rm sgn}(\sigma)\!:\!u_{\sigma 1}(n_{\sigma 1})
        \cdots u_{\sigma k}(n_{\sigma k})\!:
\end{equation}
where $\sigma$ is again any permutation of the index set
$\{1,\ldots,k\}$ such that $n_{\sigma 1}\leq\cdots\leq n_{\sigma
k}$, and we extend the definition multilinearly to
$\Cl(\hat{\gt{u}}_{\ogi})$.

\medskip

For $u\in\gt{u}$ we now define the generating function, denoted
$u(z)$, of operators on $A(\gt{u})_{\Theta}=A(\gt{u})\oplus
A(\gt{u})_{\ogi}$ by setting
\begin{gather}
     u(z)=\sum_{r\in\hZZ}u(r)z^{-r-1/2}
\end{gather}
Note that $u(r)$ acts as $0$ on $A(\gt{u})$ if $r\in \ZZ$, and acts
as $0$ on $A(\gt{u})_{\ogi}$ if $r\in \ZZh$. To an element $a\in
A(\gt{u})$ of the form $a=u_{1}(-m_1-\tfrac{1}{2})\cdots
u_{k}(-m_k-\tfrac{1}{2}){\bf 1}$ for $u_i\in \gt{u}$ and $m_i\in
\ZZ_{\geq 0}$, we associate the operator valued power series
$\overline{Y}(a,z)$, given by
\begin{gather}
    \overline{Y}(a,z)=\,:\!D_z^{(m_1)}u_{i_1}(z)\cdots
        D_z^{(m_k)}u_{i_k}(z)\!:
\end{gather}
We define the vertex operator correspondence
\begin{gather}
    Y(\cdot\,,z):A(\gt{u})\otimes A(\gt{u})_{\Theta}
        \to A(\gt{u})_{\Theta}((z^{1/2}))
\end{gather}
by setting $Y(a,z)b=\overline{Y}(a,z)b$ when $b\in A(\gt{u})$, and
by setting $Y(a,z)b=\overline{Y}(e^{\Delta_z}a,z)b$ when $b\in
A(\gt{u})_{\ogi}$, where $\Delta_z$ is the expression defined by
\begin{gather}
    \Delta_z=-\frac{1}{4}\sum_i\sum_{m,n\in\ZZ_{\geq 0}}C_{mn}
        e_i(m+\tfrac{1}{2})e_i(n+\tfrac{1}{2})z^{-m-n-1}\\
        C_{mn}=\frac{1}{2}\frac{(m-n)}{m+n+1}
            \binom{-\tfrac{1}{2}}{m}\binom{-\tfrac{1}{2}}{n}
\end{gather}
Set $\cas=\tfrac{1}{4}\sum_i e_i(-\tfrac{1}{2})
e_i(-\tfrac{3}{2})\vac \in A(\gt{u})_2$. Then one has the following
\begin{thm}[\cite{FFR}]
The map $Y$ defines a structure of self-dual VOA of rank $N$ on
$A(\gt{u})$ when restricted to $A(\gt{u})\otimes A(\gt{u})$, and the
Virasoro element is given by $\cas$. The map $Y$ defines a structure
of $\ogi$-twisted $A(\gt{u})$-module on $A(\gt{u})_{\ogi}$ when
restricted to $A(\gt{u})\otimes A(\gt{u})_{\ogi}$.
\end{thm}
Observe that $A(\gt{u})_2$ is spanned by vectors of the form
$e_I(-\tfrac{1}{2})\vac$ for $I\subset\fset$ with $|I|=4$, and by
the $e_i(-\tfrac{3}{2})e_j(-\tfrac{1}{2})\vac$ with $i,j\in\fset$.

Essentially all we need to know about the expressions $Y(a,z)b$ for
$b\in A(\gt{u})_{\ogi}$ is contained in the following
\begin{prop}\label{prop:cliffalgs:VOAs:twops}
Let $b\in A(\gt{u})_{\ogi}$.
\begin{enumerate}
\item   If $a=\vac\in A(\gt{u})_0$ then $Y(a,z)b=b$.
\item   If $a\in A(\gt{u})_1$ then $\Delta_za=0$ so that
$Y(a,z)b=\overline{Y}(a,z)b$.
\item   If $a\in A(\gt{u})_2$ and
%\begin{gather}
     $a\in\Span\left\{e_I(-\tfrac{1}{2})\vac,
          e_i(-\tfrac{3}{2})e_j(-\tfrac{1}{2})\vac\mid i\neq j\right\}$
%\end{gather}
then $\Delta_za=0$ and $Y(a,z)b= \overline{Y}(a,z)b$.
\item   For $a=e_i(-\tfrac{1}{2})e_i(-\tfrac{3}{2})\vac$ we have
$\Delta_za=\tfrac{1}{4}z^{-2}$ and $\Delta_z^2a=0$ so that
$Y(a,z)b=\overline{Y}(a,z)b+\tfrac{1}{4}bz^{-2}$ in this case.
\end{enumerate}
\end{prop}
As a corollary of Proposition~\ref{prop:cliffalgs:VOAs:twops} we
have that
\begin{gather}
     Y(\omega,z){\bf 1}_{\ogi}=\frac{N}{8}
          {\bf 1}_{\ogi} z^{-2}
\end{gather}
and consequently the $L(0)$-grading on $A(\gt{u})_{\Theta}$ is given
by
\begin{gather}
    A(\gt{u})=\coprod_{n\in\hZZ_{\geq 0}}A(\gt{u})_n,\quad
    A(\gt{u})_{\ogi}=\coprod_{n\in\hZZ_{\geq 0}}
        (A(\gt{u})_{\ogi})_{n+N/8}.
\end{gather}
Given a specific choice of $E$, the embedding of $\Cl(\gt{u})$ in
$\Cl(\hat{\gt{u}}_{\ogi})$ gives rise to an isomorphism of
$\Cm(\gt{u})_E$ with $(A(\gt{u})_{\ogi,E})_{N/8}$, and it will be
convenient to consider these spaces as identified.

\medskip

The group $\Sp(\gt{u})$ acts naturally on $A(\gt{u})_{\Theta}$, and
this action is generated by the exponentials of the operators
$x_{(0)}$ for $x\in A(\gt{u})_1$. In particular, any
$a\in\Sp(\gt{u})\subset\Cl(\gt{u})$ may be regarded as a VOA
automorphism of $A(\gt{u})$, and as an equivariant linear
isomorphism of the $A(\gt{u})$-module $A(\gt{u})_{\ogi,X}$. Recall
the element $\ogz=e_{\fset}$ of \S\ref{sec:cliffalgs:spin}. When
constructing a realization of the twisted module $A(\gt{u})_{\ogi}$
we will always choose $E$ so that $e_{\fset}\in E$, and thus $\ogz$
will be the unique preimage in $\Sp(\gt{u})$ of $-\Id_{\gt{u}}$ that
fixes the vector $\vac_{E}$ in $A(\gt{u})_{\ogi,E}$. Equivalently,
$\ogz$ is the unique element of $\Sp(\gt{u})$ whose action on
$A(\gt{u})_{\Theta}$ coincides with that of the parity involution
$\ogi$.

\subsection{Hermitian structure}\label{sec:cliffalgs:Herm}

We will now take special interest in the case that $\FF=\CC$ and
$\gt{u}$ is of the form $\gt{u}=\gt{a}\oplus\gt{a}^*$ for $\gt{a}$
some complex vector space with non-degenerate Hermitian form, and
$\gt{a}^*$ the dual space to $\gt{a}$. We shall denote the Hermitian
form by $(\cdot\,,\cdot)$, and our convention will be that a
Hermitian form be anti-linear in the second slot. In this instance
we insist that the bilinear form $\lab\cdot\,,\cdot\rab$ on $\gt{u}$
be $1/2$ the symmetric $\CC$-bilinear form on $\gt{u}$ induced by
the natural pairing between $\gt{a}$ and $\gt{a}^*$. That is, we set
$\lab a,f\rab =\lab f,a\rab =\tfrac{1}{2}f(a)$ for $a\in\gt{a}$ and
$f\in\gt{a}^*$. We suppose that $\hfset$ is some set with
cardinality $N$, and that $\{a_i\}_{i\in\hfset}$ is a basis for
$\gt{a}$ satisfying $(a_i,a_j)=\delta_{ij}$ for $i,j\in\hfset$. We
suppose also that $\{a_i^*\}_{i\in\hfset}$ is the dual basis for
$\gt{a}^*$, so that $\lab a_i,a_j^*\rab= \lab a_j^*,a_i\rab=
\tfrac{1}{2}\delta_{ij}$. With this convention we have
$a_ia_j^*+a_j^*a_i+\delta_{ij}=0$ in $\Cl(\gt{u})$.

The Hermitian form on $\gt{a}$ entails an anti-involution $\ainv$
defined on $\gt{u}$ in such a way that $\ainv$ interchanges $\alpha
a_i$ with $\bar{\alpha}a_i^*$ for $\alpha\in\CC$ and $i\in\hfset$.
We then have the identity $(u,v)=2\lab u,\ainv v\rab$ for
$u,v\in\gt{u}$, and our definition of $\ainv$ is invariant with
respect to the choice of orthonormal basis for $\gt{a}$. The fixed
points under $\ainv$ form a real vector space of dimension $2N$ in
$\gt{u}$. We will write $_{\RR}\gt{u}$ for this space, keeping in
mind the dependence upon the choice of Hermitian structure for
$\gt{a}$. With such $\gt{u}$ we make the convention that $\mc{E}$
shall denote the orthonormal basis for $_{\RR}\gt{u}$ obtained in
the following way. We let $\hfset'$ be another copy of the set
$\hfset$ with the natural identification
$\hfset\leftrightarrow\hfset'$ denoted $i\leftrightarrow i'$, and we
set $\mc{E}=\{e_i,e_{i'}\}_{i\in\hfset}$ where
\begin{gather}
     e_i=a_i+a_i^*,\quad e_{i'}=\ii(a_i-a_i^*),
\end{gather}
for $i\in\hfset$. We let $\fset=\hfset\cup\hfset'$ and we agree to
order $\fset$ in this case by first choosing an ordering on
$\hfset$, and then setting $i<j'$ for all $i,j\in\hfset$.

In our present situation there is a standard choice of group to play
the role of $E$ in \S\ref{sec:cliffalgs:mods}. We reserve the
notation $X$ for the subgroup of $\Sp(\gt{u})$ generated by elements
of the form $\ii e_ie_{i'}$ for $i\in\fset$, and we may then take
$E=X$ in \S\ref{sec:cliffalgs:mods}, so as to obtain the
$\Cl(\gt{u})$-module $\Cm(\gt{u})_X$. We note here that ${\rm
CM}(\gt{u})_X$ may be characterized as the left-module over ${\rm
Cliff}(\gt{u})$ spanned by a vector ${1}_{X}$ such that
$e_{i'}{1}_X=\ii e_i{1}_X$ for all $i\in \hfset$, or equivalently,
such that $a_i^*1_X=0$ for all $i\in\hfset$. In particular,
$\Cm(\gt{u})_X$ may be naturally identified with the space
$\bigwedge(\gt{a})$. When no confusion will arise we regard
$\bigwedge(\gt{a})$ as a $\Cl(\gt{u})$ module via this
identification. Note also the identity $e_{i_1}\cdots
e_{i_k}1_X=a_{i_1}\cdots a_{i_k}1_X$ in $\Cm(\gt{u})_X$ for
$I=\{i_1,\ldots,i_k\}\subset\hfset$.

%%%%%%%%%%%%%%%%%%%%%%%%%%%%%%%%%%%%%%%%%%%%%%%%%%%%%%%%%%%%%%%%%

\section{Linear groups}\label{sec:GLnq}

In this section we define families of enhanced VOAs with symmetry
groups related to the linear groups $\GL_N(\CC)$ for $N$ a positive
integer. In defining these objects we are laying some of the ground
work for the construction of the enhanced VOA for the Rudvalis group
in \S\ref{sec:Ru}.

\subsection{Untwisted construction}\label{sec:GLnq:untw}

Let $N$ be a positive integer and let $\gt{a}$ be a complex vector
space of dimension $N$ equipped with a positive definite Hermitian
form denoted $(\cdot\,,\cdot)$. As in \S\ref{sec:cliffalgs:Herm} we
set ${\gt{u}}=\gt{a}\oplus\gt{a}^*$, we extend the Hermitian form to
$\gt{u}$, and equip this space also with the symmetric bilinear form
$\lab\cdot\,,\cdot\rab$ which is just the form arising naturally
from the pairing between $\gt{a}$ and $\gt{a}^*$ scaled by a factor
of $1/2$.

Let us set $\cgU=\{\cas,\cej\}$ where $\cas$ and $\cej$ are given by
\begin{align}
     \cas&=\frac{1}{4}\sum_{i\in\Sigma\cup\Sigma'}
     e_i(-\tfrac{1}{2})e_i(-\tfrac{3}{2})\vac\in A(\gt{u})_2,\\
     \cej&=\frac{1}{2}\sum_{i\in\Sigma}
     e_i(-\tfrac{1}{2})e_{i'}(-\tfrac{1}{2})\vac\\
     &=\ii\sum_{i\in\Sigma}a_i^*(-\tfrac{1}{2})a_i(-\tfrac{1}{2})\vac
          \in A(\gt{u})_1,
\end{align}
so that $\cas$ is just the usual Virasoro element for $A(\gt{u})$,
and $\cej$ is some element dependent upon the Hermitian structure on
$\gt{a}$. %Note that $\cas$ and $\cej$ are fixed by the action of
%$\ainv$, so that they lie even in $_{\RR}A(\gt{u})$.
We define operators $L(n)$ and $\ceJ(n)$ for $n\in\ZZ$ by setting
\begin{align}
     Y(\cas,z)&=\;L(z)=\sum L(n)z^{-n-2},\\
     Y(\cej,z)&=\;\ceJ(z)=\sum \ceJ(n)z^{-n-1}.
\end{align}
We then have
\begin{prop}\label{prop:GLnq:untw:UnOPEs}
The vertex operators $L(z)$ and $\ceJ(z)$ satisfy the following OPEs
\begin{gather}
     L(z)L(w)=%\frac{1}{4}
          \frac{N/2}{(z-w)^4}+\frac{2L(w)}{(z-w)^2}
          +\frac{D_wL(w)}{(z-w)}
          +{\rm reg.}\\
     L(z)J(w)=\frac{J(w)}{(z-w)^2}+\frac{D_wJ(w)}{(z-w)}
          +{\rm reg.}\\
     J(z)J(w)=-\frac{N}{(z-w)^2}
          +{\rm reg.}
\end{gather}
\end{prop}
We recalled in \S\ref{sec:VOAs:nbrackets} how the OPE of two vertex
operators encodes the commutation relations of their component
operators. Proposition~\ref{prop:GLnq:untw:UnOPEs} entails the
following
\begin{cor}\label{cor:GLnq:untw:LJcomrels}
The operators $L(m)$ and $J(m)$ satisfy
\begin{gather}
     \begin{split}
     \left[L(m),L(n)\right]&=(m-n)L(m+n)
          +\frac{m^3-m}{12}\delta_{m+n,0}{N}\Id\\
     [L(m),J(n)]&=-nJ(m+n),\quad
     [J(m),J(n)]%&
     =-m\delta_{m+n,0}N\Id
     \end{split}
\end{gather}
for all $m,n\in\ZZ$.
\end{cor}
The objects $(\au,\cgU)$ we have constructed provide our first
family of examples of enhanced VOAs beyond ordinary VOAs and $N=1$
VOAs. In fact they are examples of $\U(1)$--VOAs (see
\S\ref{sec:eVOAs}), as is shown by the following
\begin{prop}\label{prop:GLnq:untw:u(1)struc}
For $N$ a positive integer, the quadruple $(\au,Y,\vac,\cgU)$ is a
self-dual $\U(1)$--VOA of rank $N$.
\end{prop}
\begin{proof}
We will show that $(\au,\cgU)$ is a $\U(1)$--VOA since the other
claims have been verified already. From Proposition~
\ref{prop:GLnq:untw:UnOPEs} we see that $\cgU$ has defect $0$, so
that the vertex Lie algebra generated by $\cgU$ is just
$[\cgU]=\Span\{\vac,T^k\cej,T^k\cas\mid k\in\ZZ\}$. We have seen
that the element $\cas\in\cgU$ is such that $(\au,Y,\vac,\cas)$ is a
VOA. From the explicit spanning set given we see that the only other
Virasoro elements in $[\cgU]$ are of the form
$\cas_{\alpha}=\cas+\alpha T\cej$ for some $\alpha\in \CC$. Using
Corollary \ref{cor:GLnq:untw:LJcomrels} we check that the Fourier
coefficients of $Y(\cas_{\alpha},z)=L_{\alpha}(m)z^{-m-2}$ furnish a
representation of the Virasoro algebra of rank $N(1+12\alpha^2)$,
and we have $L_{\alpha}(n)=L(n)+\alpha(-n-1)J(n)$ where
$Y(\cas,z)=\sum L(m)z^{-m-2}$ and $Y(\cej,z)=\sum J(m)z^{-m-1}$. In
particular, $L_{\alpha}(-1)=L(-1)=T$ for all $\alpha$, and
$L_{\alpha}(0)=L(0)-\alpha J(0)$. Thus if $u\in [\cgU]$ and
$L_{\alpha}(0)u=u$ then $u\in\Span\{\cej\}$. Now we compute
$L_{\alpha}(1)\cej= \alpha(-2)J(1)\cej= 2\alpha N$, and it follows
that $(\au,Y,\vac,\cgU)$ is an enhanced VOA with the Virasoro
element given by $\cas=\cas_{0}$. From Corollary
\ref{cor:GLnq:untw:LJcomrels} we see that $\cej$ is the unique (up
to sign) choice of element in $[\cgU]\cap (\au)_1$ for which the
component operators of $Y(\cej,z)$ satisfy the required relations
(\ref{eqn:eVOAs:UVirRels}).
\end{proof}
The proof of Theorem~\ref{prop:GLnq:untw:u(1)struc} shows also that
$(A(\gt{u})_{\bar{0}},Y,\vac,\cgU)$ is a $\U(1)$--VOA.

We would like to compute the automorphism group of
$(A(\gt{u}),\cgU)$. First we will compute the automorphism groups of
the underlying VOA, and of its even subVOA.
\begin{prop}\label{prop:GLnq:untw:VOAauts}
Let $N>4$. Then the automorphism group of $A(\gt{u})$ is
$\Or(\gt{u})$, and the automorphism group of $A(\gt{u})_{\bar{0}}$
is $\Or(\gt{u})/\lab\!\pm\Id\!\rab$.
\end{prop}
\begin{proof}
Let $\gt{g}=\{x_0\mid x\in A(\gt{u})_1\}\subset \End(A(\gt{u}))$.
Then $\gt{g}$ is a simple complex Lie algebra of type $D_N$. We set
$G=\Aut(A(\gt{u}))$ and $G_0=\Aut(A(\gt{u})_{\bar{0}})$, and we
write $S$ for the subgroup of $G$ generated by the exponentials
$e^X$ for $X\in\gt{g}$. Note that any automorphism of $A(\gt{u})$
must preserve the even subspace $A(\gt{u})_{\bar{0}}$, so there is a
natural map $\phi:G\to G_0$. The kernel of this map is generated by
the canonical automorphism of $U$; viz. the automorphism that fixes
the even subspace and negates the odd subspace. In particular,
$\ker(\phi)$ is contained in $S$. Note also that there is a natural
isomorphism $G_0\cong\Aut(\gt{g})$ since $A(\gt{u})_{\bar{0}}$ is
generated by its subspace of elements of degree $1$.
\begin{gather}
     A(\gt{u})_{\bar{0}}=\Span\left\{u^{1}_{-n_1}\cdots u^k_{-n_k}\vac\mid
          u^i\in A(\gt{u})_1,\,n_1\geq \cdots\geq n_k>0\right\}
\end{gather}
Now the image of $\phi(S)$ in $\Aut(\gt{g})$ is the subgroup
$\Inn(\gt{g})$ of inner automorphisms of $\gt{g}$. At least for $N>
4$ then, $\phi(S)$ has index two in $G_0$ since $\gt{g}$ is of type
$D_N$, so $G_0=\phi(S)\cup \bar{x}\phi(S)$ for some $\bar{x}\in
G_0\setminus \phi(S)$.

Observe now that $A(\gt{u})$ is generated by its subspace of degree
$1/2$,
\begin{gather}
     A(\gt{u})=\Span\left\{u^{1}_{-n_1}\cdots u^k_{-n_k}\vac\mid
          u^i\in A(\gt{u})_{1/2},\,n_1\geq \cdots\geq n_k>0\right\}
\end{gather}
which is naturally identified with $\gt{u}$ under the correspondence
$u\leftrightarrow u(-\tfrac{1}{2})\vac$ for $u\in\gt{u}$. Further,
this correspondence identifies the bilinear form on $\gt{u}$ with
the restriction to $A(\gt{u})_{1/2}$ of the invariant form on
$A(\gt{u})$. We see then that any automorphism of $A(\gt{u})$
induces an element of $\Or(\gt{u})$ by restriction to the degree
$1/2$ subspace, and conversely, any element of $\Or(\gt{u})$ extends
to an automorphism of $A(\gt{u})$, and thus we may identify
$\Aut(A(\gt{u}))$ as $\Or(\gt{u})$. We see also, that the map
$\phi:G\to G_0$ is surjective, since $\Or(\gt{u})$ acts
non-trivially on the even subVOA $A(\gt{u})_{\bar{0}}$. The group
$S$ must then be $\SO(\gt{u})$, and we may regard the $\bar{x}$ in
$G_0=\phi(S)\cup \bar{x}\phi(S)$ as the image under $\phi$ of an
element $x$ say, in $\Or(\gt{u})\setminus \SO(\gt{u})$. Since the
kernel of $\phi$ is the canonical automorphism group of $A(\gt{u})$
(generated by $-\Id$), we see that $G_0$ is the group
$\Or(\gt{u})/\lab\!\pm\Id\!\rab$.
\end{proof}

\begin{prop}\label{prop:GLnq:untw:eVOAauts}
The automorphism group of $(A(\gt{u}),\cgset_U)$ is $\GL(\gt{a})$.
The automorphism group of $(A(\gt{u})_{\bar{0}},\cgset_U)$ is
$\GL(\gt{a})/\lab\!\pm\Id\!\rab$.
\end{prop}
\begin{proof}
Recall that $\cej= \sum\ii
a_i^*(-\tfrac{1}{2})a_i(-\tfrac{1}{2})\vac$. From the previous
proposition $\Aut(A(\gt{u}),\{\cas\})=\Or(\gt{u})$. The group
$\Aut(A(\gt{u}),\cgset_U)$ is just the subgroup of $\Aut(A(\gt{u}))$
that fixes $\cej$. Consider the automorphism of $A(\gt{u})$ obtained
by setting $\ogi^{1/2}=\exp(\pi J(0)/2)$. It is the automorphism of
$A(\gt{u})$ induced by the orthogonal transformation of $\gt{u}$
which is multiplication by $\ii$ on $\gt{a}$, and multiplication by
$-\ii$ on $\gt{a}^*$. Clearly, any element of
$\Aut(A(\gt{u}),\cgset_U)$ commutes with $\ogi^{1/2}$. On the other
hand, if $g\in\Or(\gt{u})$ commutes with $\ogi^{1/2}$ then $g$
preserves the decomposition $\gt{u}=\gt{a}\oplus\gt{a}^*$, and with
respect to the basis $\{a_1,\ldots,a_1^*,\ldots\}$ is represented by
a block matrix
\begin{gather}\label{blckdecomp}
     g\sim
          \left(
              \begin{array}{cc}
                T_g^t & 0 \\
                0 & T_g^{-1} \\
              \end{array}
            \right)
\end{gather}
for some invertible $N\times N$ matrix $T_g$, where $T_g^t$ denotes
the transpose of $T_g$. Evidently $\Aut(A(\gt{u}),\cgset_U)$ is the
centralizer of $\ogi^{1/2}$ in $\Or(\gt{u})$, and this group is all
matrices of the form $(\ref{blckdecomp})$. Since the natural map
$\Aut(A(\gt{u}))\to\Aut(A(\gt{u})_{\bar{0}})$ is surjective, the
natural map $\Aut(A(\gt{u}),\cgset_U)\to
\Aut(A(\gt{u})_{\bar{0}},\cgset_U)$ is also surjective. The kernel
in each case is $\pm \Id$, and this completes the proof.
\end{proof}
Note that the embedding of $\GL(\gt{a})$ in $\Or(\gt{u})$ given by
(\ref{blckdecomp}) is actually an embedding of $\GL(\gt{a})$ in
$\SO(\gt{u})$.

We record some of the results of this section in
\begin{thm}\label{thm:GLnq:untw:main}
For $N$ a positive integer, the quadruple $(\au,Y,\vac,\cgU)$ is a
self-dual $\U(1)$--VOA of rank $N$. The full automorphism group of
$(\au,\cgU)$ is $\GL(\gt{a})$.
\end{thm}

\subsection{Twisted construction}\label{sec:GLnq:tw}

In practice we may wish to take advantage of the fact that the
larger group $\Sp(\gt{u})$ acts non-trivially on $A(\gt{u})_{\ogi}$.
Also, we may wish to be able to realize groups no larger than
$SL(\gt{a})$ as symmetry groups of enhanced VOAs. For these goals it
is useful to consider the following twisted analogue of the Clifford
module construction of VOAs.

Let $\gt{a}$ and $\gt{u}$ be as in the previous section, but
restrict now to the case that $N$ is a positive integer divisible by
four. Recall from \S\ref{sec:cliffalgs:VOAs} that $A(\gt{u})$ is a
VOA with superspace decomposition given by
\begin{gather}
     A(\gt{u})=A(\gt{u})^0\oplus A(\gt{u})^1,
\end{gather}
coinciding with the decomposition into eigenspaces for the action of
the parity involution $\ogi$. Recall also that the $\ogi$-twisted
module $A(\gt{u})_{\ogi}$ may be realized as $A(\gt{u})_{\ogi,X}$
where $X$ is as in \S\ref{sec:cliffalgs:Herm}, and the parity
involution acts naturally also on $A(\gt{u})_{\ogi,X}$ with
eigenspaces $A(\gt{u})_{\ogi,X}^0$ and $A(\gt{u})_{\ogi,X}^1$. We
define the space $\atw{\gt{u}}$ by setting
\begin{gather}\label{eqn:GLnq:VOAstruc:defatw}
    \atw{\gt{u}}=A(\gt{u})^0\oplus A(\gt{u})_{\ogi,X}^0.
\end{gather}
and we claim that $\atw{\gt{u}}$ admits a structure of self-dual VOA
of rank $N$. In the case that $N$ is divisible by $8$, this object
is even a VOA, but for now we are more interested in the super case,
so we will henceforth assume that $N$ is congruent to $4$ modulo
$8$.

The superspace decomposition of $\atw{\gt{u}}$ coincides with the
decomposition in (\ref{eqn:GLnq:VOAstruc:defatw}) so that the even
subVOA is $A(\gt{u})^0$. We require to exhibit the vertex operators
on $\atw{\gt{u}}$, and this may be done as follows. The required
vertex operator correspondence
$Y:\atw{\gt{u}}\otimes\atw{\gt{u}}\to\atw{\gt{u}}((z))$ is defined
already on $\atw{\gt{u}}_{\bar{0}}\otimes\atw{\gt{u}}_{\bar{0}}$ and
on $\atw{\gt{u}}_{\bar{0}}\otimes\atw{\gt{u}}_{\bar{1}}$ curtesy of
\S\ref{sec:cliffalgs:VOAs}. For $u\otimes v\in
\atw{\gt{u}}_{\bar{1}}\otimes\atw{\gt{u}}_{\bar{0}}$ we define
$Y(u,z)v$ by setting
\begin{gather}\label{YdefAtwOnA}
  Y(u,z)v=e^{zL(-1)}Y(v,- z)u
\end{gather}
Suppose now that $u\otimes v\in\atw{\gt{u}}_{\bar{1}}\otimes
\atw{\gt{u}}_{\bar{1}}$. Then we %are motivated by
%\S\ref{sec:VOAs:adj} to
define $Y(u,z)v$ by requiring that for any $w\in
\atw{\gt{u}}_{\bar{0}}$ we should have
\begin{gather}\label{YdefAtwOnAtw}
    \langle Y(u,z)v\mid w\rangle
    =(-1)^n\langle e^{z^{-1}L(1)}v\mid
    Y(w,- z^{-1})e^{zL(1)} z^{-2L(0)}
    u\rangle
\end{gather}
whenever $u\in \atw{\gt{u}}_{n-1/2}$ for some $n\in\ZZ$.

The proof of the following proposition is almost identical to that
of Proposition~4.1 in \cite{Dun_VACo} (see also
\cite{HuaXtnMoonVOA}).
\begin{prop}\label{prop:GLnq:tw:VOAstruc}
The map $Y:\atw{\gt{u}}\otimes\atw{\gt{u}}\to\atw{\gt{u}}((z))$
defines a structure of self-dual VOA of rank $N$ on $\atw{\gt{u}}$.
\end{prop}

%The action of $\ainv$ on $A(\gt{u})\oplus A(\gt{u})_{\ogi,X}$ (see
%\S\ref{sec:cliffalgs:Herm}) preserves $\atw{\gt{u}}$. We write
%$\atwr{\gt{u}}$ for the $\ainv$-fixed points, and $\atwr{\gt{u}}$ is
%then a real form for the self-dual VOA $\atw{\gt{u}}$. As promised
%at the end of \S\ref{sec:cliffalgs:VOAs}, the element $\gt{z}$ of
%\S\ref{sec:cliffalgs:spin} lies in $X$, so that the automorphism
%group of the VOA $\atw{\gt{u}}$ is the group $\Spql{\gt{u}}$. The
%automorphism group of the real form $\atwr{\gt{u}}$ is
%$\Spql{{_{\RR}\gt{u}}}$, and the action of $\Sp({_{\RR}\gt{u}})$ on
%$\atwr{\gt{u}}$ is generated by the exponentials of operators
%$x_{(0)}$ for $x$ in $\LieSOr={_{\RR}A(\gt{u})_1}$.

Recall that the even part of $\atw{\gt{u}}$ is
$A(\gt{u})^0=A(\gt{u})_{\bar{0}}$ (so long as $N/4$ is odd).
\begin{prop}\label{prop:GLnq:tw:VOAauts}
For $N>4$ the automorphism group of the VOA $\atw{\gt{u}}$ is
$\Sp(\gt{u})/\Kumz$.
\end{prop}
\begin{proof}
The proof is very similar to that of Proposition~
\ref{prop:GLnq:untw:VOAauts}. We set $G=\Aut(\atw{\gt{u}})$ and
$G_0=\Aut(\atw{\gt{u}}_{\bar{0}})$ (so that $G_0$ is the same as in
Proposition~\ref{prop:GLnq:untw:VOAauts}) and we let $\phi:G\to G_0$
be the natural map. $S$ is the subgroup of $G$ generated by
exponentials $\exp(x_{(0)})$ for $x\in\atw{\gt{u}}_1$, and we find
as before that $S$ contains the kernel of $\phi$, and
$G_0=\phi(S)\cup \bar{x}\phi(S)$. We claim that the map $\phi:G\to
G_0$ is in this case not surjective. We have seen already that
$G_0=\Or(\gt{u})/\Kumgp$. If $\phi$ were surjective then $G$ would
be a double cover of $\Or(\gt{u})/\Kumgp$ containing $S$, and the
group $S$ is $\Sp(\gt{u})/\Kumz$ in this case. The only possibility
is that $G=\Pin(\gt{u})/\Kumz$, and this group has no irreducible
representation of dimension $2^{N-1}=\dim \atw{\gt{u}}_{N/8}$. We
conclude that $\phi(G)=\phi(S)$, so that $\Aut(\atw{\gt{u}})=S$ is
the group $\Sp(\gt{u})/\Kumz$.
\end{proof}
Note that $\cej$ belongs to $\atw{\gt{u}}$. It is straightforward to
check, as in the proof of Proposition~
\ref{prop:GLnq:untw:u(1)struc}, that $(\atw{\gt{u}},\cgU)$ is
another example of a $\U(1)$--VOA.
\begin{prop}\label{prop:GLnq:tw:main}
Let $N$ be a positive integer congruent to $4$ modulo $8$. Then the
quadruple $(\atw{\gt{u}},Y,\vac,\cgU)$ is a self-dual $\U(1)$--VOA
of rank $N$.
\end{prop}

\subsection{Special linear groups}\label{sec:GLnq:SLgps}

Suppose we wish to restrict the symmetry of $\atw{\gt{u}}$ further
so as to obtain an action by a group no bigger than the preimage of
$\SL(\gt{a})$. We can achieve this by including elements such as
$\vemp=\vac_{X}$ and $\vful=e_{\hfset}\vac_X$ in some new enhanced
conformal structure on $\atw{\gt{u}}$.
\begin{prop}\label{prop:GLnq:SLgps:Autsjvpm}
The subgroup of $\Aut(\atw{\gt{u}},\cas)$ fixing $\cej$, $\vemp$,
and $\vful$ is the group $\SL(\gt{a})/\Kumgp$.
\end{prop}
\begin{proof}
Let $H$ be the group of the statement of the proposition, and recall
from Proposition~\ref{prop:GLnq:tw:VOAauts} that
$\Aut(\atw{\gt{u}})=\Sp(\gt{u})/\Kumz$. Let
$G=\Aut(\atw{\gt{u}},\cgset_U)$ and
$G_0=\Aut(\atw{\gt{u}}_{\bar{0}},\cgset_U)$. Then the natural map
$\phi:G\to G_0$ is surjective, since by the remark following the
proof of Proposition~\ref{prop:GLnq:untw:eVOAauts} every
automorphism of $(A(\gt{u})_{\bar{0}},\cgset_U)$ is inner (and
$A(\gt{u})_{\bar{0}}=\atw{\gt{u}}_{\bar{0}}$). In fact, $G$ is the
quotient by $\Kumz$ of the group generated by automorphisms
$\pm\exp(tX_{ij})$ where $t\in\CC$ and $X_{ij}$ is the residue of
$Y(x_{ij},z)$ for
$x_{ij}=a_i(-\tfrac{1}{2})a_j^*(-\tfrac{1}{2})\vac$. The
exponentials $\exp(tX_{ij})$ generate a copy of $\GL(\gt{a})/\Kumgp$
in $\Sp(\gt{u})/\Kumz$. Allowing $\pm\exp(tX_{ij})$ we obtain a
double cover of this group. Now consider the action of $G$ on
$\vemp,\vful\in\atw{\gt{u}}_{N/8}$. The subgroup of $G$ fixing
$\vac_X$ is just the group $\GL(\gt{a})/\Kumgp$. This latter group
then preserves all the subspaces $\bigwedge^{2k}(\gt{a})\vac_X$ of
$\atw{\gt{u}}_{N/8}$, and the action is given by
\begin{gather}
     \bar{g}:a_{i_1}\cdots a_{i_{2k}}\vac_X\mapsto
          g(a_{i_1})\cdots g(a_{i_{2k}})\vac_X
\end{gather}
for $g$ a preimage in $\GL(\gt{a})$ of
$\bar{g}\in\GL(\gt{a})/\Kumgp$. Evidently, the coefficient of
$a_{\Delta}\vac_X$ in $\bar{g}a_{\Delta}\vac_X$ is the determinant
of either preimage of $\bar{g}$ in $\GL(\gt{a})$. The claim follows.
\end{proof}
To understand the local algebras (see \S\ref{sec:eVOAs}) arising
from the inclusion of vectors like $\vemp$ and $\vful$ in an
enhanced conformal structure on $\atw{\gt{u}}$ we should calculate
the OPEs involving the corresponding vertex operators. This
calculation is the content of the following proposition. Let
$\vp=\vemp+\vful$ and $\vm=\ii(\vemp-\vful)$. Then we have
\begin{prop}\label{prop:GLnq:SLgps:vpmOPEs}
The operators $\Vp(z)=Y(\vp,z)$ and $\Vm(z)=Y(\vm,z)$ satisfy the
following OPEs
\begin{align}
     \Vpm(z)\Vpm(w)&=\;\pm(-1)^{N/8+1/2}\sum_{k=0}^{N/4-1}
          \frac{P_{k,-\ii\ceJ}(w)+P_{k,\ii\ceJ}(w)}{(z-w)^{N/4-k}}
               +{\rm reg.}\\
     \Vp(z)\Vm(w)&=\;\ii(-1)^{N/8+1/2}\sum_{k=0}^{N/4-1}
          \frac{P_{k,-\ii\ceJ}(w)-P_{k,\ii\ceJ}(w)}{(z-w)^{N/4-k}}+{\rm reg.}
\end{align}
where $P_{k,\pm\ii\ceJ}(z)$ is defined to be the coefficient of
$X^k$ in
\begin{gather}
     \tcolon\exp\left(\sum_{m>0} \pm\ii\ceJ^{(m)}(z)X^m\right)\tcolon
\end{gather}
and $\ceJ^{(m)}(z)=\tfrac{1}{m}D_z^{(m-1)}J(z)$ for $m\in\ZZ_{>0}$.
\end{prop}
We see from Proposition~\ref{prop:GLnq:SLgps:vpmOPEs} that (at least
when $N>4$) if $\Omega\subset\atw{\gt{u}}$ contains $\vp$ and $\vm$
then $\Omega$ and $\Omega\cup\{\cej\}$ determine the same enhanced
conformal structure on $\atw{\gt{u}}$. More precisely, the vertex
Lie subalgebra of $\Sing \atw{\gt{u}}$ generated by $\Omega$
coincides with that generated by $\Omega\cup\{\cej\}$ (see
\S\ref{sec:eVOAs}).

We record some of the observations of this section in the following
\begin{thm}\label{thm:GLnq:SLgps:SL28}
For $N$ congruent to $4$ modulo $8$ and greater than $4$, and for
$\cgset=\{\cas,\vp,\vm\}$, the quadruple
$(\atw{\gt{u}},Y,\vac,\cgset)$ is a self-dual enhanced $\U(1)$--VOA
of rank $N$. The full automorphism group of $(\atw{\gt{u}},\cgset)$
is the group $\SL(\gt{a})/\Kumgp$.
\end{thm}

%%%%%%%%%%%%%%%%%%%%%%%%%%%%%%%%%%%%%%%%%%%%%%%%%%%%%%%%%%%%%%%%%%

\section{The Rudvalis group}\label{sec:Ru}

In this section we realize the sporadic simple group of Rudvalis as
symmetries of an enhanced VOA $\aru$. In \S\ref{sec:Ru:symms} we
show that the full automorphism group of $\aru$ is the direct
product of a cyclic group of order seven with the Rudvalis group.

%\medskip

Our plan for realizing the Rudvalis group as symmetry of an enhanced
VOA is to consider first the enhanced VOA for ${\SL}_N(\CC)$ for
suitable $N$, constructed in \S\ref{sec:GLnq}, and then to find a
single extra (super)conformal generator with which to refine the
conformal structure, just to the point that the Rudvalis group
becomes visible.

We give two constructions of this conformal generator. The first
construction arises directly from the geometry of the Conway--Wales
lattice \cite{ConRu}, and is given in \S\ref{sec:Ru:geom}. The
second construction, which is more convenient for computations, is a
description in terms of a particular maximal subgroup of a double
cover of the Rudvalis group, and is given in \S\ref{sec:Ru:mnml}.
Sections \ref{sec:Ru:geom} and \ref{sec:Ru:mnml} are independent,
and the reader my safely skip one in favor of the other. The
approach of \S\ref{sec:Ru:geom} is certainly more conceptual, and
more brief.

The enhanced VOA structure for $\aru$ is described in
\S\ref{sec:Ru:cnst}. In \S\ref{sec:Ru:symms} we determine its
symmetry group, and conjecture a characterization.

\subsection{Geometric description}\label{sec:Ru:geom}

Let $\LL$ denote the Conway--Wales lattice \cite{ConRu}, viewed as a
module of rank $28$ over the Gaussian integers $\ZZ[\ii]$, and
equipped with a non-degenerate Hermitian form denoted
$(\cdot\,,\cdot)$. (An explicit description of this lattice is given
in the sequel \cite{DunVARuII}.) We write $\widetilde{R}$ for
$\Aut(\LL)$, the subgroup of the unitary group of the space
$\gt{r}=\CC\otimes_{\ZZ[\ii]}\LL$ that preserves the set of vectors
in $\LL$. Then $\widetilde{R}$ is a four-fold cover of the Rudvalis
group $\Ru$, and may be written as a central product $4\circ
(2.\!\Ru)$ where the $4$ is generated by multiplication by $\ii$,
and the non-trivial central element of the perfect cover $2.\!\Ru$
is multiplication by $-1$  \cite{ConRu}. Let $\LL_2$ denote the set
of type $2$ vectors in $\LL$; i.e. the vectors $\lambda\in\LL$ with
$(\lambda,\lambda)=4$. These are also called the {\em sacred
vectors} \cite{ConRu}. Counting projectively, there are $4060$
sacred vectors in $\LL$, and $u\lambda$ is a sacred vector whenever
$\lambda$ is sacred and $u$ is a unit in $\ZZ[\ii]$. %We write
%$\gt{r}=\CC\otimes_{\ZZ[\ii]}\LL$ for the ambient Hermitian space
%containing $\LL$.

For $X\subset\LL_2$ we define $\Gamma(X)$ to be the (undirected)
graph with vertices $\{v_{\lambda}\}_{\lambda\in X}$ indexed by the
$\lambda$ in $X$, and edges such that $v_{\lambda}$ is joined to
$v_{\lambda'}$ just when $(\lambda,\lambda')=1$. We define
$\Delta(X)$ to be the directed graph with vertices
$\{v_{\lambda}\}_{\lambda\in X}$ and a directed edge from
$v_{\lambda}$ to $v_{\lambda'}$ just when $(\lambda,\lambda')=\ii$.

Let $B\subset\LL_2$ be an orbit of size $13$ for some element of
order $13$ in $\Rtd$ such that $\Delta(B)$ has $13$ (directed)
edges. (All elements of order $13$ are conjugate in $\Rtd$
\cite{ATLAS}, so it matters not which one we choose.) Then
$\Delta(B)$ is in fact an oriented triangulation of the circle, and
for each $\lambda$ in $B$ there is a unique $\lambda'\in B$ such
that $(\lambda, \lambda')=\ii$. In addition to this, the graph
$\Gamma(B)$ has $26$ (undirected) edges, and is a triangulation of a
M\"obius band, or at least, the graph one would obtain by taking a
triangulated M\"obius band (with $13$ triangles) and placing a
(graph) edge along the ``edge'' of each $2$-cell, and vertices at
the intersection of each (graph) edge.
\begin{defn}
We call $B\subset\LL_2$ an {\em M-set} in case $\Gamma(B)$ is a
triangulation of a M\"obius band (in the above sense), and
$\Delta(B)$ is an oriented triangulation of a circle, each
triangulation having $13$ triangles.
\end{defn}

It turns out that any M-set in $\LL_2$ is an orbit for some element
of order $13$ in $\Rtd$. Considering the M-sets that are orbits for
a fixed element of order $13$ in $\Rtd$, we can easily find elements
of $\Rtd$ sending one to another. Recalling that there is a single
conjugacy class of elements of order $13$ in $\Rtd$ we then have
\begin{prop}
The M-sets in $\LL_2$ constitute a single orbit under $\Rtd$.
\end{prop}

For $B$ an M-set, we say that a labeling $B=\{\beta^1, \ldots,
\beta^{13}\}$ of the elements of $B$ is {\em oriented} if we have
$(\beta^i,\beta^{i+1})=\ii$ for each $i\in\{1,\ldots,12\}$. Then all
oriented labelings of $B$ are equivalent under cyclic permutations
of the $\beta^i$, and the expression
$\beta^1\wedge\beta^2\wedge\cdots\wedge\beta^{13}$ is a well defined
element of $\bigwedge^{13}(\gt{r})$ that is independent of the
choice of oriented labeling. In Table \ref{tab:bstar} we give an
example of an M-set; we denote this particular example $B_*$. Our
notation for vectors in $\LL_2$ is similar to that used in
\cite{WilRu}: we view these vectors as $7$-tuples of complex
quaternions as in \cite{ConRu}, and we present such an object as a
$4\times 7$ array, but with the coefficients of $1$, ${\bf j}$,
${\bf k}$ and ${\bf l}$ (rather than $\ii$, ${\bf j}$, ${\bf k}$ and
${\bf l}$) appearing in the respective columns. We obtain a
particular oriented labeling for $B_*$ by reading the vectors in
Table \ref{tab:bstar} from left to right within rows, from the top
row to the bottom. We denote this particular labeling by
$B_*=\{\beta_*^1,\ldots,\beta_*^{13}\}$.
\begin{table}
   \centering
   \caption{The M-set $B_*$}
   \label{tab:bstar}
%\begin{small}
\begin{gather*}
\begin{array}{|rrrr|}%\hline
     0 & 0 & 0 & 0\\
     0 & 0 & 0 & 0\\
     2 & 0 & 0 & 0\\
     0 & 0 & 0 & 0\\
     0 & 0 & 0 & 0\\
     0 & 0 & 0 & 0\\
     0 & 0 & 0 & 0\\%\hline
\end{array}\\
\frac{1}{2}
\begin{array}{|rrrr|}%\hline
     0 & 0 & 0 & 0\\
     0 & 0 & 0 & 0\\
     \bar{\ii} & \bar{1} & \ii & \bar{1}\\
     0 & 0 & 0 & 0\\
     1 & 1 & \bar{\ii} & \ii\\
     1 & \ii & \ii & \bar{1}\\
     \bar{1} & \bar{1} & \bar{1} & \bar{1}\\%\hline
\end{array}\quad
\frac{1}{2}
\begin{array}{|rrrr|}%\hline
     \bar{1} & \ii & 0 & 0\\
     0 & 0 & \ii & \bar{\ii}\\
     0 & \ii & 0 & \bar{1}\\
     0 & \ii & \bar{\ii} & 0\\
     0 & \bar{1} & 0 & 1\\
     \bar{1} & 1 & 1 & \bar{\ii}\\
     0 & \ii & \bar{1} & 0\\%\hline
\end{array}\quad
\frac{1}{2}
\begin{array}{|rrrr|}%\hline
     0 & 0 & \ii & \bar{1}\\
     0 & 0 & \bar{1} & \bar{1}\\
     0 & 1 & 0 & \ii\\
     \ii & 0 & 0 & \ii\\
     0 & \ii & 0 & \bar{\ii}\\
     \ii & \bar{1} & \ii & \ii\\
     0 & \ii & \bar{1} & 0\\%\hline
\end{array}\quad
\frac{1}{2}
\begin{array}{|rrrr|}%\hline
     0 & \ii & \ii & 0\\
     0 & \bar{1} & 0 & 1\\
     1 & 1 & \ii & 1\\
     0 & \ii & \bar{1} & 0\\
     0 & 0 & \bar{\ii} & \bar{1}\\
     \bar{\ii} & \ii & 0 & 0\\
     0 & 1 & 0 & \ii\\%\hline
\end{array}\\
\frac{1}{2}
\begin{array}{|rrrr|}%\hline
     1 & 0 & 0 & \bar{1}\\
     0 & 1 & 0 & 1\\
     1 & \bar{1} & 1 & \bar{\ii}\\
     0 & 1 & \ii & 0\\
     1 & \ii & 0 & 0\\
     0 & 0 & 1 & \bar{1}\\
     \ii & 0 & \bar{1} & 0\\%\hline
\end{array}\quad
\frac{1}{2}
\begin{array}{|rrrr|}%\hline
     \bar{1} & 0 & 0 & \ii\\
     0 & 0 & 1 & \bar{\ii}\\
     0 & 0 & 1 & \bar{1}\\
     0 & \bar{1} & 0 & \ii\\
     0 & \ii & \bar{\ii} & 0\\
     0 & \ii & 0 & \bar{\ii}\\
     1 & \bar{1} & \ii & \bar{1}\\%\hline
\end{array}\quad
\frac{1}{2}
\begin{array}{|rrrr|}%\hline
     \ii & \ii & \ii & \bar{1}\\
     0 & \ii & \bar{1} & 0\\
     0 & 0 & \bar{\ii} & \bar{1}\\
     0 & 0 & \bar{\ii} & \ii\\
     0 & 1 & 0 & \bar{\ii}\\
     0 & 1 & 1 & 0\\
     \bar{\ii} & 0 & \bar{\ii} & 0\\%\hline
\end{array}\quad
\frac{1}{2}
\begin{array}{|rrrr|}%\hline
     1 & 0 & 0 & \ii\\
     0 & 0 & \bar{1} & \bar{\ii}\\
     1 & 1 & 0 & 0\\
     0 & \bar{1} & 0 & \ii\\
     0 & \bar{1} & \bar{1} & 0\\
     1 & 0 & \bar{1} & 0\\
     \bar{1} & 1 & \bar{1} & \ii\\%\hline
\end{array}\\
\frac{1}{2}
\begin{array}{|rrrr|}%\hline
     1 & 0 & \ii & 0\\
     0 & \bar{\ii} & \ii & 0\\
     1 & 0 & 1 & 0\\
     \bar{1} & \ii & \bar{\ii} & \bar{\ii}\\
     0 & 1 & \bar{\ii} & 0\\
     0 & 0 & \ii & \bar{1}\\
     \ii & \bar{\ii} & 0 & 0\\%\hline
\end{array}\quad
\frac{1}{2}
\begin{array}{|rrrr|}%\hline
     \bar{\ii} & 0 & 1 & 0\\
     0 & \bar{\ii} & \bar{\ii} & 0\\
     0 & \bar{1} & 0 & \bar{1}\\
     1 & \ii & \bar{1} & \bar{1}\\
     0 & \ii & \bar{1} & 0\\
     0 & 0 & 1 & \bar{\ii}\\
     0 & 0 & \bar{1} & \bar{1}\\%\hline
\end{array}\quad
\frac{1}{2}
\begin{array}{|rrrr|}%\hline
     \bar{1} & \bar{1} & 1 & \ii\\
     0 & \ii & \bar{1} & 0\\
     0 & 0 & 1 & \ii\\
     0 & 0 & \bar{\ii} & \ii\\
     \ii & 0 & 1 & 0\\
     \bar{1} & 0 & 0 & \bar{1}\\
     1 & 0 & \bar{1} & 0\\%\hline
\end{array}\quad
\frac{1}{2}
\begin{array}{|rrrr|}%\hline
     0 & \ii & \bar{1} & 0\\
     0 & 0 & \bar{1} & \bar{\ii}\\
     \ii & \bar{\ii} & 0 & 0\\
     \ii & 0 & \bar{1} & 0\\
     \bar{\ii} & 0 & 0 & \bar{\ii}\\
     0 & \ii & 0 & \ii\\
     \bar{\ii} & \bar{\ii} & \bar{1} & \ii\\%\hline
\end{array}
\end{gather*}
%\end{small}
\end{table}

It turns out that an M-set $B=\{\beta^1,\ldots,\beta^{13}\}$ wields
much control over the other sacred vectors in $\LL$, in the sense
that, generically, a sacred vector $\lambda\in\LL_2$ is determined
uniquely among all sacred vectors by the sequence of values
$(\beta^k,\lambda)$ for $k=1,\ldots,13$. Important exceptions to
this principle are the vectors $\lambda\in\LL_2$ satisfying
$(B,\lambda)=\{u_0\}$ for some fixed unit $u_0$ in $\ZZ[\ii]$. (We
write $(B,\lambda)$ as a shorthand for the set
$\{(\beta,\lambda)\mid \beta\in B\}$.) For example, for any M-set
$B$ there are exactly two vectors $\mu,\mu'$ in $\LL_2$ such that
$(B,\mu)=(B,\mu')=\{1\}$. Nonetheless, we can tell them apart, since
the inner product $(\mu,\mu')$ is always non-zero and non-real, and
hence depends on the order. We call $\mu$ and $\mu'$ the {\em
neighbors} to $B$ if $(B,\mu)=(B,\mu')=\{1\}$ (and $\mu\neq\mu'$),
and we say that $\mu$ is the {\em first neighbor} if
$(\mu,\mu')=\ii$. The vector $\mu'$ is then called the {\em second
neighbor} to $B$. The first and second neighbors to $B_*$ will be
denoted $\mu_*$ and $\mu_*'$, respectively. They appear in Table
\ref{tab:bstarnbrs}.

\begin{table}
   \centering
   \caption{Neighbors to $B_*$}
   \label{tab:bstarnbrs}
%\begin{small}
\begin{gather*}
\mu_*=\frac{1}{2}
\begin{array}{|rrrr|}%\hline
     0 & \ii & 0 & \bar{1}\\
     \ii & 0 & 0 & \bar{\ii}\\
     1 & 0 & 1 & 0\\
     1 & \ii & \bar{\ii} & \ii\\
     \bar{\ii} & 0 & 0 & \bar{1}\\
     0 & 0 & \ii & \bar{1}\\
     0 & 0 & \bar{1} & \bar{1}\\%\hline
\end{array}\quad
\mu_*'=\frac{1}{2}
\begin{array}{|rrrr|}%\hline
     0 & 0 & \ii & \ii\\
     \ii & 0 & \bar{1} & 0\\
     1 & 0 & 0 & \bar{1}\\
     \bar{1} & 0 & \bar{1} & 0\\
     \ii & \ii & \bar{\ii} & 1\\
     0 & \ii & 1 & 0\\
     0 & 0 & \bar{1} & \ii\\%\hline
\end{array}
\end{gather*}
%\end{small}
\end{table}

Consider then the case that $u_0=0$. Counting projectively, there
are just two vectors satisfying $(B,\lambda)=\{0\}$. Again we can
distinguish them, at least up to multiplication by units, for if
$\lambda\in\LL_2$ satisfies $(B,\lambda)=\{0\}$ then for $\mu,\mu'$
the neighbors to $B$, we have either $(\lambda,\mu)=(\lambda,\mu')$
or $(\lambda,\mu)=-(\lambda,\mu')$. We say that $\lambda$ is a {\em
complement} to $B$ if $(B,\lambda)=\{0\}$, and we say that a
complement $\lambda$ is {\em positive} (resp. {\em negative}) if
$\frac{(\lambda,\mu)}{(\lambda,\mu')}$ is $+1$ (resp. $-1$). If
$\lambda$ and $\lambda'$ are positive and negative complements to an
M-set $B$, then they are not orthogonal. In particular, for a given
positive complement $\lambda$, there is a unique negative complement
$\lambda'$ such that $(\lambda,\lambda')=1$. In Table
\ref{tab:bstarcplts} we give a positive complement $\lambda_*$ for
$B_*$, and a negative complement $\lambda_*'$ satisfying
$(\lambda_*,\lambda_*')=1$.
\begin{table}
   \centering
   \caption{Complements to $B_*$}
   \label{tab:bstarcplts}
%\begin{small}
\begin{gather*}
\lambda_*=
\begin{array}{|rrrr|}%\hline
     0 & 0 & 0 & 0\\
     2 & 0 & 0 & 0\\
     0 & 0 & 0 & 0\\
     0 & 0 & 0 & 0\\
     0 & 0 & 0 & 0\\
     0 & 0 & 0 & 0\\
     0 & 0 & 0 & 0\\%\hline
\end{array}\quad
\lambda_*'=\frac{1}{2}
\begin{array}{|rrrr|}%\hline
     \bar{1} & 1 & 1 & \bar{\ii}\\
     1 & 0 & 0 & \bar{\ii}\\
     0 & 0 & \bar{\ii} & 1\\
     1 & \bar{1} & 0 & 0\\
     \bar{\ii} & 0 & 1 & 0\\
     1 & 0 & 0 & \bar{1}\\
     \ii & 0 & \bar{\ii} & 0\\%\hline
\end{array}
\end{gather*}
%\end{small}
\end{table}

If $B$ is an M-set and $\lambda$ and $\lambda'$ are positive and
negative complements to $B$ respectively, then there are M-sets for
which $\lambda$ is a negative complement, and $\lambda'$ is a
positive complement. We can single out one of these $C$ say, by
defining $C=\{\gamma^1,\ldots,\gamma^{13}\}$ to be the set of
vectors such that the inner products $(\beta^i,\gamma^j)$ take on
certain specified values. If $B$ is an M-set, then there is a unique
vector $\gamma$ such that the sequence $(\beta^i,\gamma)$ coincides
with that given in (\ref{eqn:partnercond}), and similarly for each
cyclic permutation of (\ref{eqn:partnercond}).
\begin{gather}\label{eqn:partnercond}
(-\ii, -\ii, 0, 0, \ii, -\ii, 0, 0, 0, 0, -\ii, \ii, 0)
\end{gather}
We say that $C=\{\gamma^j\}$ is a {\em right partner} to an M-set
$B$ if the $\gamma^j$ are all the vectors such that the sequence
$(\beta^i,\gamma^j)$ (with varying $i$) is a cyclic permutation of
(\ref{eqn:partnercond}). It turns out that $C$ is again an M-set,
with the same complements as $B$, and that a positive complement to
$C$ is a negative complement to $B$, and vice-versa. To prove these
statements it suffices to check one example, so we exhibit a right
partner to $B_*$ in Table \ref{tab:cstar}. This right partner to
$B_*$ will be denoted $C_*$. We write
$\{\gamma_*^1,\ldots,\gamma_*^{13}\}$ for the particular oriented
labeling of $C_*$ obtained by reading from left to right, and top to
bottom in Table \ref{tab:cstar}. If $C$ is a right partner to $B$,
then $-B$ is a right partner to $C$.
\begin{table}
   \centering
   \caption{The M-set $C_*$}
   \label{tab:cstar}
%\begin{small}
\begin{gather*}
\frac{1}{2}
\begin{array}{|rrrr|}%\hline
     0 & 0 & 1 & 1\\
     0 & \ii & 0 & \bar{1}\\
     \ii & 0 & 0 & \bar{\ii}\\
     \bar{1} & 0 & \bar{1} & 0\\
     1 & \ii & \bar{\ii} & \ii\\
     \ii & 0 & 0 & 1\\
     0 & 0 & \bar{\ii} & 1\\%\hline
\end{array}\\
\frac{1}{2}
\begin{array}{|rrrr|}%\hline
     \bar{1} & \bar{\ii} & 0 & 0\\
     0 & 0 & \bar{1} & 1\\
     \ii & 0 & \bar{1} & 0\\
     \bar{1} & 0 & 0 & 1\\
     0 & 1 & 0 & 1\\
     1 & \bar{1} & 1 & \bar{\ii}\\
     0 & 1 & \ii & 0\\%\hline
\end{array}\quad
\frac{1}{2}
\begin{array}{|rrrr|}%\hline
     0 & \bar{\ii} & 0 & \bar{1}\\
     0 & \bar{1} & 1 & 0\\
     0 & 1 & 0 & \bar{1}\\
     \ii & \bar{\ii} & \bar{1} & \bar{\ii}\\
     \ii & 0 & 0 & 1\\
     0 & 0 & \bar{\ii} & \bar{1}\\
     0 & 0 & \bar{\ii} & \ii\\%\hline
\end{array}\quad
\frac{1}{2}
\begin{array}{|rrrr|}%\hline
     0 & 0 & \ii & 1\\
     0 & 0 & \bar{\ii} & \ii\\
     0 & \bar{\ii} & 0 & \bar{1}\\
     0 & \bar{1} & 1 & 0\\
     0 & \bar{1} & 0 & 1\\
     \bar{\ii} & \ii & 1 & \ii\\
     \ii & 0 & 0 & 1\\%\hline
\end{array}\quad
\frac{1}{2}
\begin{array}{|rrrr|}%\hline
     0 & 1 & 0 & \ii\\
     0 & 1 & \bar{1} & 0\\
     \bar{\ii} & 0 & \ii & 0\\
     \bar{1} & \bar{1} & \bar{1} & \bar{\ii}\\
     0 & \ii & 1 & 0\\
     0 & 0 & \bar{\ii} & 1\\
     0 & 0 & \bar{\ii} & \bar{\ii}\\%\hline
\end{array}\\
\frac{1}{2}
\begin{array}{|rrrr|}%\hline
     \ii & \bar{\ii} & \ii & 1\\
     0 & \ii & \bar{1} & 0\\
     \ii & \bar{1} & 0 & 0\\
     0 & 0 & \ii & \bar{\ii}\\
     \bar{1} & 0 & \bar{\ii} & 0\\
     \ii & 0 & 0 & \bar{\ii}\\
     0 & \ii & 0 & \ii\\%\hline
\end{array}\quad
\frac{1}{2}
\begin{array}{|rrrr|}%\hline
     \bar{1} & \bar{1} & 0 & 0\\
     0 & \bar{\ii} & 0 & 1\\
     0 & \ii & \bar{\ii} & 0\\
     0 & \bar{\ii} & 0 & \bar{\ii}\\
     \ii & \ii & \bar{\ii} & 1\\
     1 & 0 & 0 & \ii\\
     \bar{\ii} & 1 & 0 & 0\\%\hline
\end{array}\quad
\frac{1}{2}
\begin{array}{|rrrr|}%\hline
     0 & 1 & 0 & \bar{\ii}\\
     0 & \bar{1} & \bar{1} & 0\\
     0 & 1 & 0 & \bar{1}\\
     \bar{\ii} & \bar{\ii} & \ii & \bar{1}\\
     \ii & 0 & 0 & 1\\
     \bar{1} & \ii & 0 & 0\\
     \ii & \bar{\ii} & 0 & 0\\%\hline
\end{array}\quad
\frac{1}{2}
\begin{array}{|rrrr|}%\hline
     0 & 0 & 1 & \ii\\
     0 & 0 & \ii & \ii\\
     0 & \bar{\ii} & 0 & 1\\
     \bar{1} & 0 & 0 & \bar{1}\\
     0 & \bar{1} & 0 & 1\\
     1 & \ii & 1 & 1\\
     0 & \bar{1} & \bar{\ii} & 0\\%\hline
\end{array}\\
\frac{1}{2}
\begin{array}{|rrrr|}%\hline
     0 & 0 & \bar{1} & 1\\
     0 & \ii & 0 & \bar{1}\\
     0 & \bar{1} & \bar{1} & 0\\
     \bar{1} & 0 & \bar{1} & 0\\
     \ii & \ii & 1 & \bar{\ii}\\
     1 & 0 & 0 & \bar{\ii}\\
     \ii & 1 & 0 & 0\\%\hline
\end{array}\quad
\frac{1}{2}
\begin{array}{|rrrr|}%\hline
     0 & \bar{\ii} & \ii & 0\\
     0 & 1 & 0 & 1\\
     \ii & \ii & \bar{1} & \ii\\
     0 & \bar{\ii} & \bar{1} & 0\\
     0 & 0 & \bar{\ii} & 1\\
     \bar{\ii} & \bar{\ii} & 0 & 0\\
     0 & \bar{1} & 0 & \ii\\%\hline
\end{array}\quad
\frac{1}{2}
\begin{array}{|rrrr|}%\hline
     0 & \bar{1} & 0 & \bar{\ii}\\
     0 & \bar{\ii} & \bar{\ii} & 0\\
     \bar{\ii} & 0 & \bar{\ii} & 0\\
     \bar{1} & \bar{1} & \ii & 1\\
     \ii & 0 & 0 & 1\\
     \bar{1} & \ii & 0 & 0\\
     0 & 0 & \bar{\ii} & \ii\\%\hline
\end{array}\quad
\frac{1}{2}
\begin{array}{|rrrr|}%\hline
     \bar{\ii} & 1 & \ii & \ii\\
     0 & \bar{\ii} & \bar{1} & 0\\
     0 & 0 & \bar{1} & \bar{\ii}\\
     0 & 0 & 1 & \bar{1}\\
     \ii & 0 & 1 & 0\\
     0 & \ii & \bar{\ii} & 0\\
     0 & \bar{1} & 0 & 1\\%\hline
\end{array}
\end{gather*}
%\end{small}
\end{table}

Let us define an undirected graph $\Xi$ whose vertices are pairs
$[\lambda, B]$ where $B$ is an M-set and $\lambda$ is a positive
complement to $B$, and whose edges are determined as follows: the
vertices $[\lambda_0,B_0]$ and $[\lambda_1,B_1]$ are joined if any
one of the following situations hold.
\begin{itemize}
\item[$(\Xi1)$]     $\lambda_0=u\lambda_1$ and $B_0=\bar{u}B_1$ for some unit
$u$ in $\ZZ[\ii]$.
\item[$(\Xi2)$]     $\lambda_0=\lambda_1$ and $(\lambda_0, \mu_0)= (\lambda_0,
\mu_1)$ where $\mu_0$ is the first neighbor to $B_0$, and $\mu_1$ is
the first neighbor to $B_1$.
\item[$(\Xi3)$]     $\lambda_1$ is the negative complement to $B_0$ for which
$(\lambda_0,\lambda_1)=-\ii$, and $B_1$ is the right partner to
$B_0$.
\end{itemize}
(For $u$ a unit and $B$ an M-set, $uB$ denotes the M-set
$\{u\beta\mid \beta\in B\}$.) Given a vertex $[\lambda,B]$ in $\Xi$,
we write $a_{[\lambda,B]}$ for the element of
$\bigwedge^{14}(\gt{r})$ given by
\begin{gather}\label{eqn:Ru:geom:adefn}
     a_{[\lambda,B]}=\lambda\wedge
               \beta^1\wedge \cdots\wedge\beta^{13}
\end{gather}
where $\{\beta^1,\ldots,\beta^{13}\}$ is any oriented labeling of
$B$. Note that $a_{[\lambda,B]}$ is independent of the oriented
labeling chosen. Observe also that $\Rtd$ acts naturally on the
pairs $[\lambda,B]$, and on $\bigwedge^{14}(\gt{r})$, and these
actions are compatible in the sense that
$ga_{[\lambda,B]}=a_{[g\lambda,gB]}$ for $g\in\Rtd$.

\medskip

We are now ready to define the extra superconformal generator
$\spsp$, which will determine the enhanced VOA structure that we
will utilize in \S\ref{sec:Ru:cnst}. Let $\Xi_*$ be the connected
component of $\Xi$ containing the vertex $[\lambda_*,B_*]$, and
define $\spsp\in\bigwedge^{14}(\gt{r})$ by setting
\begin{gather}\label{eqn:Ru:geom:spspdefn}
     \spsp=\sum_{[\lambda,B]\in\Xi_*}a_{[\lambda,B]}
\end{gather}
where $a_{[\lambda,B]}$ is defined in (\ref{eqn:Ru:geom:adefn}).

We claim that $\spsp$ is invariant for the natural action of the
Rudvalis group on $\bigwedge^{14}(\gt{r})$. This claim will follow
from Proposition~\ref{prop:Ru:geom:RhtPrsvsCmpts} below, and in
preparation for this we now introduce some notation and a few
lemmas. Let $R$ be the unique subgroup of $\Rtd=\Aut(\LL)$ such that
$R$ is isomorphic to $\cRu$. (If $R'$ is another such group, then
its intersection with $R$ contains all elements of odd order in $R$.
These elements generate a group that is non-central and normal in
both $R'$ and $R$, and hence is both $R'$ and $R$.) For
$\lambda\in\LL$, we write $\Fix(\lambda)$ for the subgroup of $\Rtd$
consisting of $g$ for which $g\lambda=\lambda$. Then $R$ contains
$\Fix(\lambda)$ for all $\lambda$ in $\LL_2$.

\begin{lem}\label{lem:Ru:geom:FixGpTransOnMSets}
For $\lambda\in\LL_2$ the group $\Fix(\lambda)$ acts projectively
transitively on M-sets for which $\lambda$ is a positive complement.
That is, if $B_1$ and $B_2$ are M-sets for which $\lambda$ is a
positive complement then there exists $g\in\Fix(\lambda)$ and a unit
$u$ such that $gB_1=uB_2$.
\end{lem}

\begin{lem}\label{lem:Ru:geom:IntrchgPtnrs}
Let $B$ be an $M$-set, and let $\nu_1$ and $\nu_2$ be positive and
negative complements to $B$, respectively, with
$(\nu_1,\nu_2)=-\ii$. Then there exists an element $s\in R$ with
$s^2=-\Id$ such that $s:\nu_1\mapsto \nu_2$ and $s(B)=C$ where $C$
is the right partner to $B$.
\end{lem}

\begin{lem}\label{lem:Ru:geom:TrianglesExist}
If $\nu_0,\nu_2\in\LL_2$ and $\nu_2$ is not a scalar multiple of
$\nu_0$ then there exists a vector $\nu_1\in\LL_2$ such that
$(\nu_0,\nu_1)=-\ii$ and $(\nu_1,\nu_2)=-\ii$.
\end{lem}

\begin{prop}\label{prop:Ru:geom:RhtPrsvsCmpts}
If $g\in R$ and $[\lambda_0,B_0]\in\Xi$ then $[\lambda_0,B_0]$ and
$[g\lambda_0,gB_0]$ belong to the same connected component of $\Xi$.
\end{prop}
\begin{proof}
We have several different cases to treat, depending on the value of
$(\lambda_0,g\lambda_0)$. For given $[\lambda_0,B_0]$ let us write
$\Xi_0$ for the connected component of $\Xi$ containing
$[\lambda_0,B_0]$. If we establish that $[g\lambda_0,gB_0]\in\Xi_0$
whenever $(\lambda_0,g\lambda_0)=a$ say, then the same holds also
when $(\lambda_0,g\lambda_0)=-a$, since then
$[-g\lambda_0,-gB_0]\in\Xi_0$, and $[-g\lambda_0,-gB_0]$ is joined
to $[g\lambda_0,gB_0]$ by $(\Xi1)$.

{\em Case (i): $(\lambda_0,g\lambda_0)=4$.} This is the case that
$g\in\Fix(\lambda_0)$. Let $\mu_0$ be the first neighbor to $B_0$
and observe that $g\mu_0$ is the first neighbor to $gB_0$. Then we
have $(\lambda_0,\mu_0)=(g\lambda_0,g\mu_0)=(\lambda_0,g\mu_0)$ so
that $[\lambda_0,B_0]$ is joined to
$[g\lambda_0,gB_0]=[\lambda_0,gB_0]$ curtesy of $(\Xi2)$.

{\em Case (ii): $(\lambda_0,g\lambda_0)=-\ii$.} By Lemma
\ref{lem:Ru:geom:FixGpTransOnMSets} there exists
$g_0\in\Fix(\lambda_0)$ such that $\lambda_0$ and $g\lambda_0$ are
positive and negative complements to $B_1:=g_0B_0$, respectively. By
Lemma \ref{lem:Ru:geom:IntrchgPtnrs} there exists $s_1\in R$ such
that $s_1^2=-\Id$ and $s_1\lambda_0=g\lambda_0$, and $B_2:=s_1(B_1)$
is the right partner to $B_1$. Now set $h=s_1g_0$ and $g_2=gh^{-1}$,
and observe that $g_2\in\Fix(g\lambda_0)$. The decomposition
$g=g_2s_1g_0$ corresponds in the following way to a path in $\Xi_0$.
We set $[\lambda_1,B_1]=[g_0\lambda_0,g_0B_0]$,
$[\lambda_2,B_2]=[s_1\lambda_1,s_1B_1]$, and
$[\lambda_3,B_3]=[g_2\lambda_2,g_2B_2]$. Then the vertices
$[\lambda_0,B_0]$ and $[\lambda_1,B_1]$ are joined curtesy of the
argument in case (i), the vertices $[\lambda_1,B_1]$ and
$[\lambda_2,B_2]$ are joined thanks to $(\Xi3)$, and the vertices
$[\lambda_2,B_2]$ and $[\lambda_3,B_3]=[g\lambda_0,gB_0]$ are joined
thanks again to case (i).

The remaining cases are dealt with in an analogous way to case (ii),
except that we may require a longer decomposition
$g=g_{2k}s_{2k-1}\cdots g_2s_1g_0$, where $k+1$ is the length of a
sequence $\lambda_0,\lambda_2,\ldots,\lambda_{2k} \in\LL_2$ such
that $\lambda_{2k}=g\lambda_0$ and
$(\lambda_{2m},\lambda_{2m+2})=-\ii$ for each $m$. In each case, the
existence of a such a sequence is insured by Lemma
\ref{lem:Ru:geom:TrianglesExist}.

{\em Case (iii): $(\lambda_0,g\lambda_0)\in\{0,1\}$.} Set
$\lambda_4=g\lambda_0$. By Lemma \ref{lem:Ru:geom:TrianglesExist}
there exists $\lambda_2\in\LL_2$ such that
$(\lambda_0,\lambda_2)=-\ii$ and $(\lambda_2,\lambda_4)=-\ii$. As
indicated above, we proceed just as in case (ii). Pick
$g_0\in\Fix(\lambda_0)$ such that $\lambda_0$ and $\lambda_2$ are
positive and negative complements to $B_1:=g_0B_0$, respectively,
and choose $s_1\in R$ such that $s_1^2=-\Id$ and
$s_1\lambda_0=\lambda_2$ and $B_2:=s_1(B_1)$ is the right partner to
$B_1$. Repeating this, pick $g_2\in\Fix(\lambda_2)$ such that
$\lambda_2$ and $\lambda_4$ are positive and negative complements to
$B_3:=g_2B_2$, respectively, and choose $s_3\in R$ such that
$s_3^2=-\Id$ and $s_2\lambda_2=\lambda_4$ and $B_4:=s_3(B_3)$ is the
right partner to $B_3$. Set $h=s_3g_2s_1g_0$ and observe that
$g_4:=gh^{-1}$ fixes $\lambda_4=g\lambda_0$. Consider the sequence
$[\lambda_k,B_k]$ in $\Xi$ obtained by setting
$\lambda_{2k+1}=\lambda_{2k}$ for $0\leq k\leq 2$, and $B_5=g_4B_4$,
so that
\begin{gather}
\begin{split}
     [\lambda_1,B_1]&=[g_0\lambda_0,g_0B_0]\\
     [\lambda_2,B_2]&=[s_1\lambda_1,s_1B_1]\\
          &\;\;\vdots\\
     [\lambda_5,B_5]&=[g_4\lambda_4,g_4B_4]
\end{split}
\end{gather}
and $[\lambda_5,B_5]=[g\lambda_0,gB_0]$. Then $[\lambda_{m},B_{m}]$
is joined to $[\lambda_{m+1},B_{m+1}]$ by an application of case (i)
for $m$ even, and by $(\Xi3)$ for $m$ odd.

{\em Case (iv): $(\lambda_0,g\lambda_0)=4\ii$.} In this case we
require a sequence $\lambda_0,\lambda_2,\ldots,\lambda_6$ of length
$4$. Choose $\lambda_4\in\LL_2$ such that
$(\lambda_0,\lambda_4)=-1$. By Lemma
\ref{lem:Ru:geom:TrianglesExist} there exists $\lambda_2 \in\LL_2$
such that $(\lambda_0, \lambda_2)=(\lambda_2, \lambda_4)=-\ii$. Let
$\lambda_6=-\ii\lambda_0=g\lambda_0$. Then by construction,
$(\lambda_4,\lambda_6)=-\ii$ also. This furnishes a sequence
$\lambda_0,\ldots,\lambda_6$ as required. We then proceed as in case
(iii) using a decomposition of the form $g=g_6s_5g_4s_3g_2s_1g_0$
and the corresponding length $8$ sequence of vertices in $\Xi$, each
one joined to its successor thanks to alternating applications of
case (i) and $(\Xi3)$.

We have accounted for all cases, and the proof is complete.
\end{proof}

Proposition~\ref{prop:Ru:geom:RhtPrsvsCmpts} shows that the
$a_{[\lambda,B]}$ for $[\lambda,B]$ in $\Xi_*$ --- which are the
summands in the expression for $\spsp$ in
(\ref{eqn:Ru:geom:spspdefn})
--- constitute a union of orbits for $R$ in $\bigwedge^{14}(\gt{r})$. Of
course, the central element $-\Id$ in $R$ acts trivially on
$\bigwedge^{14}(\gt{r})$, so we may view this set as a union of
orbits for a copy $R/\Kumgp$ of the simple group $\Ru$. We have
established
\begin{thm}\label{thm:Ru:geom:FhasRu}
The vector $\spsp$ is invariant for the action of the Rudvalis group
on $\bigwedge^{14}(\gt{r})$.
\end{thm}
In fact the summands of $\spsp$ are just a single orbit for the
Rudvalis group, since from the proof of Proposition~
\ref{prop:Ru:geom:RhtPrsvsCmpts} we see that two vertices of $\Xi$
are joined only if they define the same element of
$\bigwedge^{14}(\gt{r})$ (this is condition $(\Xi1)$), or if they
are joined by some element of $R$ (these are conditions $(\Xi2)$ and
$(\Xi3)$).

\subsection{Monomial description}\label{sec:Ru:mnml}

We now give a description of $\spsp$ in terms of the monomial action
of a group which will turn out to be a maximal subgroup of a double
cover of the Rudvalis group. This group has the shape $2^7.G_2(2)$,
and we will refer to the copy we construct in
\S\ref{sec:Ru:mnml:monomialgp} as {\em the monomial group}. The
quotient group $G_2(2)$ may be viewed as the symmetry group of a
certain algebraic structure arising from the $E_8$ lattice, and in
\S\ref{sec:Ru:mnml:cayleyalg} we review these
relationships. %To some extent this elucidates the significance of
%the $N$ that is chosen when we construct the VOA underlying $\aru$
%in \S\ref{sec:Ru:cnst}.
The extra superconformal generator must admit a description in terms
of invariants for the monomial group. We discuss these monomial
invariants, and conclude the construction in
\S\ref{sec:Ru:mnml:monomialinvs}.

\subsubsection{Cayley algebra}\label{sec:Ru:mnml:cayleyalg}

Let $\Pi={\rm PG}(1,7)=\{\infty,0,1,2,3,4,5,6\}$ be a copy of the
projective line over $\FF_7$. The group $L_2(7)$ acts doubly
transitively on $\Pi$ by permutations. Let $\gt{h}$ be a real vector
space of dimension $8$ with positive definite bilinear form
$\lab\cdot\,,\cdot\rab$, and let $\{h_i\}_{i\in\Pi}$ be an
orthonormal basis for $\gt{h}$ indexed by the set $\Pi$. The
$\ZZ$-lattice in $\gt{h}$ generated by the vectors of the form
$h_i\pm h_j$ is a lattice of type $D_8$, and adding the vector
$\tfrac{1}{2}\sum_ih_i$ to these generators, we obtain a copy of the
$E_8$ lattice, the unique up to isomorphism self-dual even lattice
of rank $8$, and we denote it by $\LL$.

A quadratic form $q$ on $\LL$ is defined by setting $q(\lambda)=\lab
\lambda,\lambda\rab/2$. Since $\LL$ is even, $q$ takes values in
$\ZZ$ and maps $2\LL$ to $4\ZZ$. Thus $q$ induces an $\ZZ/2=\FF_2$
valued quadratic form $\bar{q}$ on the quotient group
$\bar{\LL}=\LL/2\LL$.

\medskip

Let us set $\uu=\tfrac{1}{2}\sum_ih_i\in\LL$. Then $\LL$ supports a
structure of non-associative $\ZZ$-algebra for which $\uu$ is a
unit; this algebra is known as the integral Cayley algebra, and is
isomorphic to a maximal integral order in the Octonion algebra. The
algebra structure may be defined in the following way \cite{ATLAS}.
We impose the relations
\begin{align}
    &2h_n^2=h_n-\uu,\\
    &2h_{\infty}h_0=\uu-h_3-h_5-h_6,\\
    &2h_{0}h_{\infty}=\uu-h_2-h_1-h_4,
\end{align}
and the images of these under the action of $L_2(7)$.

The subspace $2\LL$ is an ideal for this algebra structure, and we
thus obtain a Cayley algebra over $\FF_2$ on the space
$\bar{\LL}=\LL/2\LL$. We write $\lambda\equiv\lambda'$ when
$\lambda+2\LL=\lambda'+2\LL$ in $\bar{\LL}$. Note that
$2h_i\in\LL\setminus 2\LL$ so that $2h_i$ is not zero in
$\bar{\LL}$. On the other hand, $h_i-h_j\in\LL$ so that $2h_i\cong
2h_j$ for all $i,j\in\Pi$. The full automorphism group of $\LL$ is
$W_{E_8}$, the Weyl group of $E_8$. The center of $W_{E_8}$ acts
trivially on $\bar{\LL}$, and the subgroup of $W_{E_8}/\Kumgp$
preserving the $\FF_2$ Cayley algebra structure on $\bar{\LL}$ is
isomorphic to the exceptional group of Lie type $G_2(2)$. This group
contains the simple group $G_2(2)'\cong U_3(3)$ with index two.

\medskip

Reducing the bilinear form on $\LL$ modulo $2\ZZ$ we obtain an
$\FF_2$ valued alternating bilinear form ${b}$ on $\bar{\LL}$. We
then have ${b}(x,y) =\bar{q}(x+y) +\bar{q}(x) +\bar{q}(y)$ for $x,y
\in\bar{\LL}$. One can check the following
\begin{prop}\label{prop:Ru:mnml:cayleyalg:isoimpliescomm}
We have $xy=yx$ in $\bar{\LL}$ just when $b(x,y)=0$.
\end{prop}

We have the following computations in the integral Cayley algebra
${\LL}$.
\begin{gather}
\begin{split}\label{eqn:Ru:mnml:cayleyalg:cubecomp}
    (h_{\infty}+h_0)^2&=\tfrac{1}{2}(h_{\infty}^2
    +h_{\infty}h_0+h_0h_{\infty}+h_0^2)\\
        &=\tfrac{1}{2}(h_{\infty}-\uu
        +\uu-h_3-h_5-h_6+\uu-h_2-h_1-h_4+h_0-\uu)\\
        &=-\uu+(h_{\infty}+h_0)
\end{split}\\
\begin{split}\label{eqn:Ru:mnml:cayleyalg:invlcomp}
    (h_{\infty}-h_0)^2&=\tfrac{1}{2}(h_{\infty}^2
    -h_{\infty}h_0-h_0h_{\infty}+h_0^2)\\
        &=\tfrac{1}{2}(h_{\infty}-\uu
        -\uu+h_3+h_5+h_6-\uu+h_2+h_1+h_4+h_0-\uu)\\
        &=-\uu
\end{split}
\end{gather}
The first computation (\ref{eqn:Ru:mnml:cayleyalg:cubecomp}) shows
that $(h_i+h_j)$ has order three in $\bar{\LL}$ and inverse
$\uu-(h_i+h_j)$ for all $i\neq j\in\Pi$. We obtain $28$ pairs
$\{a,a^{-1}\}$ of order three elements in $\bar{\LL}$ in this way.
When $a$ is a cube root of unity we refer to a pair of the form
$\{a,a^{-1}\}$ as a {\em cube root pair}, and we write $\Delta$ for
the set of all cube root pairs in the $\FF_2$ Cayley algebra
$\bar{\LL}$. The map $x\mapsto axa^{-1}$ is an algebra automorphism
of $\bar{\LL}$ whenever $a\in\bar{\LL}$ has order three.

From the second computation (\ref{eqn:Ru:mnml:cayleyalg:invlcomp})
we see that $(h_i-h_j)^2$ is an involution in $\bar{\LL}$, and
$(\uu+h_i-h_j)^2\equiv 0$ for all $i\neq j\in\Pi$. There are $28$
involutions of this form, but we obtain another $35$ involutions
like $\uu-(h_i+h_j+h_k+h_l)$. Thus we have $28+35=63$ pairs
$\{x,y\}\neq\{0,\uu\}$ in $\bar{\LL}$ such that $x+y\equiv \uu$ and
one element of the pair is an involution. Note that if $x$ is an
involution in $\bar{\LL}$ and $y\equiv\uu+x$ then $y^2\equiv 0$. In
each case the involution in the pair is distinguished by having
non-trivial norm and by being orthogonal to $\uu$.

The group $G_2(2)$ acts transitively on the $63$ involutions, and
also on the $28$ cube root pairs $\Delta$.

\medskip

If $x,y$ are involutions in $\bar{\LL}$ then $x+y+\uu$ is also an
involution just when $b(x,y)=0$. In this case we have
$x+y+z=\uu=x^2=y^2=z^2$ for $z=x+y+\uu$, and also that $b$ vanishes
on $\{\uu,x,y,z\}$. We refer to such a triple $\{x,y,z\}$ as an {\em
isotropic line} in $\bar{\LL}$. A minimal representative in $\LL$
for an involution in $\bar{\LL}$ either takes the form
$(1\bar{1}0^6)$ or $\tfrac{1}{2}(1^4\bar{1}^4)$, so that any
involution is orthogonal to precisely $30$ other involutions, and
there are thus $63\!\cdot\!30/6=315$ isotropic lines in $\bar{\LL}$.

Suppose that $\{x,y,z\}$ is an isotropic line in $\bar{\LL}$ and
that $xy=z$. Then $\{\uu,x,y,z\}$ is a multiplicative four group in
$\bar{\LL}$. Indeed, commutativity follows from isotropy by
Proposition~\ref{prop:Ru:mnml:cayleyalg:isoimpliescomm}, and from
$xy=z$ we deduce
\begin{gather}
    yz=y(x+y+\uu)=z+\uu+y=x,\\
    zx=z(y+z+\uu)=x+\uu+z=y.
\end{gather}
In this case we say that $\{x,y,z\}$ is an {\em isotropic ring} in
$\bar{\LL}$. Each involution belongs to just $3$ isotropic four
groups in $\bar{\LL}$, so just $63\!\cdot \!3/3=63$ of the isotropic
lines are isotropic rings.

For any given involution $x$ say there are exactly $24$ cube roots
of unity that are not orthogonal to $x$. We have
$b(x,a)=b(x,a^{-1})$ so that these $24$ cube roots constitute $12$
cube root pairs. We call the set of cube root pairs incident to an
involution $x$ the {\em dozen} associated to $x$. Let $l=\{x,y,z\}$
be an isotropic line. Then the intersection of the dozens associated
to the involutions in $l$ is a set of four inverse pairs of cube
roots of unity. We refer to the set of inverse pairs incident to an
isotropic line $l$ as the {\em quartet} associated to $l$. In the
case that $l$ is a ring, we say that the quartet associated to $l$
is a {\em ringed quartet}.

We refer to two distinct inverse pairs as a {\em couple}. Then any
couple belongs to five quartets, and thus determines five isotropic
lines. These five lines intersect in a single involution which in
turn determines a dozen. Thus we find that each couple determines a
unique dozen. On the other hand, each couple belongs to just one
ringed quartet. Six different couples determine the same dozen, just
as $\binom{4}{2}=6$ different couples determine the same ringed
quartet. The six couples corresponding to a given dozen constitute a
partition of the $12$ inverse pairs of that dozen into disjoint
couples.

\begin{rmk}
The combinatorial concepts that arise here from the action of
$G_2(2)$ on the $\FF_2$ Cayley algebra have direct analogues in the
Conway-Wales lattice \cite{ConRu}; a certain self-dual lattice of
rank $56$ whose automorphism group is a quadruple cover of the
Rudvalis group. We intentionally name these concepts so that there
is resonance with the nomenclature of \cite{ConRu}.
\end{rmk}

\subsubsection{Monomial group}\label{sec:Ru:mnml:monomialgp}

In order to specify subsets of $\Delta$ and so forth in a convenient
way, we follow \cite{WilRu} and arrange the elements of $\Delta$ in
to a $3\times 3$ grid where each block has three elements, and a
distinguished element ${\infty}$ is placed just below the center
bottom block of the $3\times 3$ grid. We suppose that $\Delta$ has
been enumerated $\Delta=\{1,\ldots,27,\infty\}$ and that the set
$\Delta$ is arranged in the following way, on what we from now on
refer to as the $\Delta$-grid.
\begin{gather}
\begin{tabular}{|c|c|c|}
    \hline
  \trip{1}{2}{3} & \trip{4}{5}{6} & \trip{7}{8}{9} \\ \hline
  \trip{10}{11}{12} & \trip{13}{14}{15} & \trip{16}{17}{18} \\ \hline
  \trip{19}{20}{21} & \trip{22}{23}{24} & \trip{25}{26}{27} \\ \hline
    \multicolumn{3}{c}{$\infty$}   \\ %\hline
\end{tabular}
\end{gather}
Now we may specify a subset of $\Delta$ by suitably placing
asterisks within a copy of the $\Delta$-grid. Each block has a
middle, left, and right element, and we sometimes write $M$, $L$, or
$R$, respectively, within a block as a short hand for
\begin{gather}
     \trip{\ast}{\cdot}{\cdot},\quad
     \trip{\cdot}{\ast}{\cdot},\quad
     \text{and}\;
     \trip{\cdot}{\cdot}{\ast},
\end{gather}
respectively. The notation $LM$, $LR$, $MR$, and $LMR$, within a
block is to be interpreted in a similar way. The symbol $\cdot$ will
denote either an empty position in the $\Delta$-grid, or an empty
block.

The arrangement of elements of $\Delta$ within the $\Delta$-grid may
be effected in such a way that the three elements within any block
are exactly the triples of inverse pairs that complete $\infty$ to a
ringed quartet. Choosing one of the blocks, the center block for
example, there are exactly three dozens that contain the ringed
quartet consisting of the pair $\infty$ and the pairs of that block.
These three dozens give us a partition of the $24$ elements outside
$\infty$ and the chosen block into three disjoint sets of size
eight. These eight sets necessarily each have one element from each
block outside the chosen one, and we may assume that the coordinates
of any one of these three dozens are all labeled by the same letter,
$M$, $L$ or $R$, so that the three dozens in question may be written
as follows.
\begin{gather}
\begin{tabular}{|c|c|c|}
  \hline
  $M$ & $M$ & $M$ \\ \hline
  $M$ & $LMR$ & $M$ \\ \hline
  $M$ & $M$ & $M$ \\ \hline
    \multicolumn{3}{c}{$\ast$}   \\ %\hline
\end{tabular}\quad
\begin{tabular}{|c|c|c|}
  \hline
  $L$ & $L$ & $L$ \\ \hline
  $L$ & $LMR$ & $L$ \\ \hline
  $L$ & $L$ & $L$ \\ \hline
    \multicolumn{3}{c}{$\ast$}   \\ %\hline
\end{tabular}\quad
\begin{tabular}{|c|c|c|}
    \hline
  $R$ & $R$ & $R$ \\ \hline
  $R$ & $LMR$ & $R$ \\ \hline
  $R$ & $R$ & $R$ \\ \hline
    \multicolumn{3}{c}{$\ast$}   \\ %\hline
\end{tabular}
\end{gather}

\medskip

Let $\gt{r}$ be a complex vector space of dimension $28$ with
positive definite Hermitian form $(\,\cdot\,,\cdot)$ and an
orthonormal basis $\{a_i\}_{i\in\Delta}$ indexed by $\Delta$. Let
$\gt{r}^*$ be the dual space to $\gt{r}$ with dual basis
$\{a^*_i\}_{i\in\Delta}$, and set $\gt{s}=\gt{r}\oplus\gt{r}^*$. We
extend the Hermitian form on $\gt{r}$ in the natural way to
$\gt{s}$. Similar to \S\ref{sec:cliffalgs:Herm}, we equip $\gt{s}$
with the symmetric bilinear form denoted
$\langle\cdot\,,\cdot\rangle$ which is just the natural form arising
from the pairing between $\gt{r}$ and $\gt{r}^*$ scaled by a factor
of $1/2$. We also suppose that $\mc{E}=\{e_i,e_{i'}\}_{i\in\Delta}$
is the orthonormal basis, with respect to the bilinear form
$\lab\cdot\,,\cdot\rab$, on $\gt{s}$ determined as in
\S\ref{sec:cliffalgs:Herm}, and that $X$ is the subgroup of
$\Sp(\gt{s})$ generated by the expressions $\ii e_ie_{i'}$. Later in
\S\ref{sec:Ru:cnst} we will define the superspace underlying $\aru$
by setting
\begin{gather}
    \aru=A(\gt{s})^0\oplus A(\gt{s})^0_{\ogi,X}
\end{gather}
so that the VOA underlying $\aru$ coincides with that constructed
for the group $\SL_{28}(\CC)/\Kumgp$ in \S\ref{sec:GLnq}.

\medskip

We now commence the construction of the monomial group; a subgroup
of ${\SU}(\gt{r})$ of the shape $2^7.G_2(2)$.

Let $A$ denote the group of sign changes on coordinates
corresponding to dozens in $\Delta$ and their complements. That is,
for each dozen $D\subset\Delta$ there is an element $\epsilon_D\in
A$ such that $\epsilon_D$ acts as $-1$ on the $a_i$ such that $i\in
D$, and fixes the other basis vectors. Generators for the dozens
(the dozens and their complements constitute a doubly even code
$\mc{D}$ on $\Delta$) are given in Table
\ref{tab:Ru:mnml:monomialgp:gendozs}.
\begin{table}
     \centering
     \caption{Dozens generating $\mc{D}$}\label{tab:Ru:mnml:monomialgp:gendozs}
\begin{tabular}{|c|c|c|}
  \hline
  $M$ & $M$ & $M$ \\ \hline
  $M$ & $L$$M$$R$ & $M$ \\ \hline
  $M$ & $M$ & $M$ \\ \hline
    \multicolumn{3}{c}{$\ast$}   \\ %\hline
\end{tabular}\,
\begin{tabular}{|c|c|c|}
  \hline
  $L$ & $L$ & $L$ \\ \hline
  $L$ & $L$$M$$R$ & $L$ \\ \hline
  $L$ & $L$ & $L$ \\ \hline
    \multicolumn{3}{c}{$\ast$}   \\ %\hline
\end{tabular}\,
\begin{tabular}{|c|c|c|}
  \hline
  $M$$R$ & $\cdot$ & $L$$M$ \\ \hline
  $L$$R$ & $\cdot$ & $L$$R$ \\ \hline
  $L$$M$ & $\cdot$ & $M$$R$ \\ \hline
   \multicolumn{3}{c}{$\cdot$}    \\ %\hline
\end{tabular}\,
\begin{tabular}{|c|c|c|}
  \hline
  $L$$M$ & $\cdot$ & $L$$R$ \\ \hline
  $M$$R$ & $\cdot$ & $M$$R$ \\ \hline
  $L$$R$ & $\cdot$ & $L$$M$ \\ \hline
   \multicolumn{3}{c}{$\cdot$}    \\ %\hline
\end{tabular}\,
\begin{tabular}{|c|c|c|}
  \hline
  $L$$R$ & $L$$M$ & $\cdot$  \\ \hline
  $L$$R$ & $L$$R$ & $\cdot$  \\ \hline
  $L$$R$ & $M$$R$ & $\cdot$  \\ \hline
   \multicolumn{3}{c}{$\cdot$}    \\ %\hline
\end{tabular}\,
\begin{tabular}{|c|c|c|}
  \hline
  $M$$R$ & $L$$R$ & $\cdot$  \\ \hline
  $M$$R$ & $M$$R$ & $\cdot$  \\ \hline
  $M$$R$ & $L$$M$ & $\cdot$  \\ \hline
   \multicolumn{3}{c}{$\cdot$}    \\ %\hline
\end{tabular}
\end{table}
These involutions together with the symmetry which is $-1$ on all
coordinates generate a group $A$ of the shape $2^7$.

Let $P$ denote the group consisting of just those coordinate
permutations corresponding to the elements of $G_2(2)$ that fix the
cube root pair labeled $\infty$. Then $P$ has the shape
$3^{1+2}.Q_8.2$, and $A$ and $P$ together generate a split extension
$2^7:3^{1+2}.Q_8.2$. Conjugation by an element $a$ of the $\FF_2$
Cayley algebra $\bar{\LL}$ is an automorphism of $\bar{\LL}$ when
$a$ is a cube root of unity. Let $a_{\infty}$ be one of the cube
roots in the pair labeled $\infty$. Then conjugation by $a_{\infty}$
fixes $\infty$ and stabilizes each block, and we may assume that it
permutes the points within each block as $(MLR)$. In keeping with
the notation of \cite{ConRu} and \cite{WilRu} we denote this
permutation by $Q$. Remaining generators are given by the
permutations
\begin{gather}
     \begin{split}
     N_0=&\;(1,19,25,7)(2,20,26,8)(3,21,27,9)\\
          &\qquad(4,10,22,16)(5,11,23,17)(6,12,24,18),
     \end{split}\\     %#N0
     \begin{split}
     N_{356}=&\;(1,22,25,4)(2,23,26,5)(3,24,27,6)\\
          &\qquad(7,16,19,10)(8,17,20,11)(9,18,21,12),\\     %#N356
     \end{split}\\
     \begin{split}
     F_{03}=&\;(1,19)(4,22)(7,25)(2,21)(3,20)(5,24)\\
          &\qquad(6,23)(8,27)(9,26)(11,12)(14,15)(17,18),
     \end{split}\\        %#F03
     \begin{split}
     V=&\;(19,11,3)(20,12,1)(21,10,2)
          (22,13,4)(23,14,5)\\
          &\qquad(24,15,6)
          (25,18,8)(26,16,9)(27,17,7),
     \end{split}               %#V
\end{gather}
which again are labeled in keeping with the notation of \cite{ConRu}
and \cite{WilRu}. Note that $N_0$ and $N_{356}$ are two order four
elements that fix the center block $\{13,14,15\}$, and generate a
$Q_8$ subgroup of $P$. Adjoining $F_{03}$ we obtain a $Q_8.2$
subgroup preserving but not fixing the center block.

To extend the group $A\tcolon P$ to a group of the shape
$2^7.G_2(2)$ we include any symmetry obtained via the following
general method. First pick a ringed quartet containing $\infty$, and
call it the chosen quartet. For example, we may take the quartet
containing $\infty$ and the points of the center block. A given
quartet is contained in exactly three dozens, and we choose one of
those that contains the chosen quartet. For example, we may take the
dozen containing the chosen quartet and having $M$'s in all blocks
other than that contained in the chosen quartet. There is a unique
partition of any dozen into three disjoint ringed quartets, and
there is a unique partition of the points of a dozen into six
couples in such a way that the union of any two couples is a
quartet. In particular the coupling on the chosen dozen refines each
of the three ringed quartets it contains. We now choose one of these
ringed quartets that is not the already chosen one, and call it the
fixed quartet. We call the remaining quartet the stabilized quartet.
We begin to define an element $m'\in {\SU}(\gt{r})$ by decreeing
that $m'$ fix the points $e_i$ for $i$ in the fixed quartet. We
decree also that $m'$ transpose $e_i$ with $e_{i'}$ for $\{i,i'\}$ a
couple in either of the chosen or stabilized ringed quartets, except
for the couple $\{\infty,\infty'\}$ containing $\infty$, for which
we insist that $m'$ transpose $e_{\infty}$ with $-e_{\infty'}$. Now
consider the dozens that contain the fixed quartet. There are just
three: the chosen dozen and two others, and we call these two the
fixed dozens. The points of the fixed dozens complementary to those
of the fixed quartet exhaust the $16$ points of $\Delta$ outside of
the chosen dozen, and the couplings on each of the fixed dozens
yield couplings on these $16$ points. These couplings have the
property that exactly one point of each couple is contained in a
block containing a point of the fixed quartet, and the other point
lies in a block containing a point of the stabilized quartet. We ask
now that $m'$ transpose $e_i$ with $f_{i'}$ whenever $\{i,i'\}$ is a
couple outside of the chosen dozen and $i$ is in the fixed quartet.
Then the action of $m'$ on all of $\gt{r}$ is determined once we
decree that $m'$ commute with multiplication by $\ii$. Taking
\begin{gather}
\begin{tabular}{|c|c|c|}
  \hline
  $\cdot$ & $\cdot$ & $\cdot$ \\ \hline
  $\cdot$ & $LMR$ & $\cdot$ \\ \hline
  $\cdot$ & $\cdot$ & $\cdot$ \\ \hline
    \multicolumn{3}{c}{$\ast$}   \\ %\hline
\end{tabular}\quad
\begin{tabular}{|c|c|c|}
  \hline
  $M$ & $M$ & $M$ \\ \hline
  $M$ & $LMR$ & $M$ \\ \hline
  $M$ & $M$ & $M$ \\ \hline
    \multicolumn{3}{c}{$\ast$}   \\ %\hline
\end{tabular}\quad
\begin{tabular}{|c|c|c|}
    \hline
  $M$ & $\cdot$ & $M$ \\ \hline
  $\cdot$ & $\cdot$ & $\cdot$ \\ \hline
  $M$ & $\cdot$ & $M$ \\ \hline
    \multicolumn{3}{c}{$\cdot$}   \\ %\hline
\end{tabular}
\end{gather}
to be the chosen quartet, the chosen dozen, and the fixed quartet,
respectively, we obtain $m'=m$ where $m$ is described in the
following way as a coordinate permutation followed by scalar
multiplications on coordinates.
\begin{gather}
\begin{tabular}{|c|c|c|}
                            \hline
  \trip{1}{\ii}{\ii} & \trip{1}{-\ii}{-\ii} & \trip{1}{\ii}{\ii}
  \\ \hline
  \trip{1}{-\ii}{-\ii} & \trip{-1}{1}{1} & \trip{1}{-\ii}{-\ii}
  \\ \hline
  \trip{1}{\ii}{\ii} & \trip{1}{-\ii}{-\ii} & \trip{1}{\ii}{\ii}     \\ \hline
      \multicolumn{3}{c}{$-1$}          \\ %\hline
\end{tabular}\;\circ
     %\\
\begin{array}{c}
     (2,5)(3,12)(4,22)(6,9)\\
     (8,17)(10,16)(11,20)(13,\infty)\\
     (14,15)(18,27)(21,24)(23,26)
\end{array}
\end{gather}
The symmetry $m$ is denoted $M_3E_{01}$ in \cite{WilRu}.

We set $M$ to be the subgroup of ${\SU}(\gt{r})$ generated by $A$,
$P$ and $m$, and $M$ is then a non-split extension of the form
$2^7.G_2(2)$. When written with respect to the basis $\{a_i\}$, the
matrices representing the action of $M$ on $\gt{r}$ all have the
property that any row or column has exactly one entry, and that
entry lies in $\{\pm 1,\pm \ii\}$. We call such a matrix monomial,
and we call $M$ the monomial group. By replacing the non-zero
entries of the monomial matrices representing $M$ with $1$'s, each
monomial matrix becomes a permutation matrix, and we obtain a
homomorphism of groups $M\to \overline{M}$ where $\overline{M}$ acts
as permutations on the $28$ coordinates $a_i$. The kernel of this
homomorphism is the $2^7$ subgroup of $M$, and $\overline{M}$ is
isomorphic to $G_2(2)$. We write $\bar{m}$ for the image of $m\in M$
under the map $M\to\overline{M}$.

\subsubsection{Monomial invariants}\label{sec:Ru:mnml:monomialinvs}

In this section we consider the invariants of the monomial group $M$
in the space $\Cm(\gt{s})_X^0$ (see \S\ref{sec:cliffalgs:mods}). In
\S\ref{sec:Ru:cnst} we will find that this space is just the
subspace of $\aru$ with degree $7/2$.

\medskip

Recall that the set $\Delta$ has been enumerated in
\S\ref{sec:Ru:mnml:monomialgp}, and let us agree to write $e_I$ for
$e_{i_1}\cdots e_{i_k}$ whenever $I=\{i_1,\ldots,i_k\}\subset
\Delta$ and $i_1<\cdots<i_k$. (We decree that $i<\infty$ for all
$i\in\Delta\setminus\infty$.) An orthonormal basis for
$\Cm(\gt{s})_X$ is then given by the set $\{e_I{1}_X\}$ where $I$
ranges over all subsets of $\Delta$.

Suppose that $t$ is an $M$ invariant vector in $\Cm(\gt{s})_X$ and
that $t$ has non-zero projection onto $e_I$ for some
$I\subset\Delta$. Then $I$ must lie in $\mc{D}^{\circ}$, the dual to
the dozens code $\mc{D}$, since otherwise some element of $A$ would
negate $e_I$. Indeed, $t$ must have non-zero projection also onto
$e_{\bar{m}(I)}$ for every $\bar{m}\in\overline{M}$ since
$m(e_I1_X)\in \CC e_{\bar{m}(I)}1_X$ for any $m\in M$ with image
$\bar{m}$ in $\overline{M}$. (In fact one can show that
$m(e_I1_X)\in \{\pm e_{\bar{m}(I)}1_X\}$ for all $m\in M$ and $I\in
\mc{D}^{\circ}$.) In this way we see that $M$-invariant vectors may
be obtained by specifying an orbit $\mc{O}$ of $\overline{M}\cong
G_2(2)$ on $\mc{D}^{\circ}$ and a map $\gamma:\mc{O}\to\CC$ such
that $m(\gamma_{I}e_I1_X)=\gamma_{\bar{m}(I)}e_{\bar{m}(I)}1_X$ for
$m\in M$ and $I\in\mc{O}$, and any $M$-invariant vector in
$\Cm(\gt{s})_X$ must be a linear combination of vectors obtained in
this way from orbits of $\overline{M}$ in $\mc{D}^{\circ}$.

It turns out that a suitable map $\gamma$ exists only for certain
orbits of $\overline{M}$ in $\mc{D}^{\circ}$. For example, explicit
calculation reveals that there are $80$ orbits of $\overline{M}$ on
words of weight $14$ in $\mc{D}^{\circ}$, and just $68$ of these
orbits give rise to $M$-invariant vectors in the weight $14$
subspace of $\Cm(\gt{s})_X$. We collect these $68$ orbits on weight
$14$ words in $\mc{D}^{\circ}$ for which suitable maps $\gamma$
exist into a set $\gt{O}=\{\mc{O}\}$, so that the vectors
$t_{\mc{O}}=\sum_{I\in\mc{O}}\gamma_{I}e_I1_X$ for $\mc{O}$ in
$\gt{O}$ span the point-wise $M$-invariant subspace of
$\bigwedge^{14}(\gt{r})\subset\Cm(\gt{s})_{X}$. We refer to the
orbits of $\overline{M}$ in $\gt{O}$ as the {\em relevant orbits}.

The functions $\gamma$ and the vectors $t_{\mc{O}}$ are determined
only up to scalar factors, and for a given $I$ in a relevant orbit
$\mc{O}$, we may assume if we wish that $\gamma$ is normalized so
that $\gamma_{I}=1$. Indeed, the $\overline{M}$ invariant vector
$t_{\mc{O}}$ is completely determined once we specify an $I$ in
$\mc{O}$, and insist that $\gamma_I=1$.
\begin{table}
  \centering
  \caption{Data for $\spsp$}
  \label{tab:Ru:mnml:monomialinvs:rhodata}
%\begin{small}
%\begin{footnotesize}
  \begin{tabular}{c|ll}
     %\hline
     $\overline{M}$ orbit representatives & $r_{\mc{O}}$ & $r_{\mc{O}'}$\\
    \hline
    $\{3, 5, 6, 8, 11, 18, 19, 20, 21, 22, 25, 26, 27, \infty\}$  & $4-3\ii$ & $4+3\ii$\\
    $\{2, 3, 4, 6, 7, 8, 10, 12, 13, 14, 17, 18, 21, \infty\}$ & $-4\ii$ & $4\ii$\\
    $\{2, 4, 5, 6, 7, 8, 9, 13, 14, 15, 16, 22, 23, \infty\}$  & $-2$ &
    $2$\\
    $\{1, 5, 6, 10, 12, 13, 14, 15, 16, 17, 20, 22, 25, 27\}$ & $-2$ &
    $-2$\\
    $\{1, 3, 9, 10, 11, 12, 18, 19, 20, 21, 24, 26, 27, \infty\}$ & $-2\ii$ &
    $2\ii$\\
    $\{1, 3, 7, 8, 9, 10, 12, 17, 19, 20, 21, 25, 27, \infty\}$  & $2\ii$ &
    $-2\ii$\\
    $\{2, 3, 4, 5, 6, 13, 14, 15, 16, 17, 19, 21, 22, 25\}$ & $-2+4\ii$ &
    $2+4\ii$\\
    $\{1, 2, 4, 5, 7, 8, 9, 10, 11, 12, 13, 18, 21, \infty\}$ & $2$ &
    $2$\\
    $\{4, 6, 8, 10, 13, 15, 18, 20, 21, 22, 24, 26, 27, \infty\}$ & $-\ii$ & $-\ii$\\
    $\{3, 4, 7, 8, 9, 10, 11, 12, 13, 17, 19, 23, 25, 26\}$ & $-\ii$ & $\ii$\\
    $\{1, 2, 3, 5, 6, 9, 12, 13, 15, 16, 20, 22, 24, \infty\}$ & $\ii$ & $\ii$\\
    $\{1, 3, 5, 6, 7, 9, 11, 12, 13, 15, 20, 22, 24, \infty\}$ & $2+\ii$ & $-2+\ii$\\
    $\{1, 3, 5, 6, 9, 12, 14, 15, 17, 19, 20, 22, 24, \infty\}$ & $-2-\ii$ & $-2+\ii$\\
    $\{1, 2, 5, 6, 7, 11, 13, 14, 15, 16, 20, 22, 24, \infty\}$ & $2\ii$ & $2\ii$\\
    $\{1, 2, 3, 5, 7, 9, 11, 13, 14, 15, 16, 17, 19, \infty\}$ & $-\ii$ & $-\ii$\\
    $\{1, 4, 5, 8, 10, 12, 13, 20, 21, 22, 23, 25, 26, 27\}$ & $-2$ & $-2$\\
    $\{1, 3, 5, 6, 9, 10, 13, 17, 18, 19, 21, 24, 26, 27\}$  & $2\ii$ & $-2\ii$\\
    $\{4, 8, 10, 12, 15, 18, 20, 21, 22, 23, 25, 26, 27, \infty\}$ & $\ii$ & $-\ii$\\
    $\{1, 3, 4, 5, 7, 8, 9, 11, 13, 14, 17, 19, 21, 27\}$ & $2$ & $2$\\
    $\{8, 9, 14, 15, 16, 17, 18, 19, 21, 23, 24, 25, 27, \infty\}$ & $-3-2\ii$ & $3-2\ii$\\
    $\{8, 9, 14, 15, 16, 17, 18, 19, 20, 22, 23, 26, 27, \infty\}$ & $-3-\ii$ & $3-\ii$\\
    $\{8, 9, 14, 15, 16, 17, 18, 19, 20, 21, 22, 24, 27, \infty\}$  & $-2-\ii$ & $2-\ii$\\
    $\{8, 9, 14, 16, 18, 19, 20, 21, 22, 23, 25, 26, 27, \infty\}$ & $-2$ & $2$\\
    $\{9, 13, 14, 16, 17, 18, 19, 20, 21, 23, 24, 25, 27, \infty\}$ & $-1+2\ii$ & $1+2\ii$\\
    $\{9, 13, 14, 15, 16, 18, 19, 21, 22, 23, 25, 26, 27, \infty\}$ & $1+2\ii$ & $-1+2\ii$\\
    $\{9, 13, 14, 15, 16, 18, 19, 20, 21, 22, 23, 24, 26, 27\}$ & $-1+\ii$ & $-1-\ii$\\
    $\{8, 9, 13, 14, 15, 18, 19, 20, 23, 24, 25, 26, 27, \infty\}$ & $-1$ & $-1$\\
    $\{8, 9, 13, 14, 15, 18, 19, 20, 21, 22, 23, 24, 26, \infty\}$ & $-1+\ii$ & $1+\ii$\\
    $\{8, 9, 13, 15, 16, 18, 19, 20, 21, 22, 23, 25, 26, 27\}$ & $1-\ii$ & $-1-\ii$\\
    $\{9, 12, 14, 15, 16, 17, 18, 20, 21, 23, 24, 25, 26, 27\}$ & $-1-\ii$ & $-1+\ii$\\
    $\{9, 12, 14, 15, 16, 17, 18, 19, 20, 21, 22, 23, 25, 26\}$ & $1-\ii$ & $1+\ii$\\
    $\{9, 12, 15, 17, 18, 19, 20, 21, 22, 23, 24, 26, 27, \infty\}$ & $1+\ii$ & $1-\ii$\\
    $\{8, 9, 12, 14, 16, 17, 20, 21, 22, 23, 24, 25, 26, 27\}$ & $1+3\ii$ & $-1+3\ii$\\
    $\{8, 9, 12, 14, 15, 16, 19, 20, 21, 24, 25, 26, 27, \infty\}$ & $2-2\ii$ & $-2-2\ii$\\
    %\hline
  \end{tabular}
%\end{small}
%\end{footnotesize}
\end{table}

If $\mc{O}$ is an orbit of $\overline{M}$ on weight $14$ words, then
$\mc{O}'=\{I'=\Delta\setminus I\mid I\in\mc{O}\}$ is also an orbit
for $\overline{M}$ on weight $14$ words, and the $68$ relevant
orbits $\gt{O}$ constitute $34$ pairs of the form
$\{\mc{O},\mc{O}'\}$. In the first column of Table
\ref{tab:Ru:mnml:monomialinvs:rhodata} we provide a list of $34$
weight $14$ subsets of $\Delta$, chosen so that the corresponding
orbits under $\overline{M}$ constitute one orbit each from the $34$
pairs $\{\mc{O},\mc{O}'\}$ in $\gt{O}$. The corresponding
$\overline{M}$ invariants $t_{\mc{O}}$ and $t_{\mc{O}'}$ are
determined once we insist that $\gamma_I=\gamma_{I'}=1$ for each $I$
in Table \ref{tab:Ru:mnml:monomialinvs:rhodata}, where
$I'=\Delta\setminus I$.

From the $68$ dimensional space spanned by the $t_{\mc{O}}$ and
$t_{\mc{O}'}$ we require to pick a single vector that will allow us
to realize the Rudvalis group. We chose $\spsp$ by setting
\begin{gather}
     \spsp=\frac{1}{\sqrt{C}}
          \sum r_{\mc{O}}t_{\mc{O}}+r_{\mc{O}'}t_{\mc{O}'}
\end{gather}
where the vectors $t_{\mc{O}}$ and $t_{\mc{O}'}$ are determined by
the entries in the first column of Table
\ref{tab:Ru:mnml:monomialinvs:rhodata} as explained above, the
coefficients $r_{\mc{O}}$ and $r_{\mc{O}'}$ are given in the second
column of Table \ref{tab:Ru:mnml:monomialinvs:rhodata}, and the
constant $C$ is given by $C=86272=2^8.337$.

Let $F$ be the subgroup of ${\SL}(\gt{r})$ that fixes $\spsp$, and
recall that $\gt{z}$ denotes the element $e_{\Delta}e_{\Delta'}\in
{\SL}(\gt{r})$, with the ordering chosen so that $\gt{z}$ lies in
$F$. Then $F$ contains the monomial group $M$. We have the following
\begin{thm}\label{thm:Ru:mnml:monomialinvs:FhasR}
The group $F$ contains a group isomorphic to $\cRu$.
\end{thm}
Recall that the monomial group $M$ is a maximal subgroup of $\cRu$
so we require just one more generator fixing $\spsp$ in order to
extend $M$ to the desired group. Following \cite{WilRu} we present
here an element $z$ (denoted there by $E_{24}$) that satisfies our
objective, and whose action is determined as follows. The element
$z$ acts with order two, fixing a $24$ dimensional space in $\gt{r}$
and negating the complementary space. The fixed $24$ dimensional
space is spanned by the vectors of Table \ref{tab:Ru:symms:zfixed}
along with their images under multiplication by $\ii$ and the
$Q_8.2$ group of coordinate permutations stabilizing the center
block. (Recall from \S\ref{sec:Ru:mnml:monomialgp} that this $Q_8.2$
group is generated by $N_0$, $N_{356}$ and $F_{03}$.)
\begin{table}
     \centering
     \caption{Vectors fixed by $z$}\label{tab:Ru:symms:zfixed}
\begin{footnotesize}
\begin{tabular}{|c|c|c|}\hline
  \trip{\cdot}{\cdot}{\cdot} & \trip{\cdot}{\cdot}{\cdot}
     & \trip{\cdot}{\cdot}{\cdot} \\ \hline
  \trip{\cdot}{\cdot}{\cdot} & \trip{1}{\cdot}{\cdot}
     & \trip{\cdot}{\cdot}{\cdot} \\ \hline
  \trip{\cdot}{\cdot}{\cdot} & \trip{\cdot}{\cdot}{\cdot}
     & \trip{\cdot}{\cdot}{\cdot} \\ \hline
     \multicolumn{3}{c}{$1$} %\hline
\end{tabular}\;
\begin{tabular}{|c|c|c|}\hline
  \trip{\cdot}{\cdot}{\cdot} & \trip{\cdot}{\cdot}{\cdot}
     & \trip{\cdot}{\cdot}{\cdot} \\ \hline
  \trip{\cdot}{\cdot}{\cdot} & $\trip{-1}{\ii}{\ii}$
     & \trip{\cdot}{\cdot}{\cdot}\\ \hline
  \trip{\cdot}{\cdot}{\cdot} & \trip{\cdot}{\cdot}{\cdot}
     & \trip{\cdot}{\cdot}{\cdot} \\ \hline
     \multicolumn{3}{c}{$1$} %\hline
\end{tabular}\;
\begin{tabular}{|c|c|c|}\hline
  \trip{\cdot}{\cdot}{\cdot} & \trip{1}{\cdot}{\cdot}
     & \trip{\cdot}{\cdot}{\cdot} \\ \hline
  \trip{1}{\cdot}{\cdot} & \trip{\cdot}{\cdot}{\cdot}
     & \trip{1}{\cdot}{\cdot} \\ \hline
  \trip{\cdot}{\cdot}{\cdot} & \trip{1}{\cdot}{\cdot}
     & \trip{\cdot}{\cdot}{\cdot} \\ \hline
     \multicolumn{3}{c}{$\cdot$} %\hline
\end{tabular}\;
\begin{tabular}{|c|c|c|}\hline
  \trip{\cdot}{\cdot}{\cdot} & \trip{\cdot}{\cdot}{\cdot}
     & \trip{\cdot}{\cdot}{1} \\ \hline
  \trip{\cdot}{\cdot}{-1} & \trip{\cdot}{\cdot}{\cdot}
     & \trip{\cdot}{\cdot}{-1}\\ \hline
  \trip{\cdot}{\cdot}{1} & \trip{\cdot}{\cdot}{\cdot}
     & \trip{\cdot}{\cdot}{\cdot}\\ \hline
     \multicolumn{3}{c}{$\cdot$} \\ %\hline
\end{tabular}\;
\begin{tabular}{|c|c|c|}\hline
  \trip{\cdot}{\cdot}{\cdot} & \trip{\cdot}{\cdot}{\ii}
     & \trip{\cdot}{\cdot}{\cdot}\\ \hline
  \trip{\cdot}{\cdot}{-1} & \trip{\cdot}{\cdot}{\cdot}
     & \trip{\cdot}{\cdot}{1} \\ \hline
  \trip{\cdot}{\cdot}{\cdot} & \trip{\cdot}{\cdot}{-\ii}
     & \trip{\cdot}{\cdot}{\cdot} \\ \hline
     \multicolumn{3}{c}{$\cdot$} \\ %\hline
\end{tabular}
\end{footnotesize}
\end{table}
It is a computation to verify that the action of $z$ fixes $\spsp$.
Once again, the computer algebra package \cite{GAP4} is a very
useful tool for carrying out the necessary calculations.

\begin{rmk}
For the interested reader we remark that other elements like $z$ may
be defined in the following way. Similar to the case with $m$ in
\S\ref{sec:Ru:mnml:monomialgp} we choose a ringed quartet containing
$\infty$ and a dozen containing that quartet. We then take $N$ and
$N'$ to be two coordinate permutations of order $4$ generating the
$Q_8$ group in $H$ that fixes the chosen quartet. We may assume that
$N$ stabilizes the ringed quartets contained in the chosen dozen,
and that $N'$ interchanges the two quartets other than the chosen
quartet. We take $F$ to be an order two permutation extending the
$Q_8$ to $Q_8.2$, and stabilizing the chosen quartet. As above we
take $Q$ to be the order three coordinate permutation in $H$ that
stabilizes all the ringed quartets containing $\infty$. We now
define an element $z'\in \SO(\gt{r})$ by first decreeing that $z'$
fix the elements $e_{\infty}+e_{\infty'}$ and
$e_{\infty}-e_{\infty'}+f_x+f_{x'}$ if $\{\infty,\infty'\}$ and
$\{x,x'\}$ are the couples of the chosen dozen contained in the
chosen quartet. We then decree that $z'$ fix
$e_{i}+e_{i'}+e_{i''}+e_{i'''}$ when $\{i,i',i'',i'''\}$ is one of
the ringed quartets in the chosen dozen other than the chosen
quartet. We insist also that $z'$ fix
$Q(e_{i}+e_{i'})-N'Q(e_i+e_{i'})$ and
$Q(e_{i}-e_{i'})+NQ(f_i-f_{i'})$ when $\{i,i'\}$ is a couple of the
chosen dozen not in the chosen quartet. Finally, the action of $z'$
on $\gt{r}$ is determined completely once we decree that $z'$
commute with $I$, $N$, $N'$ and $F$. With the chosen quartet and the
chosen dozen given by
\begin{gather}
\begin{tabular}{|c|c|c|}
  \hline
  $\cdot$ & $\cdot$ & $\cdot$ \\ \hline
  $\cdot$ & $LMR$ & $\cdot$ \\ \hline
  $\cdot$ & $\cdot$ & $\cdot$ \\ \hline
    \multicolumn{3}{c}{$\ast$}    %\hline
\end{tabular}\quad\text{ and }\quad
\begin{tabular}{|c|c|c|}
  \hline
  $M$ & $M$ & $M$ \\ \hline
  $M$ & $LMR$ & $M$ \\ \hline
  $M$ & $M$ & $M$ \\ \hline
  \multicolumn{3}{c}{$\ast$} %\hline
\end{tabular}
\end{gather}
respectively, we may take $N=N_0$, $N'=N_{356}$, and $F=F_{03}$, and
we then obtain $z'=z$.
\end{rmk}

\subsection{Construction}\label{sec:Ru:cnst}

Let $\gt{r}$ be as in \S\ref{sec:Ru:geom} or \S\ref{sec:Ru:mnml} so
that $\gt{r}$ is a complex vector space of dimension $28$ with
positive definite Hermitian form $(\,\cdot\,,\cdot)$ and an
orthonormal basis $\{a_i\}_{i\in\Delta}$ indexed by $\Delta$. We
write $\gt{r}^*$ for the dual space to $\gt{r}$ with dual basis
$\{a^*_i\}_{i\in\Delta}$, and we set $\gt{s}=\gt{r}\oplus\gt{r}^*$.
We extend the Hermitian form on $\gt{r}$ in the natural way to
$\gt{s}$, and similarly to \S\ref{sec:cliffalgs:Herm}, we equip
$\gt{s}$ with the symmetric bilinear form denoted
$\langle\cdot\,,\cdot\rangle$ which is just the natural form arising
from the pairing between $\gt{r}$ and $\gt{r}^*$ scaled by a factor
of $1/2$. We also suppose that $\mc{E}=\{e_i,e_{i'}\}_{i\in\Delta}$
is the orthonormal basis, with respect to the bilinear form
$\lab\cdot\,,\cdot\rab$, on $\gt{s}$ determined as in
\S\ref{sec:cliffalgs:Herm}, and that $X$ is the subgroup of
$\Sp(\gt{s})$ generated by the expressions $\ii e_ie_{i'}$.

The space underlying $\aru$ is defined by setting
$\aru=A(\gt{s})^0\oplus A(\gt{s})^0_{\ogi,X}$. Then $\aru$ admits a
natural structure of VOA, as demonstrated in \S\ref{sec:GLnq:tw}.
Let us set $\cgset_{\Ru}=\{\cas,\vp,\vm,\spsp\}$, where $\spsp$ is
as in \S\ref{sec:Ru:geom} or \S\ref{sec:Ru:mnml}. (By character
table computations we know that there is a unique invariant for the
Rudvalis group in $\bigwedge^{14}(\gt{r})$, so the two constructions
can differ only up to a constant, and this does not effect the
enhanced conformal structure.) Recalling Proposition~
\ref{prop:GLnq:SLgps:vpmOPEs} and the remarks following, we then
have
\begin{thm}\label{thm:Ru:cnst:eVOA}
The quadruple $(\aru,Y,\vac,\cgset_{\Ru})$ is a self-dual enhanced
$\U(1)$--VOA of rank $28$.
\end{thm}

\subsection{Symmetries}\label{sec:Ru:symms}

We now consider the automorphism group of the enhanced VOA $\aru$
with enhanced conformal structure determined by $\cgset_{\Ru}$.

Let $F$ denote the subgroup of $\SL(\gt{r})$ that fixes $\spsp$. By
Theorem~\ref{thm:Ru:geom:FhasRu} or Theorem~
\ref{thm:Ru:mnml:monomialinvs:FhasR}, a group isomorphic to $\cRu$
is contained in $F$. We now show
\begin{prop}\label{prop:Ru:symms:Fisfinite}
The group $F$ is finite.
\end{prop}
\begin{proof}
The group $F$ is the stabilizer of a subspace of the natural module
for the algebraic group $\SL_{28}(\CC)$, and is thus also algebraic.
Since $F$ contains the group $R$ which acts irreducibly on the
natural module, we have that $F$ is reductive. Let $\gt{k}$ be the
Lie algebra of the connected component of the identity in $F$. Then
$\gt{k}$ embeds in the Lie algebra of
$\Aut\left(\aru,\{\cas,\vp,\vm\}\right)$ which we denote $\gt{g}$,
and which may be identified with the subspace of $(\aru)_{1}$
spanned by the $x_{ij}=a_i(-\tfrac{1}{2})a_j^*(-\tfrac{1}{2})\vac$
for $i\neq j$, and by the $y_{ij}=x_{ii}-x_{jj}$ for $i\neq j$ (c.f.
Proposition~\ref{prop:GLnq:SLgps:Autsjvpm}). Then for all
$x\in\gt{k}$ we have $\exp(x_{(0)})\spsp=\spsp$, and this implies
$x_{(0)}\spsp=0$ for some non-trivial $x\in\gt{k}$ so long as
$\gt{k}$ is not trivial. Consider now the image of $\gt{g}$ in
$\bigwedge^{14}(\gt{a})\vac_X\subset \atw{\gt{u}}_{7/2}$ under the
map $\phi:x\mapsto x_{(0)}\spsp$. The group $R$ acts irreducibly on
$\gt{g}$, so $\phi$ is either the zero map, or is an embedding of
$R$ modules. In the former case $\spsp$ is invariant for the action
of $\SL(\gt{a})$ on $\bigwedge^{14}(\gt{a})\vac_X$, but we know that
this group acts irreducibly on this space. In the latter case no
non-trivial element $x$ of $\gt{g}$ can satisfy $x_{(0)}\spsp=0$,
and thus we have that $\gt{k}$ is trivial. We conclude that
$\dim(F)=0$, and thus $F$ is finite.
\end{proof}

We complete the determination of $\Aut(\aru,\cgRu)$ by showing that
$R\cong\cRu$ is a maximal finite subgroup of $\SL_{28}(\CC)$ up to a
group of scalar matrices. The techniques we employ in order to
establish this result are borrowed directly from
\cite{NebRaiSloInvtsCliffGps}, including the following
\begin{lem}\label{lem:Ru:symms:psubscentral}
Suppose $X$ is a finite subgroup of $\SL(\gt{r})$ containing $R$.
Then any normal $p$-subgroup of $X$ is central.
\end{lem}
\begin{proof}
Observe that $R$ acts {\em primitively} on $\gt{r}$, meaning that
there is no decomposition
$\gt{r}=\gt{r}_1\oplus\cdots\oplus\gt{r}_k$ into (non-trivial)
subspaces that is preserved by $R$; similarly then for $X$. If $N$
is a normal subgroup of $X$, then $X$ permutes the isotypic
components of $\gt{r}|_N$ (which is $\gt{r}$ viewed as an $N$
module) and it follows that $\gt{r}|_N$ is isotypic (a direct sum of
copies of a single irreducible representation) for $N$ whenever $N$
is normal in $X$. Since an irreducible representation of an abelian
group is faithful only on a cyclic subgroup (and $X$ acts faithfully
on $\gt{r}$ by definition), any abelian normal subgroup of $X$ is
cyclic.

Let $N$ be a normal $p$-subgroup of $X$, and consider the group
$C=C_X(N$), the centralizer of $N$ in $X$. The intersection $C\cap
R$ is normal in $R$, and contains the center $Z(R)$ of $R$, and thus
$C\cap R$ is either $Z(R)$ or $R$. In the former case we obtain an
embedding $\Ru\hookrightarrow\Aut(N)$. In the latter case $N$ can
consist only of scalar matrices, since $R$ acts absolutely
irreducibly on $\gt{r}$. So we may assume that $N$ is a normal
$p$-subgroup such that $\Aut(N)$ admits an embedding of $\Ru$. An
irreducible representation of $N$ in $\gt{r}|_N$ has degree that is
a power of $p$ and divides $\dim(\gt{r})=28$. So if $p$ is not $2$
or $7$, then $N$ is abelian, and we have seen above that such a
group is cyclic. Evidently, this contradicts
$\Ru\hookrightarrow\Aut(N)$. For the remaining cases that $p=2$
(resp. $p=7$) and $N$ admits a faithful embedding in $\SL_4(\CC)$
(resp. $\SL_7(\CC)$), we consult Hall's classification of $p$-groups
of symplectic type (c.f. \cite[(29.3)]{AscFGT}). A $p$-group is said
to be of {\em symplectic type} if it has no non-cyclic
characteristic abelian subgroups. $N$ is evidently a group of
symplectic type since a characteristic abelian subgroup of $N$ would
be normal in $X$, and we have seen above that such a group is
cyclic. Hall's classification reveals that none of the remaining
cases for $N$ admit a non-trivial action by $\Ru$ as automorphisms.
This completes the proof.
\end{proof}

The proof of the following proposition owes a great deal to the
proof of Theorem~6.5 in \cite{NebRaiSloInvtsCliffGps}.
\begin{prop}\label{prop:Ru:symms:Rumaxfinite}
If $X$ is a finite subgroup of $\SL_{28}(\CC)$ containing $R$ then
$X=\lab R,\xi\Id\rab$ for $\xi$ a root of unity in $\CC$.
\end{prop}
\begin{proof}
Let $K$ be a minimal abelian number field containing $\ii$ such that
$X$ is conjugate to a subgroup of $\SL_{28}(K)$, and let $\gt{O}$ be
the ring of integers in $K$. The group $X$ must preserve a Hermitian
lattice in the natural module (which we denote $K^{28}$). Any
Hermitian $\ZZ[\ii]R$ lattice in $\QQ(\ii)^{28}$ is isometric to the
Conway--Wales lattice $\LL_{\Ru}$ (c.f. \cite{TieGIRReps}); it
follows that any $\gt{O}X$ lattice in $K^{28}$ is of the form
$I\otimes_{\ZZ[\ii]}\LL_{\Ru}$ for $I$ a fractional ideal of
$\gt{O}$, and if $X$ preserves any such lattice, it must also
preserve $\gt{O}\otimes_{\ZZ[\ii]}\LL_{\Ru}$. Using the action of
$X$ on $\gt{O}\otimes_{\ZZ[\ii]}\LL_{\Ru}$ we may regard $X$ as a
group of matrices with entries in $\gt{O}$, and we then may let the
Galois group $\Gamma=\Gal(K/\QQ(\ii))$ act on $X$, by acting
componentwise on the corresponding matrices in $\SL_{28}(\gt{O})$.

If $\Gamma$ is trivial then $K=\QQ(\ii)$ and $X=\lab R,\ii\Id\rab$,
and we are done. Otherwise, let $\gt{p}$ be a prime ideal of
$\gt{O}$ that ramifies in $K/\QQ(\ii)$, and let $\sigma$ be a
non-trivial element in the inertia group $\Gamma_{\gt{p}}$.
\begin{gather}
     \Gamma_{\gt{p}}=\{\sigma\in\Gamma\mid
     \sigma(a)\equiv a\pmod{\gt{p}},\;\forall a\in\gt{O}\}
\end{gather}
Observe that for arbitrary $g\in X$ we have that $g^{-1}\sigma(g)$
lies in the group $X_{\gt{p}}$ consisting of elements $g\in X$ such
that $g\equiv \Id\pmod{\gt{p}}$. The group $X_{\gt{p}}$ is a normal
$p$-subgroup of $X$, for $p$ the rational prime divisible by
$\gt{p}$. By Lemma \ref{lem:Ru:symms:psubscentral} such a group is
central, and we see that the map $g\mapsto g^{-1}\sigma(g)$ is a
homomorphism of $X$ into an abelian group. This shows that the
commutator subgroup $X^{(1)}$ of $X$ is fixed by $\sigma$. Now
$X^{(1)}$ is properly contained in $X$ for otherwise $K$ is not
minimal, and the argument thus far shows that any finite subgroup of
$\SL(\gt{a})$ containing $R$ either is realized over $\QQ(\ii)$ (and
is then contained in $\lab R,\ii\Id\rab$), or properly contains its
own commutator subgroup. Consider now the chain $X\geq X^{(1)}\geq
X^{(2)}\geq \cdots$. Since $X$ is finite and every term $X^{(k)}$
contains $R=R^{(1)}$, the sequence stabilizes and we have
$X^{(k)}=X^{(k+1)}$ for some $k$, and $R\leq X^{(k)}\leq\lab
R,\ii\Id\rab$. This shows that $R$ is normal in $X$, and we then
obtain a map $X\to\Aut(R)$ by letting $x\in X$ act by conjugation on
$R$. The action of $x$ on $R$ is trivial if and only if $x$ commutes
with $R$, and since $R$ acts absolutely irreducibly on $K^{28}$,
such an $x$ must be a scalar matrix in $\SL_{28}(K)$. Clearly, the
scalar matrices in $X$ constitute the center of $X$. Let us write
$T$ for the group of scalar matrices in $X$. Then we have an
injective map $X/T\to\Ru$, which is also surjective since $X$
contains $R$. We conclude that $X=TR$, as required.
\end{proof}
We are now ready to identify $\Aut(\aru,\cgRu)$.
\begin{thm}\label{thm:Ru:symms:F}
The group $\Aut(\aru,\cgRu)$ is a direct product $7\times\Ru$.
\end{thm}
\begin{proof}
From Proposition~\ref{prop:Ru:symms:Rumaxfinite} we have that $F$ is
generated by $R\cong \cRu$ together with $\xi\Id$ for $\xi$ some
root of unity. If $\xi\Id\in \SL(\gt{r})$ preserves any non-zero
element of $\bigwedge^{14}(\gt{a})\vac_X\subset (\aru)_{7/2}$ then
$\xi^{14}=1$. Conversely, $\xi\Id$ preserves all elements of
$\cgRu=\{\cas,\vp,\vm,\spsp\}$ when $\xi^{14}=1$. So $F$ is a
central product $14\circ R$ with center generated by $\ogz$. The
image of $F$ in $\SL(\gt{r})/\Kumz$ is evidently $7\times \Ru$, as
required.
\end{proof}

We summarize the results of this section with the following
\begin{thm}\label{thm:Ru:main}
The quadruple $(\aru,Y,\vac,\cgRu)$ is a self-dual enhanced
$\U(1)$--VOA of rank $28$. The full automorphism group of
$(\aru,\cgRu)$ is a direct product of a cyclic group of order seven
with the sporadic simple group of Rudvalis.
\end{thm}

Ultimately we would like to characterize $\aru$ in a fashion
analogous to that applied in \cite{Dun_VACo} to the enhanced VOA
associated to the Conway group. If $U$ is a self-dual nice rational
VOA of rank $28$ with no odd vectors of degree less than $7/2$, then
the character of $U$ is very restricted, and the method of
Proposition~5.7 in \cite{Dun_VACo} can be used to show that the
character must coincide with that of $\aru$. Probably, the
techniques of \S5.1 in \cite{Dun_VACo} can be applied to show that
the VOAs underlying $U$ and $\aru$ are isomorphic. We conjecture
that self-duality and the above vanishing condition, together with
the local algebra (see \S\ref{sec:eVOAs}) determined by $\cgRu$ are
sufficient to determine $\aru$ uniquely (among nice rational
enhanced VOAs).
\begin{conj}
Suppose that $(U,Y,\vac,\cgset)$ is a nice rational enhanced VOA
with $\cfalg(\cgset)\cong\cfalg(\cgRu)$ such that
\begin{enumerate}
\item     $U$ has rank $28$
\item     $U$ is self-dual
\item     $U_{\bar{1}}$ has no non-trivial vectors of degree less than $7/2$
\end{enumerate}
Then $(U,Y,\vac,\cgset)$ is isomorphic to $(\aru,Y,\vac,\cgRu)$.
\end{conj}
This conjectural characterization of $\aru$ (and hence the sporadic
group $\Ru$) is reminiscent of the uniqueness results that exist for
the Golay code, the Leech lattice (see \cite{Con69}), and the
enhanced VOA for Conway's group (see \cite{Dun_VACo}), and those
which are conjectured to hold for the Moonshine VOA (c.f.
\cite{FLM}) and the Baby Monster VO(S)A (c.f. \cite{HohnPhD}). All
these objects have sporadic automorphism groups. The object $\aru$
is a first example with non-Monstrous sporadic automorphism group.

%%%%%%%%%%%%%%%%%%%%%%%%%%%%%%%%%%%%%%%%%%%%%%%%%%%%%%%%%%%%%%%%%%

\section{McKay--Thompson series}\label{sec:series}

In this section we consider the McKay--Thompson series arising from
the enhanced VOAs constructed in \S\S\ref{sec:GLnq},\ref{sec:Ru}.
These series furnish analogues of Monstrous Moonshine for the
sporadic group of Rudvalis, and in this section we will derive
explicit expressions for them.

The main tool for expressing the series explicitly is the notion of
Frame shape, and this is reviewed in \S\ref{sec:series:Frames}. We
then define the McKay--Thompson series associated to the action of a
group on a VOA in \S\ref{sec:series:ord}, and present explicit
expressions for all the McKay--Thompson series arising in our
examples.

For $\U(1)$--VOAs (see \S\ref{sec:eVOAs}) one may consider
McKay--Thompson series in two variables, and we define these in
\S\ref{sec:series:2var}. The enhanced VOA for the Rudvalis group is
a $\U(1)$--VOA, and so this notion applies also to our main example.
As in the ordinary case, one is able to provide explicit expressions
for all the two variable McKay--Thompson series arising for the
Rudvalis group, and we derive these expressions also in
\S\ref{sec:series:2var}. For this we employ the notion of weak Frame
shape, and this is introduced in \S\ref{sec:series:wkFrames}.

Since the Rudvalis group is a non-Monstrous group, the
McKay--Thompson series arising are of particular interest. In
\S\ref{sec:series:MbeyM} (and the Appendix) we provide some further
information about these series.

\subsection{Frame shapes}\label{sec:series:Frames}

Suppose that $\gt{u}$ is a complex vector space with non-degenerate
symmetric bilinear form. Let us set $N$ to be the dimension of
$\gt{u}$, and suppose that $g$ is a finite order element of
$\SO(\gt{u})$ satisfying $g^n=\Id_{\gt{u}}$ for some $m$. Then $g$
has eigenvalues $\{\xi_i\}_{i=1}^N$ say, where each $\xi_i$ is an
$m^{\text{th}}$-root of unity, and the polynomial $\det({\rm
Id}_{\gt{u}}-gx)$ satisfies
\begin{gather}
    \det({\rm Id}_{\gt{u}}-gx)
        =\prod_{i=1}^N(1-\xi_ix)
\end{gather}
Suppose now that the action of $g$ can be written over $\QQ$; that
is, suppose that there is some basis with respect to which $g$ is
represented by a matrix with entries in $\QQ$. (This holds for each
element in $\cRu$, for example, since this group preserves the an
integral lattice --- the Conway--Wales lattice.) Then all primitive
$n^{\text{th}}$-roots of unity in $\{\xi_i\}_{i=1}^{N}$ appear with
the same multiplicity, and there are uniquely determined integers
$m_k\in\ZZ$ for each $k$ dividing $m$ such that we have
\begin{gather}\label{eqn:series:Frames|Frameshape}
    \det({\rm Id}_{\gt{u}}-gx)=\prod_{k|n} (1-x^k)^{m_k}
\end{gather}
Note also that we have $\sum_{k|n}km_k=N$ in this case. For such $g$
in ${\SO}(\gt{u})$ we may express the data $\{(k,m_k)\}$ by writing
\begin{gather}
     g|_{\gt{u}}\sim \prod_{k|n}k^{m_k}
\end{gather}
and the formal expression $\prod_{k|n}k^{m_k}$ is called the {\em
Frame shape} for the action of $g$ on $\gt{u}$. We say that $g$ {\em
admits a Frame shape over $\gt{u}$} if the characteristic polynomial
of $g^{-1}$ has an expression of the form
(\ref{eqn:series:Frames|Frameshape}). %Evidently, $g$ admits a Frame
%shape over $U$ just when the action of $g$ on $U$ can be written
%over $\QQ$.

\subsection{Ordinary McKay--Thompson series}\label{sec:series:ord}

Let $U$ be a VOA of rank $c$ and let $g$ be an (VOA) automorphism of
$U$. Then the {\em McKay--Thompson series associated to the action
of $g$ on $U$} is the $q$-series defined by
\begin{gather}
     \tr_Ugq^{L(0)-c/24}=\sum (\tr_{U_n}g)q^{n-c/24}
\end{gather}
In the special case that $g$ is the identity we recover, what we
call, the {\em character} of $U$.

For most of the examples of enhanced VOAs that arise in this article
the underlying VOA is of the form $\atw{\gt{u}}=A(\gt{u})^0\oplus
A(\gt{u})^0_{\ogi,E}$ (see \S\ref{sec:GLnq:tw}), and the
automorphisms all lie in the corresponding group $\Spql{\gt{u}}$.
Given $g$ in $\Spql{\gt{u}}$ then, we would like to compute the
trace of the operator $gq^{L(0)-c/24}$ on $\atw{\gt{u}}$.

Let $g\in\Spql{\gt{u}}$, and note that $g$ has preimages $\hat{g}$
and $\gt{z}\hat{g}$ say, in $\Sp(\gt{u})$, and these in turn have
two images in $\SO(\gt{u})$; one the negative of the other. We write
$g\leftrightarrow \{\pm \bar{g}\}$ for this correspondence between
elements of $\Spql{\gt{u}}$ and pairs in $\SO(\gt{u})$, and we
consider only $g\in\Spql{\gt{u}}$ for which $\pm \bar{g}$ admit
Frame shapes over $\gt{u}$. (This will hold for all elements in
$\cRu$, thanks to its action on the Conway--Wales lattice.) Then for
$\bar{h}$ either one of $\pm \bar{g}$, suppose that $\bar{h}$ admits
a Frame shape $\prod_{k|n}k^{m_k}$ say, where $n$ is the order of
$\bar{h}$. Recall $\eta(\tau)$, the Dedekind eta function
(\ref{eqn:intro:notation:eta}). Then for such $\bar{h}$ we set
\begin{gather}
    \eta_{\bar{h}}(\tau)=\prod_{k|n}\eta(k\tau)^{m_k}
\end{gather}
For $\hat{h}$ in $\Sp(\gt{u})$ we write $\chi_{\hat{h}}$ for the
trace of $\hat{h}$ on the finite dimensional $\Sp(\gt{u})$-module
$\Cm(\gt{u})_E\hookrightarrow A(\gt{u})_{\ogi}$ (see
\S\ref{sec:cliffalgs:VOAs}). We will also convene that $\hat{g}$ be
the preimage of $g$ in $\Sp(\gt{u})$ whose image in $\SO(\gt{u})$ is
$\bar{g}$, so that $\gt{z}\hat{g}$ is the preimage of $g$ whose
image in $\SO(\gt{u})$ is $-\bar{g}$.

We now have the following
\begin{thm}\label{thm:series:ord|chars}
Let $g\in\Spql{\gt{u}}$, with $\hat{g}\in\Sp(\gt{u})$ and
$\bar{g}\in\SO(\gt{u})$ as above. Then the McKay--Thompson series
associated to the action of $g$ on $\atw{\gt{u}}$ admits the
following expression.
\begin{gather}\label{eqn:series:ord|chars}
    \mathsf{tr}|_{\atw{\gt{u}}}gq^{L(0)-c/24}=
        \frac{1}{2}\left(
          \frac{\eta_{\bar{g}}(\tau/2)}{\eta_{\bar{g}}(\tau)}
          +\frac{\eta_{-\bar{g}}(\tau/2)}{\eta_{-\bar{g}}(\tau)}
          +\chi_{\gt{z}\hat{g}}\eta_{\bar{g}}(\tau)
            +\chi_{\hat{g}}\eta_{-\bar{g}}(\tau)\right)
\end{gather}
\end{thm}
\begin{proof}
Let $g$, $\bar{g}$ and $\hat{g}$ be as in the statement of the
theorem. Then $\bar{g}$ has Frame shape $\prod_{k|n}k^{m_k}$ say,
and the inverse transformation $\bar{g}^{-1}$ has the same Frame
shape (since all primitive $k^{{\rm th}}$-roots of unity appear with
the same multiplicity for any given $k$). Let $\{f_i\}_{i=1}^{24}$
be a basis for $\gt{u}$ consisting of eigenvectors of $\bar{g}$ with
eigenvalues $\{\xi_i\}_{i=1}^{24}$. Then we have
\begin{gather}\label{eqn:series:ord|cpxferevs}
    \det(\Id_{\gt{u}}-\bar{g}x)
    =\prod_i(1-\xi_ix)= \prod_{k|n}(1-x^k)^{m_k}
\end{gather}
Recall that $\atw{\gt{u}}$ may be described as
$\atw{\gt{u}}=A(\gt{u})^0\oplus A(\gt{u})^0_{\ogi}$ (see
\S\ref{sec:cliffalgs:VOAs}, \S\ref{sec:GLnq:tw}). We have
\begin{align}
    {\sf tr}|_{A(\gt{u})}\gt{z}\hat{g}q^{L(0)-c/24}
        =\;&q^{-c/24}\prod_{n\geq 0}\prod_i
            (1-\xi_iq^{n+1/2})\\
     \tr_{A(\gt{u})_{\ogi}}\gt{z}\hat{g}q^{L(0)-c/24}
          =\;&q^{c/12}\chi_{\gt{z}\hat{g}}\prod_{n\geq 1}\prod_i
               (1-\xi_iq^n)
\end{align}
(Recall that the image of $\gt{z}\hat{g}$ in $\SO(\gt{u})$ is
$-\bar{g}$.) Substituting $q^r$ for $x$ in
(\ref{eqn:series:ord|cpxferevs}) we obtain
\begin{gather}
    \prod_i(1-\xi_iq^r)=\prod_{k|m}(1-(q^{k})^r)^{p_k}
\end{gather}
so that for ${\sf tr}|_{A(\gt{l})}\gt{z}\hat{g}q^{L(0)-c/24}$ for
example, we have
\begin{gather}
    \begin{split}
    {\sf tr}|_{A(\gt{l})}\gt{z}\hat{g}q^{L(0)-c/24}
        &=q^{-c/24}\prod_{n\geq 0}\prod_{i}
            (1-\xi_iq^{n+1/2})\\
        &=\prod_{k|n}\left(
            q^{-km_k/48}\prod_{n\geq 0}
            (1-(q^k)^{n+1/2})^{m_k}\right)\\
        &=q^{-c/24}\frac{\eta_{\bar{g}}(\tau/2)}
          {\eta_{\bar{g}}(\tau)}
    \end{split}
\end{gather}
since the rank $c$ is half the dimension of $\gt{u}$, and thus half
the degree of the polynomial $\det(\Id_{\gt{u}}-\bar{g}x)$, so that
$2c=\sum_{k|n}km_k$. A similar argument shows that
\begin{gather}
     {\sf tr}|_{A(\gt{u})_{\ogi}}
     \gt{z}\hat{g}q^{L(0)-c/24}=
     q^{c/12}\chi_{\gt{z}\hat{g}}\eta_{\bar{g}}(\tau)
\end{gather}
and taking $-\bar{g}$ in place of $\bar{g}$, we find
\begin{align}
     \tr_{A(\gt{u})}\hat{g}q^{L(0)-c/24}
          &=q^{-c/24}\frac{\eta_{-\bar{g}}(\tau/2)}
          {\eta_{-\bar{g}}(\tau)}\\
     \tr_{A(\gt{u})_{\ogi}}\hat{g}q^{L(0)-c/24}
          &=q^{c/12}\chi_{\hat{g}}\eta_{-\bar{g}}(\tau)
\end{align}
Finally, for $g\in\Spql{\gt{u}}$ as above, we have
\begin{align}
     \tr_{A(\gt{u})^0}gq^{L(0)-c/24}
     &=\frac{1}{2}\left(
          \tr_{A(\gt{u})}\hat{g}q^{L(0)-c/24}
          +\tr_{A(\gt{u})}\gt{z}\hat{g}q^{L(0)-c/24}
               \right)\\
     \tr_{A(\gt{u})^0_{\ogi}}gq^{L(0)-c/24}
     &=\frac{1}{2}\left(
          \tr_{A(\gt{u})_{\ogi}}\hat{g}q^{L(0)-c/24}
          +\tr_{A(\gt{u})_{\ogi}}\gt{z}\hat{g}
               q^{L(0)-c/24}\right)
\end{align}
since the action of $\gt{z}$ on $A(\gt{u})\oplus A(\gt{u})_{\ogi}$
coincides with that of $\ogi$. This completes the proof.
\end{proof}
One may now study the McKay--Thompson series associated to the
Rudvalis group as soon as the Frame shapes (of elements in the
relevant covering group) are given. These Frame shapes may be
deduced easily from the character tables and power maps of these
groups, and this data is recorded in \cite{ATLAS}. We list the Frame
shapes for elements in the double cover $\cRu$ in Table
\ref{tab:series:Ru|frms} under the columns headed $\SO_{56}$. (See
also the comments in \S\ref{sec:series:2var}.) The character of
$\aru$ will be given in \S\ref{sec:series:2var}.

\subsection{Weak Frame shapes}\label{sec:series:wkFrames}

Now let us suppose that we are in the situation of
\S\ref{sec:cliffalgs:Herm}, so that $\gt{u}$ is of the form
$\gt{u}=\gt{a}\oplus\gt{a}^*$ for some complex vector space $\gt{a}$
with Hermitian form, and the symmetric bilinear form on $\gt{u}$ is
that induced by the natural pairing $\gt{a}\times\gt{a}^*\to \CC$,
so that any invertible linear transformation of $\gt{a}$ extends
naturally to an orthogonal transformation of $\gt{u}$. We consider a
finite order element $g$ in $\SU(\gt{a})$ with $g^n={\rm
Id}_{\gt{a}}$. Then $g$ has eigenvalues $\{\gamma_i\}$ say on
$\gt{a}$, where each $\gamma_i$ is an $n^{\text{th}}$-root of unity
and the polynomial $\det({\rm Id}_{\gt{u}}-gx)$ (corresponding to
the action of $g$ on $\gt{u}$) satisfies
\begin{gather}
    \det({\rm Id}_{\gt{u}}-gx)
        =\prod_{i=1}^N(1-\gamma_ix)(1-\gamma_i^{-1}x)
\end{gather}
For $g$ in $\cRu$ it will always be the case that all primitive
$n^{\text{th}}$-roots of unity in $\{\gamma_i,\gamma_i^{-1}\}$
appear with the same multiplicity (thanks to the action of this
group on an integral lattice --- the Conway--Wales lattice) so that
we may consider the Frame shape for the action of such an element on
$\gt{u}$, as defined in \S\ref{sec:series:Frames}. Regarding $g$ as
a unitary transformation on $\gt{a}$, we are in general not able to
find such a nice expression for $\det({\rm Id}_{\gt{a}}-gx)$ ---
corresponding to the action of $g$ on $\gt{a}$ --- as in
(\ref{eqn:series:Frames|Frameshape}), but we may consider the
following generalization: the notion of weak Frame shape.

For $m$ an integer, $k$ a nonnegative integer, and $a\in\QQ/\ZZ$ we
write $k^{m}_{a}$ to indicate a factor of the form $(1-e^{2\pi\ii
a}x^k)^m$. We write $k^m$ in place of $k^m_0$ so as to incorporate
the notation for ordinary Frame shapes, and we agree to multiply the
symbols $k^m_a$ just as we do the corresponding polynomials, thus
obtaining various identities, including the following.
\begin{gather}
    k_ak^m_a=k^{m+1}_a\\
    k^m_{1/p}k^m_{2/p}\cdots k^m_{(p-1)/p}
        =(kp)^mk^{-m}\\%\frac{(kp)^m}{k^m}\\
    k_{a+1/2}^m=(2k)_{2a}^mk_a^{-m}
\end{gather}
We say that $\prod_j (k_j)^{m_j}_{a_j}$ is a {\em weak Frame shape}
for the action of $g$ on $\gt{a}$ when
\begin{gather}
    \det({\rm Id}_{\gt{a}}-gx)=
        \prod_{j} (1-e^{2\pi\ii a_j}x^{k_j})^{m_j}
\end{gather}
Weak Frame shapes always exist, since choosing $a_i\in\QQ/\ZZ$ such
that $\gamma_i=e^{2\pi\ii a_i}$ for example, we have that $\prod_i
1_{a_i}$ is a weak Frame shape for the action of $g$ on $\gt{a}$.
%Using the obvious multiplication rules we may always find a weak
%Frame shape for $g$ of the form $\prod_{k|n}k^{m_k}_{a_k}$, and in
%this case the values $a_k$ and $m_k$ are uniquely determined.

\subsection{Two variable McKay--Thompson
series}\label{sec:series:2var}

Let $\left(U,Y,\vac,\{\cas,\cej\}\right)$ be a $\U(1)$--VOA of rank
$c$ (see \S\ref{sec:eVOAs}), and suppose that the action of
$J(0)=\cej_{(0)}$ on $U$ is diagonalizable. Set $J'(0)=\ii J(0)$. We
say that $u\in U$ has {\em charge $m$} when $J'(0)u= mu$.
\begin{defn}
For $g$ an automorphism of $(U,\{\cas,\cej\})$, the {\em two
variable McKay--Thompson series associated to the action of $g$ on
$U$} is the series in variables $p$ and $q$ given by
\begin{gather}
     \tr_{U}gp^{J'(0)}q^{L(0)-c/24}
     =\sum_{m,n}(\tr_{U_{n}^m}g)p^mq^{n-c/24}
\end{gather}
where $U_n^m$ denotes the subspace of $U$ of degree $n$ consisting
of vectors of charge $m$.
\end{defn}
In the limit as $p\to 1$ we recover the (ordinary) McKay--Thompson
series associated to the action of $g$ on $U$. In the case that $g$
is the identity we obtain what we call the {\em two variable
character of $U$}.

We would like to compute the two variable McKay--Thompson series
arising from the action of $\Ru$ on $\aru$. Observe that this action
is contained in that of $\SL(\gt{r})/\Kumgp$, and even in that of
$\SUql{\gt{r}}$ where $\gt{r}$ is the Hermitian vector space of
\S\ref{sec:Ru:geom} or \S\ref{sec:Ru:mnml}, and the VOA underlying
$\aru$ is of the form $\atw{\gt{s}}$ for $\gt{s}=\gt{r}\oplus
\gt{r}^*$.

Similar to \S\ref{sec:series:ord}, we suppose that $g$ lies in
$\SUql{\gt{a}}$, and we let $\pm\bar{g}$ be the preimages of $g$ in
$\SU(\gt{a})\subset\SO(\gt{u})$. For $\bar{h}$ one of $\pm\bar{g}$,
one has a weak Frame shape $\prod_j(k_j)_{a_j}^{m_j}$ say, for the
action of $\bar{h}$ on $\gt{a}$, curtesy of
\S\ref{sec:series:wkFrames}. Recall the Jacobi theta function
$\vartheta(z|\tau)$ from (\ref{eqn:intro:notation:theta}). For
$\bar{h}$ as above, we set
\begin{align}\label{eqn:series:2var|phipsidef}
     \phi_{\bar{h}}(z|\tau)
          &=\prod_j
          \frac{\vartheta(k_jz+\tfrac{1}{2}-a_j|k_j\tau)^{m_j}}
               {\eta(k_j\tau)^{m_j}}\\
     \psi_{\bar{h}}(z|\tau)
          &=\prod_j
          p^{k_jm_j/2}q^{k_jm_j/8}
          \frac{\vartheta(k_jz+\tfrac{1}{2}k_j\tau+\tfrac{1}{2}-a_j|k_j\tau)^{m_j}}
               {\eta(k_j\tau)^{m_j}}
\end{align}
and we then have
\begin{thm}\label{thm:series:2var|chars}
Let $g\in\SUql{\gt{a}}$ with $\pm\bar{g}$ as above. Then the two
variable McKay--Thompson series associated to the action of $g$ on
$(\atw{\gt{u}},\{\cas,\cej\})$ admits the following expression.
\begin{gather}\label{eqn:series:2var|chars}
    \mathsf{tr}|_{\atw{\gt{u}}}gp^{J'(0)}q^{L(0)-c/24}=
        \frac{1}{2}\left(
          \phi_{-\bar{g}}(z|\tau)+\phi_{\bar{g}}(z|\tau)
          +\psi_{-\bar{g}}(z|\tau)+\psi_{\bar{g}}(z|\tau)
               \right)
\end{gather}
\end{thm}
\begin{proof}
The verification is very similar to that of Theorem~
\ref{thm:series:ord|chars}. With $g$ and $\pm\bar{g}$ as above, let
$\prod_j(k_j)_{a_j}^{m_j}$ be a weak Frame shape for $\bar{g}$ say,
and let $\{\gamma_i\}$ be the eigenvalues for the action of $g$ on
$\gt{a}$. Then we have
\begin{gather}\label{eqn:series:2var|cpxfer}
     %\begin{split}
     \det(\Id_{\gt{u}}-gx)%&=
     %\prod_i(\gamma_i^{-1}-x)(\gamma_i-x)\\
     %&=
     %\prod_i(1-\gamma_ix)(1-\gamma_i^{-1}x)\\
     %&
     =
     \prod_j(1-e^{2\pi\ii a_j}x^{k_j})^{m_j}
          (1-e^{-2\pi\ii a_j}x^{k_j})^{m_j}
     %\end{split}
\end{gather}
%since $\det(\bar{g}|_{\gt{a}})=\prod_i\gamma_i=\prod_j e^{2\pi\ii
%a_jm_j}=1$.
For the action of $-\bar{g}$ on $A(\gt{u})$ and
$A(\gt{u})_{\ogi}$ we have
\begin{align}
     \begin{split}
     \tr_{A(\gt{u})}(-\bar{g})p^{J'(0)}q^{L(0)-c/24}
     &=q^{-c/24}\prod_{m\geq 0}\prod_i
          (1-\gamma_ip^{-1}q^{m+1/2})\\
          &\qquad\qquad\times
          (1-\gamma_i^{-1}pq^{m+1/2})
     \end{split}\\
     \tr_{A(\gt{u})_{\ogi}}(-\bar{g})p^{J'(0)}q^{L(0)-c/24}
     &=p^{c/2}q^{c/12}\prod_{m\geq 0}\prod_i
          (1-\gamma_ip^{-1}q^m)
          (1-\gamma_i^{-1}pq^{m+1})
\end{align}
since left multiplication by $a_i(-r)$ decreases charge by $1$, and
left multiplication by $a_i^*(-r)$ increases charge by $1$. We
substitute $p^{-1}q^r$ or $pq^r$ for $x$ in
(\ref{eqn:series:2var|cpxfer}) for various $r$, and compare with
(\ref{eqn:series:2var|phipsidef}) and the definition of
$\vartheta(z|\tau)$ (\ref{eqn:intro:notation:theta}), so as to
obtain
\begin{align}
     \begin{split}
     \tr_{A(\gt{u})}(-\bar{g})p^{J'(0)}q^{L(0)-c/24}
     &=q^{-c/24}\prod_{m\geq 0}\prod_j
          \left(1-e^{2\pi\ii a_j}p^{-k_j}q^{k_j(m+1/2)}\right)^{m_j}\\
          &\quad\times
          \left(1-e^{-2\pi\ii a_j}p^{k_j}q^{k_j(m+1/2)}\right)^{m_j}\\
     &=\phi_{\bar{g}}(z|\tau)
     \end{split}
\end{align}
and similarly,
\begin{align}
     \begin{split}
     \tr_{A(\gt{u})_{\ogi}}(-\bar{g})p^{J'(0)}q^{L(0)-c/24}
     &=p^{c/2}q^{c/12}\prod_{m\geq 0}\prod_j
          \left(1-e^{2\pi\ii a_j}p^{-k_j}q^{k_jm}\right)^{m_j}\\
          &\quad\times
          \left(1-e^{-2\pi\ii a_j}p^{k_j}q^{k_j(m+1)}\right)^{m_j}\\
     &=\psi_{\bar{g}}(z|\tau)
     \end{split}
\end{align}
where we have used the fact that $c=\sum_jk_jm_j$. Finally we note
that
\begin{align}
     \begin{split}
     \tr_{A(\gt{u})^0}gp^{J'(0)}q^{L(0)-c/24}
     &=\frac{1}{2}\left(
          \tr_{A(\gt{u})}\bar{g}p^{J'(0)}q^{L(0)-c/24}
          \right.\\
     &\qquad\left.+\tr_{A(\gt{u})}(-\bar{g})p^{J'(0)}q^{L(0)-c/24}
               \right)
     \end{split}\\
     \begin{split}
     \tr_{A(\gt{u})^0_{\ogi}}gp^{J'(0)}q^{L(0)-c/24}
     &=\frac{1}{2}\left(
          \tr_{A(\gt{u})_{\ogi}}\bar{g}p^{J'(0)}q^{L(0)-c/24}
          \right.\\
          &\qquad\left.+\tr_{A(\gt{u})_{\ogi}}(-\bar{g})p^{J'(0)}
               q^{L(0)-c/24}\right)
     \end{split}
\end{align}
and the desired result follows.
\end{proof}
One may now study the two variable McKay--Thompson series associated
to the Rudvalis group as soon as the weak Frame shapes (of elements
in the double cover) are given. Again, these Frame shapes may be
deduced easily from the character tables and power maps of these
groups, which are recorded in \cite{ATLAS}.

Since they may not be well known, we record the Frame shapes
(columns headed $\SO_{56}$) and weak Frame shapes (columns headed
$\SU_{28}$) for the Rudvalis group in Table
\ref{tab:series:Ru|frms}. There are $36$ conjugacy classes in the
group $\Ru$, and $25$ of these conjugacy classes split in two when
lifted to $\cRu$. Suppose that $g\in \Ru$ has preimages $\pm
\bar{g}$ in $\cRu< {\SU}(\gt{a})$, and consider the multiplicative
map $\iota$ on Frame shapes which fixes $k^m_a$ when $k$ is even,
and maps $k^m_a$ to $(2k)_{2a}^mk_a^{-m}$ when $k$ is odd. Then for
$\prod_j(k_j)_{a_j}^{m_j}$ a weak Frame shape for $\bar{g}$, the
image of this under $\iota$ is a weak Frame shape for $-\bar{g}$, so
we list Frame shapes and weak Frame shapes for only one of $\pm
\bar{g}\in \cRu$ as $\bar{g}$ ranges over a set of conjugacy class
representatives in $\Ru$. Our naming of the conjugacy classes
follows \cite{ATLAS}.
\begin{table}
  \centering
  \caption{Frame Shapes for $2.\!\Ru$}
  \label{tab:series:Ru|frms}
%\begin{small}
%\begin{footnotesize}
\renewcommand{\arraystretch}{1.1}
\begin{tabular}{c|cc||c|cc}
  Class & $\SO_{56}$ & ${\SU}_{28}$ & Class & $\SO_{56}$ & ${\SU}_{28}$   \\  \hline
    1A &  $1^{56}$ &  $1^{28}$ &12B & $4^112^5/2^16^1$ & $1_{1/4}^13_{3/4}^112^2$   \\
    2A &  $1^82^{24}$ &    $1^42^{12}$ &13A & $1^413^4$ & $1^213^2$ \\
    2B &  $4^{28}/2^{28}$ &    $4^{14}/2^{14}$ &14A & $28^4/14^4$ & $28^2/14^2$ \\
    3A &  $1^23^{18}$ &  $1^13^9$ &14B & $28^4/14^4$ & $28^2/14^2$  \\
    4A &  $1^84^{12}$ &  $1^44^6$ &14C & $28^4/14^4$ & $28^2/14^2$  \\
    4B &  $4^{16}/2^4$ & $1_{1/4}^44^{6}$ &15A & $1^215^4/3^2$ & $1^115^2/3^1$  \\
    4C &  $4^{16}/2^4$ & $4^8/2^2$ &16A & $1^216^4/2^18^1$ & $1^11_{1/4}^12_{3/4}^116^2/8^1$\\
    4D &  $2^44^{12}$ &  $2^24^6$ &16B & $1^216^4/2^18^1$ & $1^11_{3/4}^12_{1/4}^116^2/8^1$ \\
    5A &  $1^65^{10}$ &  $1^35^5$ &20A & $1^24^210^220^2/2^25^2$ & $1^14^110^120^1/2^15^1$ \\
    5B &  $5^{12}/1^4$ & $5^{6}/1^2$ &20B & $2^14^120^3/10^1$ & $1_{1/4}^12^15_{3/4}^120^1$ \\
    6A &  $1^23^26^8$ & $1^13^16^4$ &20C & $2^14^120^3/10^1$ & $1_{1/4}^12^15_{3/4}^120^1$\\
    7A &  $7^8$ &   $7^4$ &24A & $4^112^124^2/2^16^1$ & $1_{1/4}^13_{1/4}^124^1$ \\
    8A &  $4^48^6/2^4$ & $1_{1/4}^24^18^3/2^1$ &24B & $4^112^124^2/2^16^1$ & $1_{1/4}^13_{1/4}^124^1$ \\
    8B &  $8^8/4^2$ & $1_{1/4}^22^18^4/4^2$ &26A & $4^252^2/2^226^2$ & $4^152^1/2^126^1$\\
    8C &  $4^28^6$ & $4^18^3$ &26B & $4^252^2/2^226^2$ & $4^152^1/2^126^1$ \\
    10A & $2^45^210^4/1^2$ & $2^25^110^2/1^1$ &26C & $4^252^2/2^226^2$ & $4^152^1/2^126^1$\\
    10B & $2^220^6/4^210^6$ & $2^120^3/4^110^3$ &29A & $29^2/1^2$ & $29^1/1^1$ \\
    12A & $1^23^212^4$ & $1^13^112^2$ &29B & $29^2/1^2$ & $29^1/1^1$   \\
\end{tabular}
%\end{small}
%\end{footnotesize}
\end{table}

We record the two variable character of $\aru$ in the following
proposition.
\begin{prop}\label{prop:series:2var|charru}
For the two variable character of $\aru$ we have
\begin{gather}
     \begin{split}
     \tr_{\aru}p^{J'(0)}&q^{L(0)-c/24}=
          \;
          \frac{\vartheta(z|\tau)^{28}}{2\eta(\tau)^{28}}
          +\frac{\vartheta(z+\tfrac{1}{2}|\tau)^{28}}{2\eta(\tau)^{28}}\\
          &\;
          +p^{14}q^{7/2}
          \frac{\vartheta(z+\tfrac{1}{2}\tau|\tau)^{28}}{2\eta(\tau)^{28}}
          +p^{14}q^{7/2}
          \frac{\vartheta(z+\tfrac{1}{2}\tau+\tfrac{1}{2}|\tau)^{28}}{2\eta(\tau)^{28}}
     \end{split}
\end{gather}
\end{prop}

\subsection{Two variable modular invariance}\label{sec:series:modinv}

Recall that the Jacboi theta function $\vartheta(z|\tau)$ is a
holomorphic function on $\CC\times \hh$, with the Fourier
development
\begin{gather}\label{eqn:series:modinv|jacser}
     \vartheta(z|\tau)=\sum_{m\in\ZZ}p^mq^{m^2/2}
\end{gather}
for $p=e^{2\pi\ii z}$ and $q=e^{2\pi\ii\tau}$ (c.f.
(\ref{eqn:intro:notation:theta})). From
(\ref{eqn:series:modinv|jacser}) it follows easily that
$\vartheta(z+1|\tau)=\vartheta(z|\tau)$, and that
$\vartheta(z+\tau|\tau)=p^{-1}q^{-1/2}\vartheta(z|\tau)$.

\begin{defn}
For $m$ an integer, we set $\gt{E}_m$ to be the class of functions
$f(z|\tau)$ holomorphic on $\CC\times\hh$ and such that
$f(z+1|\tau)=f(z|\tau)$, and such that for given $\tau$ we have
$f(z+\tau|\tau)=e^{-2\pi\ii z m}c(\tau)f(z|\tau)$ for some function
$c(\tau)$, holomorphic on $\hh$.
\end{defn}
Evidently, $\vartheta(z|\tau)$ belongs to $\gt{E}_1$. From Theorem~
\ref{thm:series:2var|chars} and its proof, we see that the two
variable character of $(\atw{\gt{u}},\cgU)$ (see
\S\ref{sec:GLnq:SLgps}) lies in $\gt{E}_c$ where $c$ is the rank of
$\atw{\gt{u}}$, and similarly for $(\au,\cgU)$ (see
\S\ref{sec:GLnq:untw}).
\begin{rmk}The class $\gt{E}_m$ is
closely related to the class of Jacobi forms of weight $0$ and index
$m$ (c.f. \cite{EicZagThyJacFrms}).
\end{rmk}

Let us define $\varepsilon(z|\tau)=e^{\pi\ii z^2/\tau}$. Then we
have
\begin{prop}
The following operations generate an action of $\SL_2(\ZZ)$ on
$\gt{E}_m$.
\begin{align}
     T:f(z|\tau)&\mapsto f(z|\tau+1)\\
     S:f(z|\tau)&\mapsto \varepsilon(z|\tau)^{-m} f(z/\tau|-1/\tau)
\end{align}
\end{prop}
The Poisson summation formula implies the following identity for the
Jacobi theta function.
\begin{gather}
     \vartheta(z/\tau|-1/\tau)=
          (-\ii\tau)^{1/2}\varepsilon(z|\tau)\vartheta(z|\tau)
\end{gather}
On the other hand, one has
$\eta(-1/\tau)=(-\ii\tau)^{1/2}\eta(\tau)$ for the Dedekind eta
function, and from these observations the following proposition
quickly follows.
\begin{prop}\label{prop:series:modinv|suzruinv}
The two variable character $\aru$ belongs to $\gt{E}_{28}$, and
spans a one dimensional representation of the subgroup of
$\SL_2(\ZZ)$ generated by $S$ and $T^2$.
\end{prop}
Given the conclusion of Proposition~
\ref{prop:series:modinv|suzruinv}, it is natural to ask if Zhu's
modular theory for VOAs, and the analogous theory for VOAs due to
Dong and Zhao, may be extended to $\U(1)$--VOAs, with Jacobi forms
on $\CC\times\hh$ taking up the role played by modular forms in the
ordinary case.

\subsection{Moonshine beyond the Monster}\label{sec:series:MbeyM}

The terms of lowest charge and degree in the character of $\aru$ are
recorded in Table \ref{tab:Ruseries}. The column headed $m$ is the
coefficient of $p^m$ (as a series in $q$), and the row headed $n$ is
the coefficient of $q^{n-c/24}$ (as a series in $p$). The
coefficients of $p^{-m}$ and $p^m$ coincide, and all subspaces of
odd charge vanish.
\begin{table}
  \centering
  \caption{Character of $\aru$}
  \label{tab:Ruseries}
%  \begin{small}
  \begin{tabular}{c|lllll}
    %\hline
     & 0 & 2 & 4 & 6 & 8\\\hline
    0  & 1 &   &   &  & \\
    1/2&   &   &   &   &\\
    1  & 784 & $378$ &  & & \\
    3/2 &  &  &  & & \\
    2 & 144452 & 92512 & 20475 & & \\
    5/2 &  &  &  & & \\
    3 & 11327232 & 8128792 & 2843568 & 376740 & \\
    7/2 & 40116600 & 30421755 & 13123110 & 3108105 & 376740 \\
    4 &  490068257 & 373673216 & 161446572 & 35904960 & 3108105\\
    9/2 & 2096760960 & 1649657520 & 794670240
     & 226546320 & 35904960\\
     5 & 13668945136 & 10818453324 & 5284484352 & 1513872360 & 226546320\\
     11/2 & 56547022140 & 45624923820 & 23757475560 & 7766243940 & 1513872360 \\
    %\hline
  \end{tabular}
%  \end{small}
\end{table}

Many irreducible representations of $\Ru$ are visible in the entries
of Table \ref{tab:Ruseries}. For example, we have the following
equalities, where the left hand sides are the dimensions of
homogeneous subspaces of $\aru$, and the right hand sides indicate
decompositions into irreducibles for the Rudvalis group.
\begin{gather}
     \begin{split}\label{eqn:series:MbeyM|AnMTobs}
     378&=378\\
     784&=1+783\\
     20475&=20475\\
     92512&=(2)378+406+91350\\
     144452&=(3)1+(3)783+65975+76125\\
     376740&=27405+65975+75400+102400
     \end{split}
\end{gather}
The identities in (\ref{eqn:series:MbeyM|AnMTobs}) may be regarded
as analogues for the Rudvalis group, of the original observations,
such as those of (\ref{eqn:series:MbeyM|MTobs})
\begin{gather}
     \begin{split}\label{eqn:series:MbeyM|MTobs}
     196884&=1+196883\\
     21493760&=1+196883+21296876\\
     864299970&=(2)1+(2)196883+21296876+842609326
     \end{split}
\end{gather}
made by McKay and Thompson, connecting the Monster group with
Klein's modular invariant.

\section*{Acknowledgement}

A number of the ideas in this article appear also in the author's
dissertation \cite{Dun_Phd}. The author thanks his thesis advisor
Igor Frenkel for exceptional guidance and advice throughout its
completion. The author is grateful to Richard Borcherds, Robert
Griess, James Lepowsky, and Arne Meurman, for comments on an early
draft of this article, and thanks Jia-Chen Fu for many helpful
discussions.

% ------------------------------------------------------------------------
\newcommand{\etalchar}[1]{$^{#1}$}

% ------------------------------------------------------------------------
\end{document}